\documentclass[aap,MSNbibl,seceqn,dvips]{arximspdf}

\doi{10.1214/13-AAP960} 
\volume{24}
\issue{4}
\pubyear{2014}
\firstpage{1691}
\lastpage{1738}

\makeatletter

\newcommand{\rrvert}{\vert}
\newcommand{\llvert}{\vert}
\def\cal{\mathcal}
\newcommand{\eqref}[1]{(\ref{#1})}
\newtheorem{theorem}{Theorem}
\newproclaim{definition}[theorem]{Definition}
\newtheorem{lemma}[theorem]{Lemma}
\newtheorem{proposition}[theorem]{Proposition}
\newproclaim{remark}[theorem]{Remark}

\newcommand{\eps}{\varepsilon}
\newcommand{\mA}{{\cal A}}
\newcommand{\mL}{\mathcal L}
\newcommand{\zz}{\mathbf{z}}
\newcommand{\bfone}{\mathbf{1}}

\makeatother

\begin{document}
\begin{frontmatter}

\title{On the conditional distributions and the efficient simulations
of exponential integrals of Gaussian random fields}
\runtitle{Exponential integrals of gaussian random fields}

\begin{aug}
\author[a]{\fnms{Jingchen} \snm{Liu}\thanksref{t1}\ead[label=e1]{jcliu@stat.columbia.edu}}
\and
\author[b]{\fnms{Gongjun} \snm{Xu}\ead[label=e2]{xuxxx360@umn.edu}}
\thankstext{t1}{Supported by Institute of Education Sciences through
Grant R305D100017, NSF CMMI-1069064, NSF SES-1323977 and Army Research Office Grant W911NF-14-1-0020.}
\runauthor{J. Liu and G. Xu}
\affiliation{Columbia University and University of Minnesota}
\address[a]{Department of Statistics\\
Columbia University\\
1255 Amsterdam Avenue\\
New York, New York 10027\\
USA\\
\printead{e1}}
\address[b]{School of Statistics\\
University of Minnesota\\
367 Ford Hall, 224 Church St. SE\\
Minneapolis, Minnesota 55455\\
USA\\
\printead{e2}}
\end{aug}

\received{\smonth{5} \syear{2013}}

%
\begin{abstract}
In this paper, we consider the extreme behavior of a Gaussian random
field $f(t)$ living on a compact set $T$. In particular, we are
interested in tail events associated with the integral $\int_{T}
e^{f(t)}\,dt$. We construct a (non-Gaussian) random field whose
distribution can be explicitly stated. This field approximates the
conditional Gaussian random field $f$ (given that $\int_{T} e^{f(t)}\,dt$
exceeds a large value) in total variation. Based on this approximation,
we show that the tail event of $\int_{T} e^{f(t)}\,dt$ is asymptotically
equivalent to the tail event of $\sup_T \gamma(t)$ where $\gamma(t)$
is a Gaussian process and it is an affine function of $f(t)$ and its
derivative field.
In addition to the asymptotic description of the conditional field, we
construct an efficient Monte Carlo estimator that runs in polynomial
time of $\log b$ to compute the probability $P(\int_{T} e^{f(t)}\,dt >b)$
with a prescribed relative accuracy.
\end{abstract}

%

\end{frontmatter}

\section{Introduction}

Consider a Gaussian random field $\{f(t)\dvtx t\in T\}$ living on a
$d$-dimensional domain $T\subset R^d$ with zero mean and unit variance,
that is, for every finite subset $\{t_1,\ldots,t_n\}\subset T$,
$(f(t_1),\ldots,f(t_n))$ is a mean zero multivariate Gaussian random vector.
Let $\mu(t)$ be a (deterministic) function and $\sigma\in(0,\infty)$
be a scale factor. Define
%
\begin{equation}
\label{IntI} \mathcal{I}(T)\triangleq\int_T
e^{\sigma f(t)+\mu(t)}\,dt.
\end{equation}
In this paper, we develop a precise asymptotic description of the
conditional distribution of $f$ given that $\mathcal{I}(T)$ exceeds a
large value $b$, that is, $P(\cdot|\mathcal{I}(T)>b)$.
In particular, we provide a tractable total variation approximation (in
the sample path space) for such conditional random fields based on a
change of measure technique.
In addition to the asymptotic descriptions, we design efficient Monte
Carlo estimators that run in polynomial time of $\log b$ for computing
the tail probabilities
%
\begin{equation}
\label{Tails} v(b)= P \bigl(\mathcal{I}(T)>b \bigr) =P \biggl(\int
_T e^{\sigma
f(t)+\mu
(t)}\,dt > b \biggr)
\end{equation}
with a prescribed relative accuracy.

\subsection{The literature}
In the probability literature, the extreme behaviors of Gaussian random
fields have been studied extensively. The results range from general bounds
to sharp asymptotic approximations. An incomplete list of works
includes {\cite{Hu90,LS70,MS70,ST74,Bor75,CIS,Berman85,LT91,TA96}}.
A few lines of investigations on the supremum norm are given as
follows. Assuming locally stationary structure, the double-sum
method \cite{Pit96} provides the exact asymptotic approximation of
$\sup_T f(t)$ over a compact set $T$, which is allowed to grow as
the threshold tends to infinity. For almost surely at least twice
differentiable fields, the authors of \cite{Adl81,TTA05,AdlTay07}
derive the
analytic form of the expected Euler--Poincar\'e characteristics of
the excursion set [$\chi(A_b)$] which serves as a good approximation
of the tail probability of the supremum. The tube method
\cite{Sun93} takes advantage of the Karhune--Lo\`eve expansion and
Weyl's formula. A recent related work along this line is given by
\cite{NSY08}. The Rice method \cite{AW05,AW08,AW09} provides an
implicit description of $\sup_T f(t)$.
Change of measure based rare-event simulations are studied in \cite{ABL09}.
The discussions also go
beyond the Gaussian fields. For instance, \cite{HPZ11} discusses the
situations of Gaussian process with random variances. See also \cite
{AST09} for discussions
on non-Gaussian cases.
The distribution of $\mathcal{I}(T)$
is studied in the literature when $f(t)$ is a Brownian motion \cite
{Yor92,Duf01}. Recently, \cite{Liu10,LiuXu11} derive the asymptotic
approximations of $P(\mathcal{I}(T)>b)$ as $b\rightarrow\infty$ for
three times differentiable and homogeneous Gaussian random fields.

Besides the tail probability approximations, rigorous analysis of the
conditional distributions of stochastic processes given the occurrence
of rare events is also an important topic. In the classic large
deviations analysis for light-tailed stochastic systems, the sample
path(s) that admits the highest probability (the most likely sample
path) under the conditional distribution given the occurrence of a rare
event is central to the entire analysis in terms of determining the
appropriate exponential change of measure, developing approximations of
the tail probabilities and designing efficient simulation algorithms;
see, for instance, standard textbook \cite{DupEll97}. For heavy-tailed
systems, the conditional distributions and the most likely paths, which
typically admit the so-called ``one-big-jump'' principle, are also
intensively studied \cite{ASM00,AsmussenKlupp,BL11}. These results
not only provide intuitive and qualitative descriptions of the
conditional distribution, but also shed light on the design of
rare-event simulation algorithms \cite{BL08,BL10,BL11}---the best
importance sampling estimator of the rare-event probability uses a
change of measure corresponding to the interesting conditional
distribution. In addition, the conditional distribution (or the
conditional expectations) is also of practical interest. For instance,
in risk management, the conditional expected loss given some
rare/disastrous event is an important risk measure and stress test.

In the literature of Gaussian random fields, the exact Slepian models
[conditional field given a local maximum or level crossing of $f(t)$]
are studied intensively for twice differentiable fields. For instance,
Leadbetter, Lindgren and Rootz\'en \cite{Leadbetter83} give the Slepian
model conditioning on an upcrossing of level $u$ at time zero. Lindgren
\cite{Lindgren70} treats conditioning on a local maximum of height $u$
at time zero.
The first rigorous treatment of Slepian models for nonstationary
processes is given by Lindgren \cite{Lindgren79}. Grigoriu \cite
{Grigoriu89} extends the results of Leadbetter, Lindgren and Rootz\'en~\cite{Leadbetter83} for level crossings to the general nonstationary case.
This work is followed up by Gadrich and Adler \cite{GadrichAdler93}. In
the later analysis, we will set an asymptotic equivalence between the
conditional distribution given $\{\mathcal I (T)>b\}$ and that given the
high excursion of the supremem of $f$. The latter can be characterized
by the Slepain model.

\subsection{Contributions}
In this paper, we pursue along this line for the extreme behaviors of
Gaussian processes and begin to describe the conditional distribution
of $f$ given the occurrence of the event $\{\mathcal{I}(T)>b\}$.
In particular, we provide both quantitative and qualitative
descriptions of this conditional distribution. Furthermore, from a
computational point of view, we construct a Monte Carlo estimator that
takes a polynomial computational cost (in $\log b$) to estimate $v(b)$
for a prescribed relative accuracy.

Central to the analysis is the construction of a change of measure on
the space $C(T)$ (continuous functions living on $T$).
The application of the change of measure ideas is common in the study
of large deviations analysis for the light-tailed stochastic systems.
However, it is not at all standard in the study of Gaussian random fields.
The proposed change of measure is not of a classical
exponential-tilting form. This measure has several features that are
appealing both theoretically and computationally. First, we show that
the change of measure denoted by $Q$ approximates the conditional
measure $P(\cdot| \mathcal{I}(T)>b)$ in total variation as
$b\rightarrow\infty$. Second, the measure $Q$ is analytically
tractable in the sense that the distribution of $f$ under $Q$ has a
closed form representation and the Radon--Nikodym derivative $dQ/dP$
takes the form of a $d$-dimensional integral.
This tractability property has useful consequences. From a
methodological point of view, the measure $Q$ provides a very precise
description of the mechanism that drives the rare event $\{\mathcal
I(T) >b\}$. This result allows us to directly use the intuitive
mechanism to provide functional probabilistic descriptions that
emphasize the most important elements that are present in the
interesting rare events.
More technically, the analytical computations associated with the
measure $Q$ are easy (compared to the conditional measure), and the
expectation $E^Q[\cdot]$ is theoretically much more tractable than
$E[\cdot|\mathcal{I}(T)>b]$.
Based on this result, we show that the tail event $\{\mathcal I(T)>b\}$
is asymptotically equivalent to the tail event of $\sup_T \gamma(t)$
where $\gamma(t)$ is an affine function of $f(t)$ and its derivative
field $\partial^2 f(t)$ and $\gamma(t)$ implicitly depends on $b$.
Thus, one can further characterize the conditional measure by means of
the results on the Slepian model mentioned earlier.

Another contribution of this paper lies in the numerical evaluation of
$v(b)$. The importance sampling algorithm associated with the proposed
change of measure yields an efficient estimator for computing $v(b)$.
An important issue concerns the implementation of the Monte Carlo
method. The processes considered in this paper are continuous while
computers can only represent discrete objects. Inevitably, we will
introduce a suitable discretization scheme and use discrete (random)
objects to approximate the continuous processes. A naturally raised
issue lies in the control of the approximation error relative to the
probability $v(b)$. We will perform careful analysis and report the
overall computational complexity of the proposed Monte Carlo estimators.

A key requirement of the current analysis is the twice
differentiability of~$f$. Our change of measure is written explicitly
in the form of $f$, $\partial f$ and $\partial^2 f$. A~very interesting
future study would be developing parallel results for
nondifferentiable fields. The technical challenges are two-fold.
First, there is lack of asymptotic analysis for the exponential
integral of general nondifferentiable fields. To the author's best
knowledge, the behavior of $\mathcal I(T)$ for nondifferentiable
processes is investigated only when $f$ is a Brownian motion whose
techniques cannot be extend to general cases \cite{Yor92,Duf01}. In
addition, there is a lack of descriptive tools (such as derivatives and
the Palm model) for nondifferentiable processes. This also leads to
difficulties in describing the Slepian model for level crossing.
To the author's best knowledge, analytic description of Slepian models
for excursion of $\sup_T f(t)$ are available only for twice
differentiable fields. Despite of the smoothness limitation, the
current analysis has important applications the details of which will
be presented in the following section.

The rest of this paper is organized as follows. Two applications of
this work are given in Section~\ref{SecApp}. In Section~\ref{SecMain},
we present the main results including the change of measure, the
approximation of $P(\cdot| \mathcal{I}(T)>b)$ and the efficient Monte
Carlo estimator of $v(b)$. Proofs of the theorems are given in Sections~\ref{SecProof}--\ref{proofth5}.
A supplemental material \cite{supp} is provided
including all the supporting lemmas.

\section{Applications}\label{SecApp}
The integral of exponential functions of Gaussian random fields
plays an important role in many probability models. We present two
such models for which the conditional distribution is of interest and
the underlying random fields are differentiable.

\subsection{Spatial point process}
In spatial point process modeling, let $\lambda(t)$ be the intensity of
a Poisson point process on $T$, denoted by $\{N(A)\dvtx A\subset T\}$. In
order to build in spatial dependence structure and to account for
overdispersion, the log-intensity is typically modeled as a Gaussian
random field, that is, $\log\lambda(t) = f(t)+\mu(t)$ and then
$E[N(A)| \lambda(\cdot)] = \int_{A}e^{f(t)+\mu(t)}\,dt$, where $\mu(t)$
is the mean function, and $f(t)$ is a zero-mean Gaussian process. For
instance, Chan and Ledolter \cite{ChLe95} consider the time series
setting in which $T$ is a
one-dimensional interval, $\mu(t)$ is modeled as the observed covariate
process and $f(t)$ is an autoregressive process; see \cite
{DDW00,Camp94,Zeger88,COX55,COIS80} for more examples in
high-dimensional domains.

For the purpose of illustration, we consider a very concrete case that
the point process $N(\cdot)$ represents the spatial distribution of
asthma cases over a geographical domain $T$. The latent intensity
$\lambda(t)$ [or equivalently $f(t)$] represents the unobserved (and
appropriately transformed) metric of the pollution severity at location~$t$.
The mean function can be written as a linear combination of the
observed covariates that may affect the pollution level, that is, $\mu
(t)= \beta^\top x(t)$ is treated as a deterministic function.
It is well understood that $\lambda(t)$ is a smooth function of the
spatial parameter $t$ at the macro level as the atmosphere mixes well;
see, for example, \cite{Aberg08}. One natural question in epidemiology is the
following: upon observing an unusually high number of asthma cases,
what is their geographical distribution, that is, the conditional
distribution of the point process $N(\cdot)$ given that $N(T)>b$ for
some large $b$?

First of all, Liu and Xu
\cite{LiuXu11} show that $P(N(T)> b) \sim P(\mathcal I(T)> b)$ as
$b\rightarrow\infty$.
Following the same derivations, it is not difficult to establish the
following convergence:
\[
P\bigl( \cdot| N(T)> b\bigr) - P\bigl(\cdot| \mathcal I(T)>b\bigr) \to0 \qquad\mbox{in
total variation as $b\to\infty$. }
\]
The total count $N(T)$ is a Poisson random variable with mean $\mathcal
I(T)$. Intuitively speaking, the tail of the integral is similar to a
lognormal random variable and thus is heavy-tailed. Its overshoot over
level $b$ is $O_p(b/\log b)$. However, a Poisson random variable with
mean $\mathcal I(T) \sim b$ has standard deviation $\sqrt b \ll b/\log b$.
Thus, a~large number of $N(T)$ is mainly caused by a large value of
$\mathcal I(T)$.
The symmetric difference of the two sets $\{N(T)> b\}$ and $\{\mathcal
I(T)>b\}$ vanishes, and the probability law of the entire system
conditional upon observing that $N(T)>b$ is asymptotically the same as
that given $\mathcal I(T) > b$. Therefore, the conditional distribution
of $N(\cdot)$ given $N(T)>b$ is asymptotically another
doubly-stochastic Poisson process whose intensity is $\lambda(t) =
e^{\mu(t) + f(t)}$ where $f(t)$ follows the conditional distribution of
$P(f\in\cdot| \mathcal I(T)>b)$.

Based on the main results presented momentarily, a qualitative
description of the conditional distribution of $N(\cdot)$ is as
follows. Given $N(T)> b$, the overshoot is of order $O_p(b/\log b)$,
that is, $N(T) = b + O_p(b/\log b)$.\vadjust{\goodbreak}
The locations of the points are i.i.d. samples approximately following
a $d$-dimensional multivariate Gaussian distribution with mean $\tau
\in
T$ and variance $\Sigma/\log b$ where $\Sigma$ depends on the spectral
moments of $f$.
The distribution of $\tau$ is uniform over $T$ if $\mu(t)$ is a
constant; if $\mu(t)$ is not constant, $\tau$ has a density $l(t)$
presented in \eqref{lt}.

\subsection{Financial application}

The exponential integral can be considered as a generalization of the
sum of dependent lognormal random variables that has been studied
intensively from different aspects in the applied probability
literature (see \cite{DufPan97,Ahs78,BasSha01,GHS00,Due04,AR08,FR10}).
In portfolio risk analysis, consider a portfolio of $n$ assets
$S_{1},\ldots,S_{n}$. The asset prices are usually modeled as log-normal
random variables. That is, let $X_{i} = \log S_{i}$ and
$(X_{1},\ldots,X_{n})$ follows a multivariate normal distribution.
The total portfolio value $S= \sum_{i=1}^{n} w_i S_{i}$ is the weighted
sum of dependent log-normal random variables.

An important question is the behavior of this sum when the portfolio
size becomes large and the assets are highly correlated.
One may employ a latent space approach used in the literature of social
networks. More specifically, we construct a Gaussian process $\{f(t)\dvtx t\in T\}$ and map each asset $i$ to a latent variable $t_i \in T$, that
is, $\log S_i = f(t_i)$.
Then the log-asset prices fall into a subset of the continuous Gaussian process.
Furthermore, we construct a (deterministic) function $w(t)$ so that
$w(t_{i})= w_{i}$. Then, the unit share value of the portfolio is
$\frac{1} n \sum w_{i} S_{i} = \frac{1} n \sum w(t_{i}) e^{f(t_{i})}$.
See \cite{BLY,LiuXu11} for detailed discussions on the random field
representations of large portfolios.

In the asymptotic regime that $n\rightarrow\infty$ and the
correlations among the asset prices become close to one, the subset $\{
t_{i}\}$ becomes dense in $T$. Ultimately, we obtain the limit
\[
\frac{1} n \sum_{i=1}^{n}
w_{i} S_{i} \rightarrow\int_T w(t)
e^{f(t)} h(t) \,dt,
\]
where $h(t)$ is the limiting spatial distribution of $\{t_{i}\}$ in
$T$. Let $\mu(t)= \log w(t) + \log h(t)$. Then the (limiting) unit
share price is
$\mathcal I(T)=\int_T e^{f(t) + \mu(t)}\,dt$.

The current study provides an asymptotic description of the performance
of each asset given the occurrence of the tail event $\mathcal I(T)>b$.
This is of great importance in the study of the so-called \emph{stress
test} that evaluates the impact of shocks on and the vulnerability of a system.
For instance, consider that another investor holds a different
portfolio that has a substantial overlap with the current one, or it
has exactly the same collection of assets but with different weights.
Thus, this second portfolio corresponds to a different mean function
$\mu'(t)$. The stress test investigates the performance of this second
portfolio on the condition that a rare event has occurred to the first,
that is,
\[
P \biggl(\int_T e^{f(t) + \mu'(t)}\,dt\in\cdot \Big| \int
_T e^{f(t) +
\mu
(t)}\,dt > b \biggr).
\]
To characterize the above distribution, we need a precise description
of the conditional measure $P(f\in\cdot| \int_T e^{f(t) + \mu
(t)}\,dt> b)$.



\section{Main results}\label{SecMain}

\subsection{Problem setting and notation}

Throughout this discussion, we consider a homogeneous Gaussian random
field $\{f(t)\dvtx t\in T\}$
living on a domain $T\subset R^d$. Let the covariance function be
\[
C(t-s) = \operatorname{Cov}\bigl(f(t),f(s)\bigr).
\]
We impose the following assumptions:
\begin{longlist}[(C1)]
\item[(C1)] $f$ is stationary with $Ef(t)=0$ and $Ef^2(t) =1$.
\item[(C2)] $f$ is almost surely at least two times differentiable with
respect to~$t$.
\item[(C3)] $T$ is a $d$-dimensional compact set of $R^d$ with
piecewise smooth boundary.
\item[(C4)] The Hessian matrix of $C(t)$ at the origin is standardized
to be $-I$, where~$I$ is the $d\times d$ identity matrix. In addition,
$C(t)$ has the following expansion when $t$ is close to $0$
%
\begin{equation}
C(t)=1-\tfrac{1}{2}t^{\top}t+C_{4}(t)+R_{C}(t),
\end{equation}
where
$C_{4}(t){=}\frac{1}{24}\sum_{ijkl}\partial_{ijkl}^{4}C(0)t_{i}t_{j}t_{k}t_{l}$
and $R_{C}(t){=}O(|t|^{4+\delta_0})$ for some \mbox{$\delta_0>0$}.
\item[(C5)] For each $t\in R^d$, the function $C(\lambda t)$ is a
nonincreasing function of $\lambda\in R^+$.
\item[(C6)] The mean function $\mu(t)$ falls into either of the two cases:
\begin{longlist}[(a)]
\item[(a)] $\mu(t)\equiv0$;
\item[(b)] the maximum of $\mu(t)$ is unique and is attained in the
interior of $T$ and $\mu(t+\varepsilon)- \mu(t) = \varepsilon^\top
\partial\mu(t) + \varepsilon^\top\Delta\mu(t) \varepsilon+
O(|\varepsilon|^{2+\delta_0})$ as $\varepsilon\to0$.
\end{longlist}
\end{longlist}

We define a set of notation constantly used in the later development
and provide some basic calculations.
Let $P^{*}_{b}$ be the conditional measure given $\{\mathcal{I}(T)> b\}
$, that is,
\[
P^{*}_{b}\bigl(f(\cdot) \in A\bigr) = P\bigl(f(\cdot)\in
A|\mathcal{I}(T)> b\bigr).
\]
Let ``$\partial$'' denote the gradient and ``$\Delta$'' denote the
Hessian matrix with respect to $t$. The notation ``$\partial^2$'' is
used to denote the vector of second derivatives. The difference between
$\partial^2 f(t)$ and $\Delta f(t)$ is that $\Delta f(t)$ is a
$d\times
d$ symmetric matrix whose diagonal and upper triangle consist of
elements of $\partial^2 f(t)$.
Furthermore, let $\partial_j f(t)$ be the partial derivative with
respect to the $j$th element of $t$.
Finally, we define the following vectors:
%
\begin{eqnarray}
\label{moment} \mu_{1}(t) &=&-\bigl(\partial_{1}C(t),\ldots,
\partial_{d}C(t)\bigr),\nonumber
\\
\mu_{2}(t) &=& \bigl(\partial^2_{ii}C(t),i=1,\ldots,d;
\partial^2 _{ij}C(t),
\nonumber
\\[-8pt]
\\[-8pt]
\nonumber
&&\hspace*{6pt}{}i=1,\ldots,d-1,j=i+1,\ldots,d \bigr),
\nonumber
\\
\mu_{02}^\top&=&\mu_{20}=\mu_{2}(0).
\nonumber
\end{eqnarray}

Suppose $0\in T$. It is well known that
$ (f(0),\partial^2 f(0), \partial f(0), f(t) )$ is a
multivariate Gaussian random vector with mean zero and covariance
matrix (cf. Chapter~5.5 of \cite{AdlTay07})
\[
\pmatrix{ 1 & \mu_{20} & 0 & C(t)
\vspace*{2pt}\cr
\mu_{02} & \mu_{22} & 0 & \mu_{2}^{\top}(t)
\vspace*{2pt}\cr
0 & 0 & I & \mu_{1}^{\top}(t)
\vspace*{2pt}\cr
C(t) & \mu_{2}(t) & \mu_{1}(t) & 1},
\]
where the matrix $\mu_{22}$ is a $d(d+1)/2$-dimensional positive
definite matrix and contains the 4th order spectral moments arranged in
an appropriate order according to the order of elements in $\partial^2 f(0)$.
Let $h(x,y,z)$ be the density function of $(f(t),\partial f(t),\partial
^{2}f(t))$ evaluated at $(x,y,z)$. Then, simple calculation yields that
%
\begin{eqnarray}
\label{denhxyz}\qquad &&h(x,y,z)
\nonumber
\\[-8pt]
\\[-8pt]
\nonumber
&&\qquad=\frac{\det(\Gamma)^{-{1}/{2}}}{(2\pi)^{
{(d+1)(d+2)}/{4}}}e^{-({1}/{2}) [{y}^{\top}{y}+{(x-\mu_{20}\mu
_{22}^{-1}{z})^{2}}/{%
(1-\mu_{20}\mu_{22}^{-1}\mu_{02})}+{z}^{\top}\mu_{22}^{-1}{z} ] },
\end{eqnarray}
where $\det(\cdot)$ is the determinant of a matrix and
\[
\label{gamma} \Gamma=\pmatrix{ 1 &
\mu_{20}
\vspace*{2pt}\cr
\mu_{02} & \mu_{22}}.
\]

We define $u$ as a function of $b$ such that
%
\begin{equation}
\label{u} \biggl(\frac{2\pi}{\sigma} \biggr)^{{d}/{2}}u^{-{d}/{2}}e^{\sigma u}=b.
\end{equation}
Note that the above equation generally has two solutions: one is
approximately $\sigma^{-1}\log b $, and the other is close to zero as
$b\rightarrow\infty$. We choose $u$ to be the one close to $\sigma
^{-1}\log b$.
For $\mu(t)$ and $\sigma$ appearing in \eqref{IntI}, we define
%
\begin{equation}
\label{notation} \mu_\sigma(t)=\mu(t)/\sigma,\qquad u_t=u-
\mu_\sigma(t).
\end{equation}
Approximately, $u_t$ is the level that $f(t)$ needs to reach so that
$\mathcal I(T)>b$. Furthermore, we need the following spatially varying set:
%
\begin{eqnarray}
A_{t}&=& \bigl\{f(\cdot)\in C(T)\dvtx \alpha_{t}>u_t-
\eta u_t^{-1} \bigr\} ,\label{At}
\end{eqnarray}
where $\eta>0$ is a tuning parameter that will be eventually sent to
zero as $b\rightarrow\infty$ and $\alpha_t$
is a function of $f(t)$ and its derivative fields taking the form of
%
\begin{eqnarray}
\alpha_{t} &=& f(t)+\frac{|{\partial f(t)}|^{2}}{2u_t}+\frac{\mathbf
1^\top\bar f''_t}{2\sigma u_t}+
\frac{B_{t}}{u_t}.\label{S}
\end{eqnarray}
In the above equation \eqref{S}, $\bar f''_t$ is defined as [with the
notation in \eqref{moment}]
%
\begin{equation}
\label{tra} \bar f''_t=
\partial^2 f(t)-u_t\mu_{02}.
\end{equation}
The term $B_t$ is a deterministic function depending only on $C(t)$,
$\mu(t)$ and~$\sigma$,
%
\begin{eqnarray}
\label{B} B_{t}&=&\frac{{\mathbf{1}^{\top}\partial^{2}\mu_\sigma(t)+d\times
\mu
_\sigma(t)}}{2\sigma} +\frac{1}{8\sigma^{2}}\sum
_{i}\partial_{\mathit{iiii}}^{4}C(0)+ \bigl|
\partial\mu_\sigma(t)\bigr|^2,
\end{eqnarray}
where
$d$ is the dimension of $T$, and $\mathbf{1} =(\mathop
{\underbrace{1,\ldots,1}}\limits_{d},\mathop{\underbrace
{0,\ldots,0}}\limits_{d(d-1)/2})^{\top}$.
Note that $\alpha_t \approx f(t)$. Thus on the set $A_t$, $f(t)\approx
\alpha_t > u_t - O(u^{-1})$.
Together with the fact that $E[{\partial^2 f(t)}|f(t) = u_{t}] =
u_t\mu
_{02}$, $\bar f''_t$ is the standardized second derivative of $f$ given
that $f(t) = u_t$. In Section~\ref{SecChange},
we will show that the event $\{\mathcal{I}(T)>b\}$ is approximately
$\bigcup_{t\in T} A_{t}$.

For notational convenience, we write $a_u = O(b_u)$ if there exists a
constant $c>0$ independent of everything such that $a_u \leq cb_u$
for all $u>1$, and $a_u = o(b_u)$ if $a_u/b_u \rightarrow0$ as
$u\rightarrow\infty$, and the convergence is uniform in other
quantities. We write $a_u = \Theta(b_u)$ if $ a_u =O(b_u)$ and $b_u
= O(a_u)$. In addition, we write $a_u\sim b_u$ if $a_u / b_u
\rightarrow1$ as $u\rightarrow\infty$.

\begin{remark}
Condition C1 assumes unit variance. We treat the standard deviation
$\sigma$ as an additional parameter and consider $\int e^{\mu(t) +
\sigma f(t)}\,dt$. Condition C2 implies that $C(t)$ is at least 4 times
differentiable and the first and third derivatives at the origin are
all zero. Differentiability is a crucial assumption in this analysis.
Condition C3 restricts the results to finite horizon.
Condition C4 assumes the Hessian matrix is standardized to be $-I$,
which is to simplify notation. For any Gaussian process $g(t)$ with
covariance function $C_g(t)$ and $\Delta C_g(0) = -\Sigma$ and $\det
(\Sigma)>0$, identity Hessian matrix can be obtained by an affine
transformation by letting $g(t)= f(\Sigma^{1/2}t)$ and
\[
\int_T e^{\mu(t) + \sigma g(t)}\,dt = \det\bigl(
\Sigma^{-1/2}\bigr)\int_{\{s:
\Sigma^{-1/2}s\in T\}}e^{\mu(\Sigma^{-1/2}s)+\sigma f(s)}\,ds.
\]
%
Condition C5 is imposed for technical reasons so that we are able to
localize the integration. For condition C6, we assume that $\mu(t)$
either is a constant or attains its global maximum at one place. If
$\mu
(t)$ has multiple (finitely many) maxima, the techniques developed in
this paper still apply, but the derivations will be more tedious.
Therefore, we stick to the uni-mode case.
\end{remark}

\begin{remark}
The setting in \eqref{Tails} incorporates the case in which the
integral is with respect to other measures with smooth densities.
Then, if $\nu(dt) = \kappa(t)
\,dt$, we will have that
\[
\int_A e^{\mu(t) +\sigma f(t) }\nu(dt) = \int
_A e^{\mu(t)+ \log
\kappa
(t) +\sigma f(t)} \,dt,
\]
which shows that the density can be absorbed by the mean function.
\end{remark}

\subsection{Approximation of the conditional distribution}
\label{SecChange}

In this subsection, we propose a change of measure $Q$ on the sample
path\vadjust{\goodbreak}
space $C(T)$ that approximates $P^*_b$ in total variation. Let
$P$ be the original measure. The measure $Q$ is defined such that $P$
and $Q$
are mutually absolutely continuous. We define the measure $Q$ under two
different scenarios: $\mu(t)$ is not a constant and $\mu(t)\equiv0$.
Note that the measure $Q$ obviously will depend on $b$. To simplify the
notation, we omit the index $b$ in $Q$ whenever there is no ambiguity.


The measure $Q$ takes a mixture form of three measures, which are
weighted by $(1-\rho_{1}-\rho_{2})$, $\rho_{1}$ and $\rho_{2}$,
respectively (a natural constraint is that $\rho_1$, $\rho_2$ and
$1-\rho_1 - \rho_2\in[0,1]$).
We define $Q$ through the Radon--Nikodym derivative
%
\begin{eqnarray}
\label{LRV} \frac{dQ}{dP}&=&(1-{\rho_{1}}-{\rho_{2}})
\int_{T}l(t)\cdot \operatorname{LR}(t)\,dt+\rho _{1}\int
_{T}l(t)\cdot \operatorname{LR}_{1}(t)\,dt
\nonumber
\\[-8pt]
\\[-8pt]
\nonumber
&&{}+{\rho_{2}}
\int_{T}\frac
{\operatorname{LR}_{2}(t)}{\operatorname{mes}(T)}\,dt,
\end{eqnarray}
where $\rho_1, \rho_2$ will be eventually sent to 0 as $b$ goes to
infinity at the rate $(\log\log b)^{-1}$, $\operatorname{mes}(T)$ is the Lebesgue
measure of $T$ and
%
\begin{eqnarray}
\label{LR} \operatorname{LR}(t)&=&\frac{h_{0,t} (f(t),{\partial f(t)},{\partial^2
f(t)}
)}{h (f(t),{\partial f(t)},{\partial^2 f(t)} )},
\nonumber
\\
\operatorname{LR}_{1}(t) &=& \frac{h_{1,t} (f(t),{\partial f(t)},{\partial^2
f(t)} )}{h (f(t),{\partial f(t)},{\partial^2 f(t)}
)},
\\
\operatorname{LR}_{2}(t)&=& \frac{{1}/{\sqrt{2\pi}}e^{-({1}/{2})(f(t)-u_{t})^{2}}} {
({1}/{\sqrt{2\pi}})e^{-({1}/{2})f(t)^{2}}}.\nonumber
\end{eqnarray}
The density $h (f(t),{\partial f(t)},{\partial^2 f(t)} )$ is
defined in \eqref{denhxyz}, $l(t)$ is a density function on $T$,
$h_{0,t}$ and $h_{1,t}$ are two density functions.
Before presenting the specific forms of $l(t)$, $h_{0,t}$ and
$h_{1,t}$, we would like to provide an intuitive explanation of $dQ/dP$
from a simulation point of view.
One can generate $f(t)$ under the measure $Q$ via the following steps:
\begin{enumerate}[(1)]
\item[(1)] Generate $\imath\sim\operatorname{Bernoulli}(\rho_{2})$.

\item[(2)] If $\imath=1$, then:
\begin{enumerate}[(a)]
\item[(a)] generate $\tau$ uniformly from the index set $T$, that is, $\tau
\sim \operatorname{Unif}(T)$;
\item[(b)] given the realized $\tau$, generate $f(\tau) \sim N(u_{\tau},1)$;
\item[(c)] given $(\tau,f(\tau))$, simulate $\{f(t)\dvtx t\neq\tau\}$ from
the original conditional distribution under $P$.
\end{enumerate}

\item[(3)] If $\imath=0$:

\begin{enumerate}[(a)]
\item[(a)] simulate a random variable $\tau$ following the density function $l(t)$;

\item[(b)] given the realized $\tau$, simulate $f(\tau)=x,\partial f(\tau
)=y,\partial^{2}f(\tau)=z$ from density function
%
\begin{equation}
\label{hall} h_{\mathrm{all}}(x,y,z)=\frac{1-\rho_{1}-{\rho_{2}}}{1-\rho_{2}}h_{0,\tau
}(x,y,z)+
\frac{\rho_{1}}{1-\rho_{2}}h_{1,\tau}(x,y,z);
\end{equation}
\item[(c)] given $ (\tau,f(\tau),\partial f(\tau),\partial
^{2}f(\tau
) )$,
simulate $\{f(t)\dvtx t\neq\tau\}$ from the original conditional distribution
under $P$.
\end{enumerate}
\end{enumerate}

Thus, $\tau$ is a random index at which we twist the distribution of
$f$ and its derivatives.
The likelihood ratio at a specific location $\tau$ is given by
$\operatorname{LR}(\tau
)$, $\operatorname{LR}_1 (\tau)$ or $\operatorname{LR}_2(\tau)$ depending on the mixture component.
The distribution of the rest of the field $\{f(t)\dvtx t\neq\tau\}$ given
$(f(\tau),\partial f(\tau), \partial^2 f(\tau))$ is the same as that
under $P$.
It is not hard to verify that the above simulation procedure is
consistent with the Radon--Nikodym derivative in \eqref{LRV}.

We now provide the specific forms of the functions defining $Q$.
We first consider the situation when $\mu(t)\neq0$.
By condition C6, $\mu(t)$ admits its unique maximum at $t_{\ast}=\arg
\sup_{t\in T}\mu(t)$ in the interior of $T$. Furthermore, the Hessian
matrix $\Delta\mu_\sigma(t_*)$ is negative definite.
The function $l(t)$ is a density on $T$ such that for $t\in T$
%
\begin{equation}
\label{lt} l(t)=\bigl(1+o(1)\bigr)\det\bigl(-\Delta\mu_{\sigma}(t_{\ast})
\bigr)^{1/2} \biggl( \frac{u_{t_{\ast}}}{2\pi} \biggr) ^{d/2}e^{({u_{t_{*}}}/{2})(t-t_{\ast})^{\top}\Delta\mu_{\sigma}(t_{\ast
})(t-t_{\ast})},\hspace*{-35pt}
\end{equation}
which is approximately a Gaussian density centered around $t_{*}$.
As $l(t)$ is defined on a compact set $t$, the $o(1)$ term goes to zero
as $b$ tends to infinity. It is introduced to correct for the integral
of $l(t)$ outside the region $T$ that is exponentially small and does
not affect the current analysis.
The functions $h_{0,t}$ and $h_{1,t}$ are density functions on the
vector space where $(f(t),{\partial f(t)},\partial^2 f(t))$ lives, and
they are defined as follows (we will explain the following complicated
functions momentarily):
\begin{eqnarray*}
&&h_{0,t}\bigl(f(t),{\partial f(t)},\partial^2 f(t)\bigr)
\\
&&\qquad=\mathbb {I}_{A_{t}}\times H_{\lambda} \times u_{t}
\times e^{-\lambda u_{t} (f(t)+({\mathbf1^\top
\bar f''_t}/{(2\sigma u_t)})+{B_{t}}/{u_t}-u_{t} )}\times e^{-
{|{\partial f(t)}|^{2}}/{2}}
\\
&&\qquad\quad{}\times\exp \biggl\{ -\frac{1}{2} \biggl[ \frac{|\mu
_{20}\mu_{22}^{-1}\bar f''_{t}|^{2}}{1-\mu_{20}\mu_{22}^{-1}\mu
_{02}}%
+
\biggl\vert\mu_{22}^{-1/2}\bar f''_{t}-
\frac{\mu
_{22}^{1/2}\mathbf
{1}}{%
2\sigma} \biggr\vert^{2} \biggr] \biggr\},
\\
&&h_{1,t}\bigl(f(t),{\partial f(t)},\partial^2 f(t)\bigr)
\\
&&\qquad=\mathbb{I}_{A_{t}^{c}}\times H_{\lambda_{1}}\times u_{t}
\times e^{\lambda_{1}u_{t} (f(t)+({\mathbf1^\top\bar
f''_t}/{(2\sigma u_t)})+{B_{t}}/{u_t}-u_{t} )}\times e^{-
{|{\partial f(t)}|^{2}}/{2}}
\\
&&\qquad\quad{}\times\exp \biggl\{ -\frac{1}{2} \biggl[ \frac{%
|\mu_{20}\mu_{22}^{-1}\bar f''_{t}|^{2}}{1-\mu_{20}\mu_{22}^{-1}\mu
_{02}}+ \biggl\vert
\mu_{22}^{-1/2}\bar f''_{t}-
\frac{\mu
_{22}^{1/2}\mathbf{%
1}}{2\sigma} \biggr\vert^{2} \biggr] \biggr\},
\end{eqnarray*}
where $\mathbb{I}$ is the indicator function, $A_t= \{f(\cdot)\dvtx f(t)+\frac{|{\partial f(t)}|^{2}}{2u_t}+\frac{\mathbf1^\top\bar
f''_t}{2\sigma u_t}+\frac{B_{t}}{u_t} > u_t - \eta/u_t\}$ is defined as
in \eqref{At}, $\bar f''_t$ is defined as in \eqref{tra}, $\lambda<1$
is positive and it will be sent to 1 as $b$ goes to infinity, $\lambda
_1$ is a fixed positive constant (e.g., $\lambda_1 =1$) and the
normalizing constants are defined as
%
\begin{eqnarray}
\label{HL}
H_{\lambda}&=&\frac{e^{-\lambda\eta}(1-\lambda)^{d/2}\lambda}{(2\pi
)^{%
{d}/{2}}}\nonumber\hspace*{-35pt}\\
&&{}\times \biggl[ \int
_{R^{{d(d+1)}/{2}}}e^{-
({1}/{2} )[
{|\mu_{20}\mu_{22}^{-1}z|^{2}}/{(1-\mu_{20}\mu_{22}^{-1}\mu
_{02})}%
+ \vert\mu_{22}^{-1/2}z-{\mu_{22}^{1/2}\mathbf
{1}}/{(2\sigma)}%
\vert^{2} ] }\,dz \biggr] ^{-1},
\nonumber\hspace*{-35pt}
\\[-6pt]
\\[-6pt]
\nonumber\hspace*{-35pt}
\quad H_{\lambda_{1}}&=&\frac{e^{\lambda_{1}\eta}(1+\lambda
_{1})^{d/2}\lambda
_{1}}{(2\pi)^{{d}/{2}}}\\
&&{}\times \biggl[ \int_{R^{{d(d+1)}/{2}}}e^{-({1%
}/{2}) [ {|\mu_{20}\mu_{22}^{-1}z|^{2}}/{(1-\mu_{20}\mu
_{22}^{-1}\mu_{02})}+ \vert\mu_{22}^{-1/2}z-{\mu_{22}^{1/2}
\mathbf{1}}/{(2\sigma)} \vert^{2} ] }\,dz
\biggr] ^{-1}.
\nonumber\hspace*{-35pt}
\end{eqnarray}
The constants $H_\lambda$ and $H_{\lambda_1}$ ensure that $h_{0,t}$ and
$h_{1,t}$ are properly normalized densities.


\textit{Understanding the measure $Q$}. The measure $Q$
is designed such that the distribution of $f$ under the measure $Q$ is
approximately the conditional distribution of $f$ given $\mathcal{I}(T)
>b $.
The two terms corresponding to the probabilities $\rho_{1}$ and $\rho
_{2}$ are included to ensure the absolute continuity and to control the
tail of the likelihood ratio. Thus, $\rho_{1}$ and $\rho_{2}$ will be
sent to zero eventually.

We now provide an explanation of the leading term corresponding to the
probability $1-\rho_{1}-\rho_{2}$.
To understand $h_{0,t}$, we use the notation $\alpha_t$ in \eqref{S}
and rewrite the density function as
\begin{eqnarray*}
&&h_{0,t} \bigl(f(t),{\partial f(t)},\partial^2 f(t) \bigr)
\\
&&\qquad\propto \mathbb{I}_{A_{t}}\exp \bigl\{ -\lambda u_{t}(
\alpha _{t}-u_{t}) \bigr\} \times\exp \biggl\{-
\frac{1-\lambda}{2}\bigl|{\partial f(t)}\bigr|^{2} \biggr\}
\\
&&\qquad\quad{}\times\exp \biggl\{-\frac{1}{2} \biggl[ \frac{|\mu_{20}\mu
_{22}^{-1}\bar f''_{t}|^{2}}{1-\mu_{20}\mu_{22}^{-1}\mu_{02}} + \biggl\vert
\mu_{22}^{-1/2}\bar f''_{t}-
\frac{\mu
_{22}^{1/2}\mathbf
{1}}{%
2\sigma}\biggr \vert^{2} \biggr] \biggr\},
\end{eqnarray*}
%
which factorizes into three pieces consisting of $\alpha_t$, $\partial
f(t)$ and $\bar f''_t$, respectively.
We consider the change of variables from $(f(t), {\partial f(t)},
\partial^2 f(t))$ to $(\alpha_{t},\break  {\partial f(t)},\bar f''_t)$. Then,
under the distribution $h_{0,t}$, the random vectors $\alpha_{t}$,
${\partial f(t)}$ and $\bar f''_t$ are independent. Note that $h_{0,t}$
is defined on the set $A_{t} =\{\alpha_{t} > u_{t} - \eta u_{t}^{-1}\}$
where $\eta$ will be send to zero eventually. Then, $\alpha_{t}-u_t$ is
approximately an exponential random variable with rate $\lambda
{u_{t}}$; ${\partial f(t)}$, and $\bar f''_t$ are two independent
Gaussian random vectors. The density $h_{1,t}$ has a similar interpretation.
The only difference is that $h_{1,t}$ is defined on the set $\{\alpha_t
- u_t < - \eta u_t^{-1}\}$ and $u_t - \alpha_t$ follows approximately
an exponential distribution.
For the last piece corresponding to $\rho_{2}$, the density is simply
an exponential tilting of $f(t)$.

Under the dominating mixture component, to generate an $f(t)$ from $Q$,
a random index $\tau$ is first
sampled from $T$ following density $l(t)$, then $(f(\tau),\partial
f(\tau), \partial^{2} f(\tau))$ is sampled according to $h_{0,\tau}$.
This implies that the large value
of the integral $\int_T e^{\mu(t) + \sigma f(t)}\,dt$ is mostly caused
by the
fact that the field reaches a high level at $\tau$; more precisely,
$\alpha_{\tau}$ reaches a high level of~$u_\tau$ (with an exponential\vadjust{\goodbreak}
overshoot of rate $\lambda u_\tau$). 
Therefore, the random index $\tau$
localizes the position where the field $\alpha_t$ goes very high. The
distribution of $\tau$ given as in \eqref{lt} is very concentrated
around $t_*$. This suggests that the maximum of $\alpha_t$ [or $f(t)$]
is attained within $O_p(u^{-1/2})$ distance from~$t_*$.

We now consider the case where $\mu(t) \equiv0$.
We choose $l(t)$ to be the uniform distribution over set $T$ and have that
%
\begin{eqnarray}
\label{LRC} \frac{dQ}{dP}&=&(1-{\rho_{1}}-{\rho_{2}})
\int_{T}\frac{
\operatorname{LR}(t)}{\operatorname{mes}(T)} \,dt +{\rho_{1}}\int
_{T}\frac{ \operatorname{LR}_{1}(t)}{\operatorname{mes}(T)} \,dt
\nonumber
\\[-8pt]
\\[-8pt]
\nonumber
&&{}+{\rho_{2}}\int
_{T}\frac{ \operatorname{LR}_{2}(t)}{\operatorname{mes}(T)} \,dt,
\end{eqnarray}
where $\operatorname{mes}(\cdot)$ is the Lebesgue measure.
The following theorem states that $Q$ is a good approximation of
$P^*_b$ with appropriate choice of the tuning parameters.\vspace*{-2pt}

\begin{theorem}\label{main}
Consider a Gaussian random field $\{f(t)\dvtx t\in T\}$ living on a
domain $T$ satisfying conditions \textup{C1--C6}.
If we choose the parameters defining the change of measure $\eta= \rho
_1 = \rho_2 = 1-\lambda= (\log\log b)^{-1}$,
then we have the following approximation:
\[
\lim_{b\to\infty}\sup_{A\in\mathcal F} \bigl|Q(A) -
P^{*}_{b}(A)\bigr|=0,
\]
where $\mathcal F$ is the $\sigma$-field where the measures are defined.\vspace*{-2pt}
\end{theorem}

\begin{remark}
Theorem \ref{main} is the central result of this paper. We present its
detailed proof.
The technical developments of other theorems are all based on that of
Theorem \ref{main}. Therefore, we only layout their key steps and the
major differences from that of Theorem \ref{main}.\vspace*{-2pt}
\end{remark}

\begin{remark}
The measure $Q$ in the limit of the above theorem obviously depends on
the tuning parameters ($\eta$, $\rho_1$, $\rho_2$, and $\lambda$) and
the level $b$. To simplify the notation, we omit the indices of those
parameters when there is no ambiguity.\vspace*{-2pt}
\end{remark}

\begin{remark}
The measure corresponding to the last mixture component in~\eqref{LRV},
$\int_{T}\frac{\operatorname{LR}_{2}(t)}{\operatorname{mes}(T)}\,dt$, has been employed by \cite
{LiuXu11} to develop approximations for $v(b)$. We emphasize that the
measure constructed in this paper is substantially different. In fact,
the measure corresponding to $\operatorname{LR}_2(t)$ does not appear in the main
proof. We included it to control the tail of the likelihood ratio in
one lemma.\vspace*{-2pt}
\end{remark}

To illustrate the application of the measure $Q$,
we provide a further characterization of the conditional distribution
$P^*_b$ by presenting another approximation result which is easier to
understand at an intuitive level.
Let
%
\begin{eqnarray}
\label{gammau} \gamma_u(t)& =& f(t)+\frac{\mathbf1^\top\bar f''_t}{2\sigma
u_t}+
\frac
{B_{t}}{u_t}+\mu_\sigma(t),\qquad \beta_u(T)=\sup
_{t\in T}\gamma_u(t),
\nonumber
\\[-8pt]
\\[-8pt]
\nonumber
\qquad{\tilde P}_{b}\bigl(f(\cdot) \in A\bigr) &=& P \bigl(f(\cdot)\in A|
\beta_u(T)> u \bigr).
\end{eqnarray}
%
The process $\gamma_u(t)$ is slightly different than $\alpha_t$.
The following theorem states that the measure $Q$ also approximates the
distribution ${\tilde P}_{b}$ in total variation for $b$ large.
%
\begin{theorem}\label{sup}
Consider a Gaussian random field $\{f(t)\dvtx t\in T\}$ living on a
domain $T$ satisfying conditions \textup{C1--C6}. With the same choice of tuning
parameters as in Theorem \ref{main}, that is, $\eta= \rho_1 = \rho
_2 =
1-\lambda= (\log\log b)^{-1}$, $Q$ approximates ${\tilde P}_{b}$ in
total variation, that is,
\[
\lim_{b\to\infty}\sup_{A\in\mathcal F} \bigl|Q(A) - {\tilde
P}_{b}(A)\bigr|=0.
\]
\end{theorem}

\subsection{Some implications of the theorems}

The results of Theorems \ref{main} and \ref{sup} provide both
qualitative and quantitative descriptions of $P^*_b$.
From a qualitative point of view, Theorems \ref{main} and \ref{sup}
suggest that
%
\begin{equation}
\label{app} \sup_{A\in\mathcal F} \bigl|P^*_b(A) - {\tilde
P}_{b}(A)\bigr|\rightarrow0
\end{equation}
as $b\rightarrow\infty$. Note that $\gamma_u(t)$ itself is a Gaussian
process. Thus, the above convergence result connects the tail events of
exponential integrals to those of the supremum of another Gaussian
random field that is a linear combination of $f$ and its derivative field.
We set up this connection mainly because the distribution of Gaussian
random fields conditional on level crossing (also known as the Slepian
model) is very well studied for smooth processes \cite{GadrichAdler93}. For the purpose of illustration, we cite one result
in Chapter~6.2 of \cite{ATW09} when $\gamma_u(t)$ is stationary and
twice differentiable. Let covariance function of $\gamma_u(t)$ be
$C_\gamma(t)$. Conditional on $\gamma_u(t)$ achieving a local maximum
at location $t^*$ at level $x$, we have the following closed form
representation of the conditional field:
%
\begin{equation}
\label{ssl} \gamma_u\bigl(t^*+t\bigr) = x C_\gamma(t) -
W_x \beta(t) + g(t),
\end{equation}
where
\[
\beta(t) = \pmatrix{ 1&\mu^\gamma_{20}
\vspace*{2pt}\cr
\mu^\gamma_{02}&\mu^\gamma_{22}}^{-1}\mu^{\gamma\top}_2
(t),
\]
$\mu^\gamma_{ij}$'s are the spectral moments of $C_\gamma(t)$,
$W_x$ is a $d(d+1)/2$ dimensional random vector whose density can be
explicitly written down and $g(t)$ is a mean zero Gaussian process
whose covariance function is also in a closed form; see \cite{ATW09}
for the specific forms. If we set $x >u \to\infty$, the local maximum
is asymptotically the global maximum.
Furthermore, thanks to stationarity, the distribution of $t^*$ is
asymptotically uniform over $T$.
The overshoot $x- u$ is asymptotically an exponential random variable.
Thus, the conditional field $\gamma_u(t)$ can be written down
explicitly through representation \eqref{ssl}, the overshoot
distribution and the distribution of $t^*$. Furthermore, the
conditional distribution of $f(t)$ can be implied by~\eqref{gammau} and
conditional normal calculations.


From a quantitative point of view, Theorem \ref{main} implies that for
any bounded function $\Xi\dvtx  C(T) \rightarrow R$ the conditional
expectation $E[\Xi(f)|\mathcal I(T)>b]$ can be approximated by
$E^Q[\Xi
(f)]$, more precisely,
%
\begin{equation}
\label{appexp} E\bigl[\Xi(f)|\mathcal I(T)>b\bigr] - E^Q\bigl[\Xi(f)
\bigr]\rightarrow0
\end{equation}
as $b\rightarrow\infty$.
The expectation $E^Q[\Xi(f)]$ is much easier to compute (both
analytically and numerically) via the following identity:
%
\begin{equation}
\label{exp} E^Q\bigl[\Xi(f)\bigr] = E^Q \bigl[ E\bigl[
\Xi(f) | \imath, \tau, f(\tau), \partial f(\tau), \partial^2 f(\tau)
\bigr] \bigr].
\end{equation}
Note that the inner expectation is under the measure $P$ in that the
conditional distribution of $f$ given $(f(\tau),\partial f(\tau),
\partial^2 f(\tau))$ under $Q$ is the same as that under~$P$.
Furthermore, conditional on $(f(\tau),\partial f(\tau), \partial^2
f(\tau))$, the process $f(t)$ is also a Gaussian process and has the expansion
\[
f(t) = f(\tau)+ \partial f(\tau)^\top(t-\tau) + \tfrac{1} 2
(t-\tau )^\top\Delta f(\tau) (t-\tau) +o\bigl(|t-\tau|^2
\bigr).
\]
These results provide sufficient tools to evaluate the conditional
expectation
\[
E \bigl[\Xi(f) |\imath, \tau, f(\tau), \partial f(\tau), \partial^2
f(\tau) \bigr].
\]
Once the above expectation has been evaluated, we may proceed to the
outer expectation in \eqref{exp}. Note that the inner expectation is a
function of $(\imath, \tau, f(\tau),  \partial f(\tau), \partial^2
f(\tau
))$, the joint distribution of which is in a closed form. Thus,
evaluating the outer expectation is usually an easier task.
In fact, the proof of Theorem \ref{main} is an exercise of the above
strategy by considering that $\Xi(f) = (dP/dQ)^2$.

\begin{remark}
According to the detailed proof of Theorem \ref{main}, the
approximation~\eqref{appexp} is applicable to all the functions such
that $\sup_b E[\Xi^2(f)|\mathcal I(T)>b] < \infty$. To see that, we
need to change the statement and the proof of Lemma \ref{LemTV1}
presented in Section~\ref{SecProof}.
\end{remark}

\subsection{Efficient rare-event simulation for $\mathcal{I}(T)$}
In the preceding subsection we constructed a change of measure that
asymptotically approximates the conditional distribution of $f$ given
$\mathcal{I}(T)>b$.
In this section, we construct an efficient importance sampling
estimator based on this change of measure to compute $v(b)$ as
$b\rightarrow\infty$. We evaluate the overall computation efficiency
using a concept that has its root in the general theory of computation
in both continuous and discrete settings~\cite{MitzUpf05,TraubWW88}.
In particular, completely analogous notions in the setting of
complexity theory of
continuous problems lead to the notion of tractability of a
computational problem \cite{Wos96}.

\begin{definition}
A Monte Carlo estimator is said to be a fully polynomial randomized
approximation scheme (FPRAS) for estimating $v(b)$ if, for some
$q_{1},q_{2}$ and $d>0$, it outputs an averaged estimator that is
guaranteed to have at most $\varepsilon>0$ relative error with
confidence at least $1-\delta\in(0,1)$ in $O(\varepsilon
^{-q_{1}}\delta
^{-q_{2}}|\log v(b)|^{d})$ function evaluations.
\end{definition}
Equivalently, one needs to compute an estimator $Z_{b}$ with complexity
$O(\varepsilon^{-q_{1}}\delta^{-q_{2}}$ $|\log v(b)|^{d})$ such that
%
\begin{equation}
\label{mse} P \bigl(\bigl\llvert Z_{b}/ {v(b)}-1\bigr\rrvert >
\varepsilon \bigr)< \delta.
\end{equation}
In the literature of rare-event simulations, an estimator $L_{b}$ is
said to be \emph{strongly efficient} in estimating $v(b)$ if $EL_{b} =
v(b)$ and $\sup_{b} \operatorname{Var} L_{b}/v^{2}(b)< \infty$. Suppose that a
strongly efficient estimator $L_{b}$ has been obtained. Let $\{
L^{(j)}_{b}\dvtx j=1,\ldots,n\}$ be i.i.d. copies of $L_{b}$. The averaged estimator
\[
Z_{b} = \frac{1} n \sum_{j=1}^{n}L_{b}^{(j)}
\]
has a \emph{relative} mean squared error equal to $\sqrt{E(Z_b/v(b) -
1)^2}=\sqrt{\operatorname{Var}(L_{b})}\times n^{-1/2}v(b)^{-1}$. A simple consequence of
Chebyshev's inequlity yields
\[
P\bigl(\bigl|Z_{b}/v(b) -1 \bigr|\geq\varepsilon\bigr)\leq\frac{\operatorname{Var}(L_{b})}{\varepsilon
^{2}n v^{2}(b)}.
\]
Thus, it suffices to simulate
$n= O(\varepsilon^{-2}\delta^{-1})$
i.i.d. replicates of $L_{b}$ to achieve the accuracy in \eqref{mse}.

The so-called importance sampling is based on the identity $ P(A)
=\break
E^Q[\mathbb I_A \,dP/ dQ]$. The random variable $\mathbb I_A \,dP/dQ$ is an
unbiased estimator of $P(A)$. It is well known that if one chooses
$Q(\cdot)= P(\cdot| A)$, then $\mathbb I_A \,dP/dQ$ has zero variance.
The measure $Q$ created in the previous subsection is a good
approximation of $P^*_b$, and thus it naturally leads an estimator for
$v(b)$ with small variance.

In addition to the variance control, another issue is that the random
fields considered in this paper are continuous objects. A computer can
only perform discrete simulations. Thus we must use a discrete object
approximating the continuous field to implement the algorithms. The
bias caused by the discretization must be well controlled relative to
$v(b)$. In addition, the complexity of generating one such discrete
object should also be considered in order to control the overall
computational complexity to achieve an FPRAS.

We create a regular lattice covering $T$. Define
\[
G_{N,d} = \biggl\{ \biggl(\frac{i_{1}}{N},\frac{i_{2}}{N},\ldots,
\frac
{i_{d}}{N} \biggr)\dvtx i_{1},\ldots,i_{d}\in\mathbb Z
\biggr\}.
\]
For each $t = (t^1, \ldots, t^d)\in G_{N,d}$, define
\[
T_N(t) = \bigl\{ \bigl(s^{1},\ldots, s^{d}
\bigr)\in T\dvtx s^{j}\in\bigl( t^{j}-1/N,t^{j}\bigr]
\hbox{ for }j=1,\ldots,d \bigr\}\vadjust{\goodbreak}
\]
that is, the $\frac{1} N$-cube intersected with $T$ and cornered at $t$.
Furthermore, let
%
\begin{equation}
\label{TM} T_{N} = \bigl\{t\in G_{N,d}\dvtx
T_N(t) \neq\varnothing\bigr\}.
\end{equation}
Since $T$ is compact, $T_{N}$ is a finite set. We enumerate the
elements in $T_{N} = \{t_{1},\ldots, t_M\}$, where $M=O(N^{d})$.
We further define
\[
X=(X_1,\ldots,X_M)^\top\triangleq
\bigl(f(t_1),\ldots,f(t_M)\bigr)^\top
\]
and use
\[
v_M(b)=P \bigl(\mathcal{I}_{M}(T)>b \bigr)
\]
as an approximation of $v(b)$ where
%
\begin{equation}
\label{defIM} \mathcal{I}_{M}(T) = \sum_{i=1}^M
\operatorname{mes}\bigl(T_N(t_i)\bigr) \times e^{\sigma
X_i+\mu
(t_i)}.
\end{equation}
We have the following theorem to control the bias.
%
\begin{theorem}\label{bias}
Consider a Gaussian random field $f$ satisfying conditions in Theorem
\ref{main}. For any $\eps_0>0$, there exists $\kappa_0$ such that for
any $\varepsilon\in(0,1)$, if $N\geq\kappa_0 \eps^{-1-\eps
_0}(\log
b) ^{2+\varepsilon_{0}}$, then for $b>2$
\[
\frac{\llvert  v_M(b)- v(b)\rrvert } {
v(b)}<\varepsilon.
\]
\end{theorem}

We estimate $v_M(b)$ using a discrete version of the change of measure
proposed in the previous section. The specific algorithm is given as follows:

\begin{enumerate}[(1)]
\item[(1)] Generate a random indicator $\imath\sim\operatorname{Bernoulli}(\rho
_{2})$. If $\imath=1$, then:
\begin{enumerate}[(a)]
\item[(a)] generate $\iota$ uniformly from $\{1,\ldots,M\}$;
\item[(b)] generate $X_{\iota} \sim N(u_{t_{\iota}},1)$;
\item[(c)] given $(t_{\iota},X_{\iota})$, simulate the joint field $(f(t),
\partial f(t), \partial^2 f(t))$ on the lattice $T_N\setminus\{
t_{\iota}\}$ from the original conditional distribution under $P$.
\end{enumerate}
\item[(2)] If $\imath=0$:

\begin{enumerate}[(a)]
\item[(a)] if $\mu(t)$ is not constant, simulate a random index $\iota$
proportional to $l(t_{\iota})$, that is, $P(\iota= i) =
l(t_{i})/\kappa
$ and $\kappa= \sum_{i=1}^{M} l(t_{i})$; if $\mu(t)\equiv0$, then
$\iota$ is simulated uniformly over $\{1,\ldots,M\}$;

\item[(b)] given the realized $\iota$, simulate $f(t_{\iota} )=X_{\iota
}=x,\partial f(t_{\iota})=y,\partial^{2}f(t_{\iota} )=z$ from
density function
\[
h_{\mathrm{all}}(x,y,z)=\frac{1-\rho_{1}-{\rho_{2}}}{1-\rho_2}h_{0,t_{\iota
}}(x,y,z)+
\frac{\rho_{1}}{{1-\rho_2}}h_{1,t_{\iota}}(x,y,z);
\]

\item[(c)] given $(t_{\iota},f(t_{\iota} ),\partial f(t_{\iota}
),\partial
^{2}f(t_{\iota} ))$, simulate the joint field $(f(t),\break  \partial f(t),
\partial^2 f(t))$ on the lattice $T_N\setminus\{t_{\iota}\}$ from the
original conditional distribution under $P$.
\end{enumerate}

\item[(3)] Output
%
\begin{eqnarray}
\label{est1} \tilde L_b&=&\mathbb{I}_{\{\mathcal{I}_{M}(T)>b\}}\Big/\Biggl(\frac
{1-{\rho
_1}-{\rho_2}}{\kappa} \sum_{i=1}^M l(t_i) \operatorname{LR}(t_i)+\frac{\rho
_1}{\kappa}
\sum_{i=1}^M l(t_i) \operatorname{LR}_1(t_i)
\nonumber
\\[-8pt]
\\[-8pt]
\nonumber
&&\hspace*{200pt}{}+ {\rho_2} \sum_{i=1}^M \frac{ \operatorname{LR}_2(t_i)}{M}\Biggr).
\end{eqnarray}
\end{enumerate}

Let $Q_M$ be the measure induced by the above simulation scheme. Then
it is not hard to verify that $\tilde L_b = \mathbb{I}_{\{\mathcal
{I}_{M}(T)>b\}}\,dP/dQ_M$, and thus $\tilde L_b$ is an unbiased estimator
of $v_M(b)$. The next theorem states the strong efficiency of the above
algorithm.

\begin{theorem}\label{variance}
Suppose $f$ is a Gaussian random field satisfying conditions in Theorem
\ref{main}. If $N$ is chosen as in Theorem \ref{bias} and all the other
parameters are chosen as in Theorem \ref{main}, then there exists some
constant $\kappa_1 >0$ such that
\[
\sup_{b>1} \frac{E^{Q_M}\tilde L_b^2} {
v^{2}_M(b)} \leq\kappa_1.
\]
\end{theorem}

Let $Z_b$ be the average of $n$ i.i.d. copies of $\tilde L_b$.
According to the results in Theorem \ref{bias}, we have that
\begin{eqnarray*}
\biggl\llvert \frac{Z _b}{v(b)} - 1\biggr\rrvert &\leq&\biggl\llvert
\frac{Z
_b}{v_M(b)}\bigl(v_M(b)/v(b) -1\bigr)\biggr\rrvert +\biggl
\llvert \frac{Z _b}{v_M(b)}-1\biggr\rrvert\\
& \leq&\eps\biggl\llvert
\frac{Z _b}{v_M(b)}\biggr\rrvert + \biggl\llvert
\frac{Z_b
}{v_M(b)}-1\biggr\rrvert .
\end{eqnarray*}
The results of Theorem \ref{variance} indicate that
\[
P\bigl(\bigl|Z_{b}/v_M(b) -1 \bigr|\geq\varepsilon\bigr)\leq
\frac{\kappa
_1}{\varepsilon
^{2}n }.
\]
If we choose $n = \kappa_1 \eps^{-2} \delta^{-1}$, then
\[
P\bigl(\bigl|Z_b/ v(b) - 1 \bigr|\geq3\eps\bigr)\leq\delta.
\]
Thus, the accuracy level as in \eqref{mse} has been achieved. Note that
simulating one $\tilde L_{b}$ consists of generating a multivariate
Gaussian random vector of dimension $M\times(d+1)(d+2)/2=O(N^d) =
O
((\log b)^{(2+\varepsilon_{0})d}\eps^{-(1+\eps_0)d} )$. The
complexity of generating such a vector is at the most $O(N^3)$. Thus
the overall complexity is $O (\varepsilon^{-2-(3+3\varepsilon
_0)d}\delta^{-1}(\log b)^{(6+ 3\varepsilon_{0})d} )$.
The proposed estimator in \eqref{est1} is a FPRAS.

\begin{remark}
The proposed algorithm can also be used to compute conditional
expectations via the representation $E[\Xi(f) | \mathcal I(T)>b] =
E[\Xi
(f);\break \mathcal I(T)>b]/v(b)$, where $E[\Xi(f);  \mathcal I(T)>b]$ can
be estimated by $\Xi(f)\,dP/dQ_M$ and $v(b)$ can be estimated by
$\mathbb
I_{\{\mathcal I(T)>b\}} \,dP/dQ_M$.\vadjust{\goodbreak}
\end{remark}

\section{Proof of Theorem \texorpdfstring{\protect\ref{main}}{3}}\label{SecProof}

We use the following simple yet powerful lemma to prove Theorem \ref{main}.

\begin{lemma}
\label{LemTV1}Let $Q_{0}$ and $Q_{1}$ be probability measures defined
on the
same $\sigma$-field $\mathcal{F}$ such that $dQ_{1}=r^{-1}\,dQ_{0}$ for a
positive random variable $r$. Suppose that for some $\varepsilon>0$,
$E^{Q_{1}} [r^{2} ]=E^{Q_{0}}[r]\leq1+\varepsilon$. Then
\[
\sup_{|X|\leq1}\bigl\llvert E ^{Q_{1}}(X) -
E^{Q_{0}}(X) \bigr\rrvert \leq\varepsilon^{1/2}.
\]
\end{lemma}

\begin{pf}
\begin{eqnarray*}
\bigl\llvert E ^{Q_{1}}(X) - E^{Q_{0}}(X) \bigr\rrvert & =&\bigl
\llvert E^{Q_{1}} \bigl[ (1-r) X \bigr] \bigr\rrvert
\\
& \leq& E^{Q_{1}} \llvert r-1\rrvert \leq\bigl[E^{Q_{1}}(r-1)^{2}
\bigr]^{1/2} \\
&= &\bigl( E^{Q_{1}}\bigl[r^{2}\bigr]-1
\bigr) ^{1/2}%
\leq\varepsilon^{1/2}.
\end{eqnarray*}
\upqed\end{pf}

We also need the following approximations for the tail probability
$v(b)$. This proposition is an extension of Theorem 3.4 and Corollary
3.5 in \cite{LiuXu11}. We layout the key steps of its proof in the
supplemental material \cite{supp}.

\begin{proposition}\label{CorGRF}
Consider a Gaussian random field $\{f(t)\dvtx t\in T\}$ living on a
domain $T$ satisfying conditions \textup{C1--C6}. If $\mu(t)$ has one
unique maximum in $T$ denoted by $t_*$, then
\[
v(b) \sim(2\pi)^{d/2}\det\bigl(-\Delta\mu_\sigma(t_{\ast
})
\bigr)^{-1/2}G(t_{\ast})\cdot u^{d/2-1}\exp \biggl\{ -
\frac{(u-\mu
_{\sigma}
(t_{\ast}) )^{2}}{2} \biggr\},
\]
where
$u$ is as defined in \eqref{u}, and
$G(t)$ is defined as
\begin{eqnarray*}
&&\frac{\det(\Gamma)^{-{1}/{2}}}
{(2\pi)^{{(d+1)(d+2)}/{4}}} e^{{\mathbf{1}^T\mu_{22}\mathbf{1}}/{(8\sigma^2)}+B_t}
\\
&&\qquad{}\times\int_{R^{{d(d+1)}/{2}}}\exp \biggl\{
- \frac{1}{2} \biggl[\frac{|\mu_{20}\mu_{22}^{-1}z|^2}{1-\mu_{20}\mu
_{22}^{-1}\mu_{02}} + \biggl|\mu_{22}^{-1/2}z-
\frac{\mu_{22}^{1/2}\mathbf{1}}{2\sigma
}\biggr |^2 \biggr] \biggr\}\,dz.
\end{eqnarray*}
If $\mu(t)\equiv0$, $G(t)$ is a constant denoted by $G$. Then
\[
v(b) \sim \operatorname{mes}(T)G\cdot u^{d-1}e^{-{u^{2}}/{2}}.
\]
\end{proposition}

\subsection{Case 1: \texorpdfstring{$\mu(t)$}{mu(t)} is not a constant}

To make the proof smooth, we arrange the statement of the rest
supporting lemmas in the \hyperref[SecLem]{Appendix}.
We start the proof of Theorem \ref{main} when $\mu(t)$ is not a constant.
Note that
\[
E^{Q} \biggl[ \biggl(\frac{dP^{*}_{b}}{dQ} \biggr)^{2}
\biggr] =v(b)^{-2} E^{Q} \biggl[ \biggl(\frac{dP}{dQ}
\biggr)^{2} ;\mathcal{I}(T)>b \biggr].
\]
Thanks to Lemma \ref{LemTV1}, we only need to show that for any
$\varepsilon> 0$ there exists $b_0$ such that for all $b> b_0$
\begin{eqnarray*}
E^{Q} \biggl[ \biggl(\frac{dP}{dQ} \biggr)^{2} ;
\mathcal{I}(T)>b \biggr] &=& E^Q \biggl[E_{\imath,\tau}^{Q}
\biggl[ \biggl(\frac{dP}{dQ} \biggr)^{2} ;\mathcal{I}(T)>b
\biggr] \biggr]\leq(1+ \eps)v(b)^{2},
\end{eqnarray*}
where we use the notation $E_{\imath,\tau}^{Q}[ \cdot]=E^{Q}[ \cdot
| \imath,\tau]$ to denote the conditional expectation given $\imath$
and $\tau$. $\tau\in T$ is the random index described as in the
simulation scheme admitting a density function $l(t)$ if $\imath=0$ and
$\operatorname{mes}^{-1}(T)\mathbb{I}_{T}(t)$ if $\imath=1$.
Note that
\begin{eqnarray*}
&&E_{\imath,\tau}^{Q} \biggl[ \biggl(\frac{dP}{dQ}
\biggr)^{2} ;\mathcal{I}(T)>b \biggr] \\
&&\qquad= E_{\imath,\tau}^Q
\biggl[E_{\imath,\tau}^{Q} \biggl[ \biggl(\frac{dP}{dQ}
\biggr)^{2} ;\mathcal{I}(T)>b \Big| f(\tau),\partial f(\tau),\partial
^{2}f(\tau ) \biggr] \biggr].
\end{eqnarray*}
For the rest of the proof, we mostly focus on the conditional expectation
\[
E_{\imath,\tau}^{Q} \biggl[ \biggl(\frac{dP}{dQ}
\biggr)^{2}; \mathcal{I}(T)>b \Big| f(\tau),\partial f(\tau),
\partial^{2}f(\tau ) \biggr].
\]

The rest of the discussion is conditional on $\imath$ and $\tau$.
To simplify notation, for a given $\tau$, we define
\[
f_{*}(t) = f(t) - u_{\tau} C(t-\tau).
\]
On the set $\{\mathcal{I}(T)>b\}$, $f(\tau)$ reaches a level $u_\tau$,
and $E[f(t)| f(\tau)=u_\tau]=u_\tau C(t-\tau)$. Thus, $f_*(t)$ is the
field with the conditional expectation removed. From now on, we work
with this shifted field $f_*(t)$.
Correspondingly, we have
\[
\partial f_{*}(t) = \partial f(t)- u_{\tau} \partial C(t-
\tau),\qquad \partial^{2} f_{*}(t) = \partial^{2} f(t)-
u_{\tau} \partial ^{2}C(t-\tau ).
\]
We further define the following notation:
%
\begin{eqnarray}
\label{not} w&=& f_{*}(\tau),\qquad y=\partial f_{*}(\tau),\qquad
z=\partial^2 f_{*}(\tau),\qquad {\zz }=\Delta f_{*}(
\tau),
\nonumber
\\
{\tilde y}&=&\partial f_{*}(\tau)+\partial\mu_\sigma(\tau),\qquad
{\tilde\zz } = \Delta f_{*}(\tau)+\mu_{\sigma}(\tau)I+\Delta
\mu_{\sigma
}(\tau ),
\\
\qquad w_{t}&=&f_{*}(t),\qquad y_{t}=\partial
f_{*}(t),\qquad z_{t}=\partial^2
f_{*}(t),\qquad \bar z_t = \partial^2
f_{*}(t)- u_{t}\mu_{02}.\nonumber
\end{eqnarray}
Under the measure $Q$ and a given $\tau$, if $\imath=0$, $(w,y,z)$ has
density function
%
\begin{equation}
\label{hall*} h^*_{\mathrm{all}}(w,y,z)=\frac{1-\rho_1-{\rho_2}}{1-\rho_2}h_{0,\tau
}^*(w,y,z)+
\frac{\rho_1}{1-\rho_2}h^*_{1,\tau}(w,y,z);
\end{equation}
if $\imath= 1$, then $(w,y,z)$ follows density $h^*_\tau(w,y,z)$. The
forms of the densities can be derived from $h_{0,t}$, $h_{1,t}$ and
$h$. In particular, their expressions are given as follows:
\begin{eqnarray*}
h^*_{0,\tau}(w,y,z)&\propto&\mathbb{I}_{A_{\tau}} \times\exp \biggl\{
-\lambda u_\tau \biggl(w+\frac{\mathbf1^\top z}{2\sigma u_\tau
}+\frac
{B_{\tau}}{u_\tau}
\biggr)-\frac{1}{2}|y|^{2} \biggr\}
\\
&&{}\times\exp \biggl\{- \frac{1}{2} \biggl[\frac{|\mu_{20}\mu_{22}^{-1} z|^2}{1-\mu_{20}\mu
_{22}^{-1}\mu_{02}} + \biggl|
\mu_{22}^{-1/2} z-\frac{\mu_{22}^{1/2}\mathbf{1}}{2\sigma
}\biggr |^2 \biggr]
\biggr\},
\\
h^*_{1,\tau}(w,y,z)&\propto&\mathbb{I}_{A_{\tau}^c}\times\exp \biggl\{
\lambda_1 u_\tau \biggl(w+\frac{\mathbf1^\top z}{2\sigma u_\tau
}+
\frac
{B_{\tau}}{u_\tau} \biggr)-\frac{1}{2}|y|^{2} \biggr\}
\\
& &{}\times\exp \biggl\{- \frac{1}{2} \biggl[\frac{|\mu_{20}\mu_{22}^{-1}
z|^2}{1-\mu_{20}\mu_{22}^{-1}\mu_{02}} + \biggl|
\mu_{22}^{-1/2} z-\frac{\mu_{22}^{1/2}\mathbf{1}}{2\sigma
} \biggr|^2 \biggr]
\biggr\},
\\
h^*_\tau(w,y,z)&=& h(w,y,z)= \frac{\det(\Gamma)^{-({1}/{2})}}{(2\pi
)^{{(d+1)(d+2)}/{4}}}\\
&&\hspace*{56pt}{}\times \exp \biggl\{-
\frac{1}{2} \biggl[{y}^{\top}{y} +\frac{\llvert w-\mu_{20}\mu_{22}^{-1}{z}\rrvert ^2}{1-\mu_{20}\mu
_{22}^{-1}\mu_{02}}+{z}^{\top
}
\mu_{22}^{-1}{z} \biggr] \biggr\},
\end{eqnarray*}
and $A_{\tau}= \{w+\frac{y^{\top}y}{2u_\tau}+\frac{\mathbf
1^\top
z}{2\sigma u_\tau}+\frac{B_{\tau}}{u_\tau}>-\eta u_\tau^{-1}
\}$
is defined as in \eqref{At}.

In the next step, we will compute $dQ/ dP$ in the form of $f_{*}(t)$.
Basically, we
replace $f(t)$ by $f_{*}(t)+u_\tau C(t-\tau)$, $\partial f(t)$ by
$y_t + u_\tau\partial C(t -\tau)$, $\partial^{2} f(t)$ by
$z_t+u_\tau\partial^{2} C(t-\tau)$ and $\bar f''_t=\partial^2 f(t)-
u_{t}\mu_{02}$ by
$\bar z_t+u_\tau\partial^{2} C(t-\tau)$. For the likelihood ratio
terms $\operatorname{LR}$ and $\operatorname{LR}_1$ in \eqref{LR},
note that the $|\partial f(t)|^2$ terms in $h_{0,t}$ and $h_{1,t}$
cancel with those in $h(f(t), \partial f(t),\partial^2f(t))$, that is,
\begin{eqnarray*}
\operatorname{LR}(t)&=& \mathbb{I}_{A_t} \cdot { H_\lambda\cdot u_t} \exp\biggl\{-\lambda
u_t
\biggl(f(t)+\frac{\mathbf1^\top\bar f''_t}{2\sigma u_t}+\frac{B_t}{u_t}-u_t
\biggr)\\
&&\hspace*{75pt}{}- \frac{1}{2}
\biggl[\frac{|\mu_{20}\mu_{22}^{-1}\bar f''_t|^2}{1-\mu_{20}\mu
_{22}^{-1}\mu_{02}}
+ \biggl|\mu_{22}^{-1/2}\bar f''_t-\frac{\mu_{22}^{1/2}\mathbf
{1}}{2\sigma
} \biggr|^2 \biggr]\biggr\}\\
&&{}\Big/\biggl(\frac{\det(\Gamma)^{-{1}/{2}}}{(2\pi)^{
{(d+1)(d+2)}/{4}}}\\
&&\hspace*{13pt}{}\times e^{
-({1}/{2}) [{(f(t)-\mu_{20}\mu_{22}^{-1}{\partial^2
f(t)})^2}/{(1-\mu_{20}\mu_{22}^{-1}\mu_{02})}
+{\partial^2 f(t)}^{\top}\mu_{22}^{-1}{\partial^2f(t)}
]}\biggr).
\end{eqnarray*}
We insert the notation in \eqref{not} and obtain that
%
\begin{eqnarray}
\label{lr} \operatorname{LR}(t)&=& \mathbb{I}_{A_t}\cdot u_t
H_\lambda\exp \biggl\{-\lambda u_t \biggl(w_{t}+u_\tau
C(t-\tau)+\frac{\mathbf1^\top(\bar z_t+\mu_2(t-\tau
)u_\tau)}{2\sigma u_t}\nonumber\\
&&\hspace*{238pt}{}+\frac{B_{t}}{u_t}-u_t \biggr) \biggr\}
\nonumber
\\
&&\times \exp\biggl\{-\frac{1}{2} \biggl[\frac{|\mu_{20}\mu_{22}^{-1} (\bar
z_t+\mu
_2(t-\tau)u_\tau)|^2}{1-\mu_{20}\mu_{22}^{-1}\mu_{02}}
\\
&&\hspace*{52pt}{}+ \biggl|\mu_{22}^{-1/2}\bigl(\bar z_t+\mu_2(t-\tau)u_\tau\bigr)-\frac{\mu
_{22}^{1/2}\mathbf{1}}{2\sigma} \biggr|^2 \biggr]\biggr\}
\nonumber
\\
&&{}\times h^{-1}_{x,z} \bigl(w_{t} +
u_{\tau} C(t-\tau), z_t+u_\tau
\partial^{2} C(t-\tau) \bigr),\nonumber
\end{eqnarray}
where
%
\begin{equation}
\label{denh} h_{x,z}(x,z)=\frac{\det(\Gamma)^{-{1}/{2}}}{(2\pi)^{
{(d+1)(d+2)}/{4}}}e^{-%
({1}/{2}) [{(x-\mu_{20}\mu_{22}^{-1}{z})^{2}}/{%
(1-\mu_{20}\mu_{22}^{-1}\mu_{02})}+{z}^{\top}\mu_{22}^{-1}{z}
] },\hspace*{-35pt}
\end{equation}
which is the function $h(x,y,z)$ with the $|y|^2$ term removed.
Similarly, we have that
%
\begin{eqnarray}
\label{lr1} \operatorname{LR}_1(t) &=&\mathbb{I}_{A_t^c}\cdot
u_t H_{\lambda_{1}} \exp \biggl\{\lambda_{1}
u_t \biggl(w_{t}+u_\tau C(t-\tau)+
\frac{\mathbf1^\top(\bar z_t+\mu
_2(t-\tau)u_\tau)}{2\sigma u_t}\nonumber\\
&&\hspace*{239pt}{}+\frac{B_{t}}{u_t}-u_t \biggr) \biggr\}
\nonumber
\\
&&\times \exp\biggl\{-\frac{1}{2} \biggl[\frac{|\mu_{20}\mu_{22}^{-1} (\bar
z_t+\mu
_2(t-\tau)u_\tau)|^2}{1-\mu_{20}\mu_{22}^{-1}\mu_{02}}
\\
&&\hspace*{52pt}{}+ \biggl|\mu_{22}^{-1/2}\bigl(\bar z_t+\mu_2(t-\tau)u_\tau\bigr)-\frac{\mu
_{22}^{1/2}\mathbf{1}}{2\sigma} \biggr|^2 \biggr]\biggr\}
\nonumber
\\
&&{}\times h^{-1}_{x,z} \bigl(w_{t} +
u_{\tau} C(t-\tau), z_t+u_\tau
\partial^{2} C(t-\tau) \bigr).\nonumber
\end{eqnarray}
With the analytic forms \eqref{lr} and \eqref{lr1}, we proceed to the
likelihood ratio in \eqref{LRV}
%
\begin{eqnarray}
\label{Integral} \frac{dQ}{dP}&=& (1-{\rho_1}-{
\rho_2})K+\rho_1 K_1+\rho_2
K_2,
\end{eqnarray}
where
\begin{eqnarray*}
K &=& \int_{A^{*}}l(t) \operatorname{LR}(t) \,dt, \qquad K_{1} = \int
_{(A^{*})^{c}}l(t) \operatorname{LR}_1(t) \,dt, \\
K_{2} &=&\int
_T \frac{e^{-({1}/{2})u_t^2+u_t w_{t}+u_t u_\tau
C(t-\tau
)}}{\operatorname{mes}(T)} \,dt.
\end{eqnarray*}
The set $A^{*}$ [depending on the sample path $f_*(t)$] is defined as
\begin{eqnarray*}
&&\biggl\{t\dvtx w_{t}+ C(t-\tau)u_\tau+\frac{|y_t+u_\tau\cdot\partial C(t
-\tau)|^2}{2u_t}+
\frac{\mathbf1^\top(\bar z_t+u_\tau\mu_2(t-\tau
))}{2\sigma u_t}+\frac{B_{t}}{u_t} \\
&&\hspace*{294pt}{}>u_t-\frac{\eta}{u_t} \biggr
\}.
\end{eqnarray*}
We may equivalently define $A^* = \{t\dvtx f \in A_t\}$.
Note that $\operatorname{LR}(t) = 0$ if $f\notin A_t$. Thus, the integral $K$ is on
the set $A^*$, and $K_1$ is on the complement of $A^*$.

Based on the above results, we have that
%
\begin{eqnarray}
\label{part2}&& E^{Q} \biggl[ \biggl(\frac{dP}{dQ}
\biggr)^2;\mathcal{I}(T)>b \biggr]\nonumber \\
&&\qquad\leq E^{Q} \biggl
\{E^{Q}_{\imath,\tau} \biggl[ \frac{1}{ [(1-\rho
_1-\rho
_2)K+\rho_1 K_1 ]^2};\mathcal{I}(T)>b
\biggr] \biggr\}
\nonumber
\\[-8pt]
\\[-8pt]
\nonumber
&&\qquad\leq E^{Q} \biggl\{E^{Q}_{\imath,\tau} \biggl[
\frac{1}{ [(1-\rho
_1-\rho
_2)K ]^2};\mathcal{I}(T)>b, \mA_\tau\geq0 \biggr] \biggr\}
\nonumber
\\
&&\qquad\quad{} +E^{Q} \biggl\{E^{Q}_{\imath,\tau} \biggl[
\frac{1}{ [(1-\rho
_1-\rho
_2)K+\rho_1 K_1 ]^2};\mathcal{I}(T)>b, \mA_\tau<0 \biggr] \biggr\} ,
\nonumber
\end{eqnarray}
where
%
\begin{equation}
\label{mA} \mA_\tau=w+\frac{y^{\top}y}{2u_\tau}+\frac{\mathbf1^\top
z}{2\sigma
u_\tau}+
\frac{B_{\tau}}{u_\tau}.
\end{equation}
Note that the term $K_2$ is not used in the main analysis. In fact,
$K_2$ is only used in Lemma \ref{LX1} for the purpose of localization
that will be presented later.
The rest of the analysis consists of three main parts.

\textit{Part} 1. Conditional on $ (\imath,\tau,f_{*}(\tau
),\partial
f_{*}(\tau),\partial
^{2}f_{*}(\tau) )$, we study the event%
%
\begin{equation}
\mathcal{E}_{b}= \bigl\{ \mathcal{I}(T)>b \bigr\}, \label{cond}
\end{equation}
and write the occurrence of this event almost as a deterministic
function of $f_{*}(\tau)$, $\partial f_{*}(\tau)$ and $\partial
^{2}f_{*}(\tau)$, equivalently, $(w,y,z)$.

\textit{Part} 2.
Conditional on $ (\imath,\tau,f_{*}(\tau),\partial f_{*}(\tau
),\partial
^{2}f_{*}(\tau) )$, we express $K$ and $K_1$ as
functions of $f_*(\tau)$, $\partial f_*(\tau)$, $\partial
^{2}f_*(\tau
)$ with small correction terms.

\textit{Part} 3. We combine the results from the first two parts and
obtain an
approximation of \eqref{part2}.

All the subsequent derivations are conditional on $\imath$ and
$\tau$.

\subsubsection{Preliminary calculations}

To proceed, we provide the Taylor expansions for $f_{*}(t)$, $C(t)$ and
$\mu(t)$:
\begin{itemize}
\item Expansion of $f_{*}(t)$ given $ (f_{*}(\tau
),\partial f_{*}(\tau),\partial^{2}f_{*}(\tau) )$. Let $t-\tau
=((t-\tau)_{1},\ldots,\break  (t-\tau)_{d})$. Conditional on $ (f_{*}(\tau
),\partial f_{*}(\tau),\partial^{2}f_{*}(\tau) )$, we first expand
the random
function%
%
\begin{eqnarray}\label{expf}
f_{*}(t) &=&E \bigl[f_{*}(t)|f_{*}(\tau),
\partial f_{*}(\tau ),\partial ^{2}f_{*}(\tau )
\bigr]+g(t-\tau) \nonumber
\\
&=&f_{*}(\tau)+\partial f_{*}(\tau)^{\top}(t-
\tau)+\tfrac
{1}{2}(t-\tau)^{\top
}\Delta f_{*}(\tau) (t-
\tau)
\\
&&{}+R_{f}(t-\tau)+g(t-\tau),
\nonumber
\end{eqnarray}
where
\[
R_{f}(t-\tau)= O \bigl(|t|^{2+\delta_0} \bigl(\bigl|f_{*}(
\tau)\bigr|+\bigl|\partial f_{*}(\tau)\bigr|+\bigl|\partial^{2}f_{*}(
\tau )\bigr|\bigr) \bigr)
\]
%
is the remainder term of the Taylor expansion of $E
[f_{*}(t)|f_{*}(\tau),\partial f_{*}(\tau),\break  \partial^{2}f_{*}(\tau
) ]$.
$g(t)$ is a mean zero Gaussian random field such that
$Eg^{2}(t)=O(|t|^{4+\delta_0})$ as $t\rightarrow0$. In addition, the
distribution of $g(t)$ is independent of $\imath$, $\tau, f_{*}(\tau),
\partial
f_{*}(\tau)$ and $\partial^2 f_{*}(\tau)$.

\item Expansion of $C(t)$:
%
\begin{equation}
C(t)=1-\tfrac{1}{2}t^{\top}t+C_{4}(t)+R_{C}(t),
\label{expc}
\end{equation}
where
$C_{4}(t)=\frac{1}{24}\sum_{ijkl}\partial_{ijkl}^{4}C(0)t_{i}t_{j}t_{k}t_{l}$
and $R_{C}(t)=O(|t|^{4+\delta_0})$.
\item Expansion of $\mu(t)$:
%
\begin{equation}\qquad
\mu(t)=\mu(\tau)+\partial\mu(\tau)^{\top}(t-\tau)+\tfrac
{1}{2}(t-
\tau)^{\top}\Delta\mu(\tau) (t-\tau)+R_{\mu}(t-\tau),
\label{expmu}
\end{equation}
where $R_{\mu}(t-\tau)=O(|t-\tau|^{2+\delta_0})$.
\end{itemize}
We write
\[
R(t)=R_{f}(t)+u_\tau R_{C}(t)+R_\mu(t)/
\sigma
\]
to denote all the remainder terms.

Choose small constants $\epsilon$ and $\delta$ such that $0<\epsilon
\ll
\delta\ll\delta_0 $. By writing
\[
x\ll y,
\]
we mean that $x/y$ is chosen sufficiently small, but $x/y$ does not
change with $b$.
Let
%
\begin{eqnarray}
\label{LQ} 
\mathcal{L}&=& \Bigl\{ |\tau-t_*|<u^{-1/2+\epsilon},
\vert w\vert\leq u^{1/2+\epsilon}, 
|y|<u^{\epsilon},|z|<u^{\epsilon},
\nonumber
\\[-8pt]
\\[-8pt]
\nonumber
&& \hspace*{7pt}\sup_{|t-\tau|<u^{-1+\delta}}|z_t-z|<u^{-\epsilon}, \sup
_{|t-\tau|<u^{-1+\delta}} \bigl|g(t)\bigr| < u^{-1-\delta} \Bigr\}.
\end{eqnarray}
%
By Lemma \ref{LX1} whose proof uses the last component $\operatorname{LR}_2(t)$, we
have that
\begin{eqnarray*}
E^{Q} \biggl[ \biggl(\frac{dP}{dQ} \biggr)^2;
\mathcal{E}_{b},{\cal L}^c \biggr]&=& o(1)v^2(b).
\end{eqnarray*}
Therefore we only need to consider the second moment on the
set ${\cal L}$, that is,
%
\begin{eqnarray}
\label{der}&& E^{Q} \biggl[ \biggl(\frac{dP}{dQ}
\biggr)^2; \mathcal {E}_{b},\mathcal {L} \biggr]\nonumber \\
&&\qquad\leq
E^{Q} \biggl[E_{\imath,\tau}^{Q} \biggl[
\frac{1}{ [(1-\rho
_1-\rho
_2)K ]^2};\mathcal{E}_{b},\mathcal{L}, \mA_\tau>0
\biggr] \biggr]
\\
&&\qquad\quad{}+E^{Q} \biggl[E_{\imath,\tau}^{Q} \biggl[
\frac{1}{ [(1-\rho
_1-\rho
_2)K+\rho_1 K_1 ]^2};\mathcal{E}_{b},\mathcal{L}, \mA_\tau <0
\biggr] \biggr],
\nonumber
\end{eqnarray}
where $K$ and $K_1$ are given as in \eqref{Integral}.
We will focus on the terms on the right-hand side of \eqref{der} in the
subsequent derivations.
Now, we start to carry out each part of the program.

\subsection{Part 1}

All the derivations in this part are conditional on specific values
of $\imath$, $\tau$, $f_{*}(\tau)$, $\partial f_{*}(\tau)$ and
$\partial^{2}f_{*}(\tau
)$, equivalently, $\imath$, $\tau$, $w$, $y$ and $z$. By definition,
\[
\mathcal{I}(T)= \int_{T}e^{\sigma f_{*}(t)+\sigma u_\tau C(t-\tau
)+\mu(t)}\,dt.
\]
We insert the expansions in (\ref{expf}), (\ref{expc}) and \eqref{expmu}
into the expression of $\mathcal{I}(T)$ and obtain that
%
\begin{eqnarray}
\label{I1} \mathcal{I}(T)& =&\int_{t\in T}\exp \biggl\{\sigma
\biggl[ w+y^{\top
}(t-\tau )+\frac{1}{2}(t-\tau)^{\top}z(t-
\tau)+R_{f}(t-\tau)\nonumber\\
&&\hspace*{226pt}{}+g(t-\tau) \biggr] \biggr\}
\\
&&\hspace*{18pt}{}\times\exp \biggl\{ \bigl(\sigma u-\mu(\tau) \bigr)
\nonumber
\\
&&\hspace*{51pt}{}\times \biggl(1-
\frac
{1}{2}(t-\tau )^{\top}(t-\tau)+C_{4}(t-
\tau)+R_{C}(t-\tau) \biggr) \biggr\}
\nonumber\\
&&\hspace*{17pt}{}\times\exp \biggl\{\mu(\tau)+\partial\mu(\tau)^{\top}(t-\tau)+
\frac{1}{2}(t-\tau)^{\top}\Delta\mu(\tau) (t-\tau)\nonumber\\
&&\hspace*{215pt}{}+R_{\mu
}(t-
\tau ) \biggr\}\,dt,\nonumber
\end{eqnarray}
where the first row corresponds to the expansion of $w_t=f_{*}(t)$, and
the second and third rows correspond to those of $C(t)$ and $\mu(t)$,
respectively.
We write the exponent inside the integral in a quadratic form of
$(t-\tau)$ and obtain that%
%
\begin{eqnarray}
\label{l4} \mathcal{I}(T) &=&\exp \biggl\{ \sigma u+\sigma w+
\frac{\sigma}{2}{\tilde y}^{\top
}(uI-{\tilde\zz})^{-1}{\tilde
y} \biggr\}
\nonumber\\
&&\times\int_{t\in T}\exp \biggl\{ -\frac{\sigma}{2}\bigl(t-
\tau -(uI-{\tilde \zz})^{-1}{\tilde y}\bigr)^{\top}(uI-{\tilde
\zz})
\nonumber
\\[-8pt]
\\[-8pt]
\nonumber
&&\hspace*{102pt}{}\times \bigl( t-\tau -(uI-{\tilde \zz})^{-1}{\tilde y} \bigr) \biggr\}
\nonumber
\\
&&\hspace*{30pt}{} \times\exp\bigl\{ \sigma u_\tau C_{4} ( t-\tau ) +
\sigma R ( t-\tau )\bigr\} \times \exp\bigl\{ \sigma g ( t-\tau )\bigr\} \,dt,
\nonumber
\end{eqnarray}
where $\tilde y$ and $\tilde{\zz}$ are defined as in \eqref{not}.
Let $a(s)$ and $b(s)$ be two generic positive functions. Then we have
the representation of the following integral:
\[
\int_{T}a(s) b(s)\,ds = E\bigl[b(S)\bigr]\int
_{T} a(s)\,ds,
\]
where $S$ is a random variable taking values\vspace*{1pt} in $T$ with density
$a(s)/\int_{T} a(t)\,dt$.
Using this representation and the change of variable that
$s=(uI-{\tilde\zz})^{1/2}(t-\tau)$,
we write the big integral in \eqref{l4} as a product of expectations
and a
normalizing constant, and obtain that
\begin{eqnarray*}
\mathcal{I}(T) &=&\det(uI-{\tilde\zz})^{-1/2}\exp\biggl\{ \sigma u+
\sigma w +\frac{\sigma}{2}{\tilde y}^{\top}(uI-{\tilde
\zz})^{-1}{\tilde y}\biggr\}
\\
&&\times\int_{(uI-\mathbf{z})^{-{1}/{2}}s+\tau\in T}
\exp \biggl\{ -\frac{\sigma}{2}
\bigl(s-(uI-{\tilde\zz})^{-1/2}{\tilde y}\bigr)^{\top}
\\
&&\hspace*{116pt}{}\times  \bigl(
s-(uI-{\tilde\zz})^{-1/2}{\tilde y} \bigr) \biggr\} \,ds
\\
&&{}\times E \bigl[ e^{ \sigma u_\tau C_{4} ( (uI-{\tilde\zz
})^{-{1}/{2}}S ) +\sigma R ( (uI-{\tilde\zz})^{-{1}/{2}}S ) } \bigr] \times E
 \bigl[ e^{ \sigma g ( (uI-{\tilde\zz})^{-{1}/{2}}\tilde
{S} ) }
\bigr].
\end{eqnarray*}
The two expectations in the above display are taken with respect to $S$ and
$\tilde{S}$ given the process $g(t)$. $S$ is a random variable
taking values in the set $\{s\dvtx (uI-{\tilde\zz})^{-1/2}s+\tau\in
T\}$ with density proportional to
%
\begin{equation}
\label{densss} e^{ -({\sigma}/{2}) (s-(uI-{\tilde\zz})^{-1/2}{\tilde
y}
)^{\top} (s-(uI-{\tilde\zz})^{-1/2}{\tilde y} )},
\end{equation}
and $\tilde{S}$ is a random variable taking values in the set $ \{
s\dvtx (uI-{\tilde\zz})^{-1/2}s+\tau\in T \}$ with density
proportional to
\begin{eqnarray*}
&&e^{-({\sigma}/{2}) (s-(uI-{\tilde\zz})^{-1/2}{\tilde y} )^{\top
} (s-(uI-{\tilde\zz})^{-1/2}{\tilde y} )+\sigma u_\tau
C_{4} ((uI-{\tilde\zz})^{-{1}/{2}}s )+\sigma
R ((uI-{\tilde\zz})^{-{1}/{2}}s )}.
\end{eqnarray*}
Together with the definition of $u$ that
$ ( \frac{2\pi}{\sigma} ) ^{d/2}u^{-d/2}e^{\sigma u}=b$,
we obtain that
$\mathcal{I}(T) > b$
if and only if
%
\begin{eqnarray}
\label{ineq} \mathcal{I}(T) &=&\det(uI-{\tilde\zz})^{-1/2}e^{\sigma u+\sigma
w+({\sigma}/{2}){\tilde y}^{\top}(uI-{\tilde\zz})^{-1}{\tilde y}}
\nonumber
\\
&&\times\int_{(uI-\mathbf{z})^{-{1}/{2}}s+\tau\in T}
 e^{-({\sigma}/{2}) ( s-(uI-{\tilde\zz})^{-1/2}{\tilde
y} )
^{\top} ( s-(uI-{\tilde\zz}%
)^{-1/2}{\tilde y} )} \,ds
\nonumber
\\[-8pt]
\\[-8pt]
\nonumber
&&\times E \bigl[e^{\sigma u_\tau C_{4} ( (uI-{\tilde\zz
})^{-{1}/{2}}S )
+\sigma R ( (uI-{\tilde\zz})^{-{1}/{2}}S )} \bigr] \cdot e^{ -u^{-1}\xi_{u}}
\\
&>& \biggl( \frac{2\pi}{\sigma} \biggr) ^{d/2}u^{-d/2}e^{\sigma u},\nonumber
\end{eqnarray}
where
%
\begin{equation}
\label{xi} \xi_{u}=-u\log \bigl\{ E\exp \bigl[ \sigma g
\bigl((uI-{\tilde\zz })^{-{1}/{2}}\tilde{S}%
\bigr) \bigr] \bigr\}.
\end{equation}
We take log on both sides, and plug in the result of Lemma \ref{LemExp}
that handles the big expectation term in \eqref{ineq}.
Then inequality \eqref{ineq} is equivalent to
%
\begin{eqnarray}
\label{A'}&& w+\frac{{\tilde y}^{\top}(uI-{\tilde\zz})^{-1}{\tilde y}}{2}
-\frac{\log\det(I-{\tilde\zz}/{u})}{2\sigma}
+\frac{\sum_{i}\partial_{\mathit{iiii}}^{4}C(0)}{8\sigma^{2}u}
\nonumber
\\[-8pt]
\\[-8pt]
\nonumber
&&\qquad>
\frac{\xi_{u}}{u\sigma} + \frac{o(|w|+|y|+|z|+1)}{u^{1+\delta_0/4}}.
\end{eqnarray}
On the set $\mathcal L$, we further simplify \eqref{A'} using the
following facts (see Lemma~\ref{LemDet}):
\begin{eqnarray*}
\partial\mu_{\sigma}(\tau) &=& O\bigl(u^{-1/2+\epsilon}\bigr),
\\
\log\det \biggl(I-\frac{\tilde\zz}{u} \biggr) &=& -\frac
{1}{u}\operatorname{Tr}(\tilde {
\zz}) +o\bigl({u}^{-1-\delta_0/4}\bigr) \\
&=& - \frac{\bfone^{\top} (z +\partial^{2}\mu_{\sigma}(\tau)) +
d\cdot
\mu_{\sigma}(\tau)}{u} + o
\bigl({u}^{-1-\delta_0/4}\bigr),
\end{eqnarray*}
where $\operatorname{Tr}$ is the trace of a matrix.
Therefore, on the set ${\cal L}$, \eqref{A'} is equivalent to
\begin{eqnarray*}
\label{A}&& w+\frac{{y}^{\top}{y}}{2u} +\frac{{\mathbf{1}^{\top}(z+\partial^{2}\mu_\sigma(\tau))+d\cdot
\mu
_\sigma(\tau)}}{2\sigma u} +\frac{\sum_{i}\partial_{\mathit{iiii}}^{4}C(0)}
{8\sigma^{2}u}
\\
&&\qquad>
\frac{\xi_{u}}{u\sigma}+\frac{o(|w|+|y|+|z|+1)}{u^{1+\delta_0/4}},
\end{eqnarray*}
and further, equivalently (by replacing $u$ with $u_\tau$),
\begin{eqnarray*}
&&w+\frac{{y}^{\top}{y}}{2u_\tau} +\frac{{\mathbf{1}^{\top}(z+\partial^{2}\mu_\sigma(\tau))+d\cdot
\mu
_\sigma(\tau)}}{2\sigma u_\tau} +\frac{\sum_{i}\partial_{\mathit{iiii}}^{4}C(0)}
{8\sigma^{2}u_\tau}\\
&&\qquad >
\frac{\xi_{u}}{u\sigma}+\frac{o(|w|+|y|+|z|+1)}{u^{1+\delta_0/4}}.
\end{eqnarray*}
Using the notation defined as in \eqref{B} and \eqref{mA}, $\mathcal
I(T) > b$ is equivalent to
\[
\mA_\tau+\frac{o(|w|+|y|+|z|+1)}{u^{1+\delta_0/4}} >\frac{\xi_{u}}{u\sigma},
\]
where $\mA_\tau$ is defined as in \eqref{mA}. Furthermore, with
$\epsilon\ll\delta_0$ and on the set $\mathcal L$,
$o(|y|+|z|)/u^{-1-\delta_0/4} = o(u^{-1-\delta_0/8})$.
For the above inequality, we absorb $o(wu^{-1-\delta_0/4})$ into
$\mathcal A_\tau$ and rewrite it as
\[
\mathcal A_\tau> \bigl(1+o\bigl(u^{-1-\delta_0/4}\bigr)\bigr) \biggl[
\frac{\xi_u}{\sigma
u} + o\bigl(u^{-1-\delta_0/8}\bigr) \biggr].
\]

\subsection{Part 2}
In part 2, we first consider $(1-\rho_1-\rho_2)K$ in the first
expectation of \eqref{part2} (which is on the set $\{\mA_\tau\geq0\}$)
and then $(1-\rho_1-\rho_2)K+\rho_1 K_1$ in the second expectation of
\eqref{part2}.

\subsubsection*{Part 2.1: The analysis of $K$ when \texorpdfstring{$\mA_{\tau}\geq0$}{A tau >= 0}}
Similarly to part 1, all the derivations are conditional
on $(\imath, \tau, w, y, z)$. We now proceed to the second part of the
proof. More
precisely, we simplify the term $K$ defined as in (\ref{Integral}), and
write it as a deterministic function of $(w,y,z)$ with a small
correction term.
Recall that
\begin{eqnarray*}
K &=& \int_{A^{*}}l(t) u_t H_\lambda\exp
\biggl\{-\lambda u_t \biggl(w_{t}+u_\tau C(t-
\tau)+\frac{\mathbf1^\top(\bar z_t+\mu_2(t-\tau
)u_\tau)}{2\sigma u_t}\\
&&\hspace*{250pt}{}+\frac{B_{t}}{u_t}-u_t \biggr) \biggr\}
\\
&&\hspace*{13pt}{}\times\exp \biggl\{-\frac{1}{2} \biggl[\frac{|\mu_{20}\mu_{22}^{-1}
(\bar z_t+\mu_2(t-\tau)u_\tau)|^2}{1-\mu_{20}\mu_{22}^{-1}\mu_{02}} \\
&&\hspace*{68pt}{}+ \biggl|
\mu_{22}^{-1/2}\bigl(\bar z_t+
\mu_2(t-\tau)u_\tau\bigr)-\frac{\mu
_{22}^{1/2}\mathbf{1}}{2\sigma}
\biggr|^2 \biggr] \biggr\}
\\
&&\hspace*{13pt}{}\times h^{-1}_{x,z} \bigl(w_{t} +
u_{\tau} C(t-\tau), z_t+u_\tau
\partial^{2} C(t-\tau) \bigr)\,dt.
\end{eqnarray*}
We plug in the forms of $h_{x,z}$ and $l(t)$ that are defined in \eqref
{denh} and \eqref{lt} and obtain that
\begin{eqnarray*}
K &=&(2\pi)^{{(d+1)(d+2)}/{4}-{d}/{2}}\det(\Gamma)^{{1}/{2}} \cdot{\det\bigl(-\Delta
\mu_\sigma(t_*)\bigr)^{1/2}u_{t_*}^{d/2}}H_\lambda
\\
&&{}\times\int_{A^{*}}\exp \biggl\{\frac{u_{t_*}\cdot(t-t_*)^{\top
}\Delta
\mu_\sigma(t_*)(t-t_*)}{2} \biggr\}
\\
&&\hspace*{25pt}{}\times u_{t}\\
&&\hspace*{25pt}{}\times \exp \biggl\{-\lambda u_t
\biggl(w_{t}+u_\tau C(t-\tau)+\frac{\mathbf
1^\top(\bar z_t+\mu_2(t-\tau)u_\tau)}{2\sigma u_t}+
\frac
{B_{t}}{u_t}-u_t \biggr)\hspace*{-0.5pt} \biggr\}
\\
&&\hspace*{26pt}{} \times\exp \biggl\{-\frac{1}{2} \biggl[\frac{|\mu_{20}\mu_{22}^{-1}
(\bar z_t+\mu_2(t-\tau)u_\tau)|^2}{1-\mu_{20}\mu_{22}^{-1}\mu_{02}}\\
&&\hspace*{82pt}{} + \biggl|
\mu_{22}^{-1/2}\bigl(\bar z_t+
\mu_2(t-\tau)u_\tau\bigr)-\frac{\mu
_{22}^{1/2}\mathbf{1}}{2\sigma}
\biggr|^2 \biggr] \biggr\}
\\
&&\hspace*{26pt}{}\times\exp \biggl\{ \frac{1}{2} \biggl[\frac{ (w_{t}+u_\tau C(t-\tau)-\mu_{20}\mu
_{22}^{-1}({z_t+\mu_2(t-\tau)u_\tau}) )^2}{1-\mu_{20}\mu
_{22}^{-1}\mu_{02}}
\\
&&\hspace*{90pt} {}+\bigl({z_t+\mu_2(t-\tau)u_\tau}
\bigr)^{\top
}\mu_{22}^{-1}{\bigl(z_t+
\mu_2(t-\tau)u_\tau\bigr)} \biggr] \biggr\} \,dt.
\end{eqnarray*}
For some $\delta'$ such that $\epsilon<\delta'<\delta$, where
$\epsilon
,\delta$ are the parameters we used to define~${\cal L}$, we further
restrict the integration region by defining
%
\begin{eqnarray}
\label{I22} \mathcal{I}_2 &=&\int_{A^{*},|t-\tau|<u^{-1+\delta'}}\exp
\biggl\{\frac
{u_{t_*}(t-t_*)^{\top}\Delta\mu_\sigma(t_*)(t-t_*)}{2} \biggr\}\nonumber
\\
&&\hspace*{6pt}{}\times u_{t}\times \exp \biggl\{-\lambda u_t
\biggl(w_{t}+u_\tau C(t-\tau)\nonumber\\
&&\hspace*{92pt}{}+\frac{\mathbf
1^\top(\bar z_t+\mu_2(t-\tau)u_\tau)}{2\sigma u_t}+
\frac
{B_{t}}{u_t}-u_t \biggr) \biggr\}
\nonumber
\\
&&\hspace*{6pt}{} \times\exp \biggl\{-\frac{1}{2} \biggl[\frac{|\mu_{20}\mu_{22}^{-1}
(\bar z_t+\mu_2(t-\tau)u_\tau)|^2}{1-\mu_{20}\mu_{22}^{-1}\mu_{02}}\\
&&\hspace*{62pt}{} +\biggl |
\mu_{22}^{-1/2}\bigl(\bar z_t+
\mu_2(t-\tau)u_\tau\bigr)-\frac{\mu
_{22}^{1/2}\mathbf{1}}{2\sigma}
\biggr|^2 \biggr] \biggr\}
\nonumber
\\
&&\hspace*{6pt}{}\times\exp \biggl\{ \frac{1}{2} \biggl[\frac{ (w_{t}+u_\tau C(t-\tau)-\mu_{20}\mu
_{22}^{-1}({z_t+\mu_2(t-\tau)u_\tau}) )^2}{1-\mu_{20}\mu
_{22}^{-1}\mu_{02}}
\nonumber
\\
&&\hspace*{71pt}{} +\bigl({z_t+\mu_2(t-\tau)u_\tau}
\bigr)^{\top
}\mu_{22}^{-1}{\bigl(z_t+
\mu_2(t-\tau)u_\tau\bigr)} \biggr] \biggr\} \,dt.
\nonumber
\end{eqnarray}
Thus
\begin{eqnarray*}
K& \geq&(2\pi)^{{(d+1)(d+2)}/{4}-{d}/{2}}\det(\Gamma )^{{1}/{2}} \\
&&{}\times {\det\bigl(-
\Delta\mu_\sigma(t_*)\bigr)^{1/2}u_{t_*}^{d/2}}H_\lambda
\cdot \mathcal{I}_2.
\end{eqnarray*}
For the rest of part 2.1, we focus on $\mathcal{I}_2$. With some
tedious algebra, Lemma \ref{LemI2} writes $\mathcal{I}_2$ in a more
manageable form; that is, $\mathcal{I}_2$ equals
%
\begin{eqnarray}
\label{I2} && \int_{A^{*},|t-\tau|<u^{-1+\delta'}} \exp \biggl\{
\frac
{u_{t_*}(t-t_*)^{\top}\Delta\mu_\sigma(t_*)(t-t_*)}{2}+\frac
{u^{2}_{t}}{2} \biggr\}\times u_{t}
\nonumber\\
&&\hspace*{6pt}{}\times\exp \biggl\{(1-\lambda) u_t \bigl[w_{t}+
u_\tau C(t-\tau )-u_t \bigr]\nonumber\\
&&\hspace*{41pt}{}+ \frac{(1-\lambda)}{2\sigma}
\mathbf1^\top \bigl(z_t-\mu_{02}u_t+
\mu_2(t-\tau)u_\tau \bigr) -\lambda B_{t}-
\frac{\mathbf1^\top\mu_{22}\mathbf1}{8\sigma
^2} \biggr\}
\nonumber
\\[-6pt]
\\[-6pt]
\nonumber
&&\hspace*{6pt}{} \times\exp \bigl\{\bigl(\bigl (w_{t}+u_\tau C(t-\tau)-u_t \bigr)^2
\\
&&\hspace*{43pt}{}-2 \bigl(w_{t}+u_\tau C(t-\tau)-u_t \bigr)\mu_{20}\mu_{22}^{-1}{
\bigl(z_t-\mu
_{02}u_t+\mu_2(t-\tau)u_\tau \bigr)}\bigr)\nonumber\\
&&\hspace*{232pt}{}/{\bigl(2 \bigl(1-\mu_{20}\mu
_{22}^{-1}\mu
_{02} \bigr)\bigr)} \bigr\}\,dt.
\nonumber
\end{eqnarray}
Lemma \ref{LemA} implies that $\{|t-\tau|<u^{-1+\delta'}\}\subset A^*$.
Thus, on the set $\{\mA_\tau>0\}$, we have $A^*\cap\{|t-\tau
|<u^{-1+\delta'}\}=\{|t-\tau|<u^{-1+\delta'}\}$ and we can remove $A^*$
from the integration region of $\mathcal I_2$.
In addition, on the set ${\cal L}$ and $|t-\tau|< u^{-1+\delta'}$, we
have that
\begin{eqnarray*}
u_\tau-u_t C(t-\tau)&=&O\bigl(u^{-1 + 2\delta'}\bigr),\qquad
\mu_2(t-\tau)=\mu _{20}+O\bigl(|t-\tau|^2
\bigr),
\\
\bigl|u_\tau\mu_2(t-\tau)-u_t
\mu_{20}\bigr|&=&O\bigl(u^{-1+2\delta'}\bigr),\qquad \bigl(u_\tau-u_t
C(t-\tau)\bigr)|z_t|=o(1).
\end{eqnarray*}
We insert the above estimates to \eqref{I2}.
Together with the fact that
\[
\exp \biggl\{\frac{u_{t_*}(t-t_*)^{\top}\Delta\mu_\sigma
(t_*)(t-t_*)}{2}+\frac{u^{2}_{t}}{2} \biggr\} =\bigl(1+o(1)
\bigr)\exp \biggl\{\frac{1}{2}u_{t_*}^2 \biggr\},
\]
we have that
\begin{eqnarray*}
\mathcal{I}_2 &\sim& u \times\exp \biggl\{\frac{1}{2}u_{t_*}^2-
\lambda B_{t_{*}}- \frac
{\mathbf1^\top\mu_{22}\mathbf1}{8\sigma^2} \biggr\}
\\
&&{}\times \int_{|t-\tau|<u^{-1+\delta'}} \exp \biggl\{(1-\lambda)\\
&&\hspace*{91pt}{}\times u_t
\bigl[w_{t}+ u_\tau\bigl(C(t-\tau)-1\bigr)+\bigl(
\mu_\sigma(t)-\mu_\sigma(\tau)\bigr) \bigr]
\\
&&\hspace*{99pt}{} + (1-\lambda)\frac{\mathbf1^\top z}{2\sigma} +\frac
{w_{t}^{2}-2w_{t}\mu_{20}\mu_{22}^{-1}z_t+o(1)w_{t}}{2 (1-\mu
_{20}\mu_{22}^{-1}\mu_{02} )} \biggr\}\,dt.
\end{eqnarray*}
Further, we have that
\[
w_{t}^{2}-2w_{t}\mu_{20}
\mu_{22}^{-1}z_t+o(1)w_{t}=o(1)+u
\cdot w\cdot O\bigl(u^{-1/2+\delta' }\bigr).
\]
Let $\zeta_{u}=O(u^{-1/2+\delta' })$, and we simplify $\mathcal I_2$ to
\begin{eqnarray*}
\mathcal I_2 &\sim& u \times\exp \biggl\{\frac{1}{2}u_{t_*}^2-
\lambda B_{t_{*}}- \frac
{\mathbf1^\top\mu_{22}\mathbf1}{8\sigma^2} \biggr\}
\\
&&{}\times \int_{|t-\tau|<u^{-1+\delta'}} \exp \biggl\{(1-\lambda)
(u_\tau +\zeta _u) \bigl[\zeta_u w +
w_{t}+ u_\tau\bigl(C(t-\tau)-1\bigr)
\\
&&\hspace*{144pt}{} +\bigl(\mu_\sigma(t)-\mu_\sigma(\tau)\bigr) \bigr]+ (1-
\lambda)\frac
{\mathbf
1^\top z}{2\sigma} \biggr\}\,dt.
\end{eqnarray*}
In what follows, we insert the
expansions in (\ref{expf}), (\ref{expc}) and \eqref{expmu} into the
expression of $\mathcal{I}_2$ and write the exponent as a quadratic
function of $t-\tau$, and we obtain that on the set $\mathcal{L}$
%
\begin{eqnarray}
\label{LRD} \mathcal{I}_2 &\sim& u \times\exp \biggl\{
\frac{1}{2}u_{t_*}^2-\lambda B_{t_{*}}-
\frac
{\mathbf1^\top\mu_{22}\mathbf1}{8\sigma^2} \biggr\}
\nonumber
\\
&&{}\times\exp \biggl\{(1-\lambda) (u_\tau+\zeta_{u})
\biggl((1+\zeta _{u})w+\frac{1}{2}{\tilde y}^{\top}({u}I-{
\tilde\zz})^{-1}{\tilde y} + \frac{\mathbf1^\top z}{2\sigma u_\tau} \biggr) \biggr\}
\nonumber
\\[-8pt]
\\[-8pt]
\nonumber
&&{}\times\int_{|t-\tau|<u^{-1+\delta'}} e^{-({1}/{2})(1-\lambda)(u_\tau+\zeta_{u}) (t-\tau
-({u}I-{\tilde\zz})^{-1}{{\tilde y}} )^{\top}
({u}I-{\tilde
\zz} )
(t-\tau-({u}I-{\tilde\zz})^{-1}{{\tilde y}} )}
\\
&&\hspace*{67pt}{}\times e^{(1-\lambda)(u_\tau+\zeta_{u}) [{u_\tau}C_{4}(t-\tau
)+R(t-\tau)+g(t-\tau)  ]} \,dt,\nonumber
\end{eqnarray}
where we recall that ${\tilde y}=y+\partial\mu_\sigma(\tau)$ and
${\tilde\zz}=\zz+u_\sigma(\tau)I+\Delta\mu_\sigma(\tau)$.
This derivation is very similar to that from \eqref{I1} to \eqref{l4}.
In the last row of the above display, on the set $\mL$ and $|t-\tau|<
u^{-1+\delta'} $,
\[
u^{2}C_{4}(t-\tau) +u R(t-\tau) = o(1).
\]
Therefore, they can be ignored.
We consider the change of variable that
\[
s=(1-\lambda)^{1/2}(u_\tau+\zeta_{u})^{1/2}({u}I-{
\tilde\zz} )^{1/2}(t-\tau)
\]
and obtain that $\mathcal{I}_2$ equals (with the terms $C_{4}$ and $R$ removed)
%
\begin{eqnarray}
\label{inter0}\label{cont} \mathcal{I}_2 
&\sim&(1-
\lambda)^{-d/2}u^{-d+1}\exp \biggl\{\frac
{1}{2}u_{t_*}^2-
\lambda B_{t_{*}}- \frac{\mathbf1^\top\mu_{22}\mathbf1}{8\sigma^2} \biggr\}
\nonumber
\\
&&\times\exp \biggl\{(1-\lambda) (u_\tau+\zeta_{u})
\biggl((1+\zeta _{u})w+\frac{1}{2}{\tilde y}^{\top}({u}I-{
\tilde\zz})^{-1}{\tilde y} + \frac{\mathbf1^\top z}{2\sigma u} \biggr) \biggr\}
\nonumber
\\[-8pt]
\\[-8pt]
\nonumber
&&\times \int_{s\in\mathcal S_{u}} e^{-({1}/{2})\llvert  s-(1-\lambda)^{1/2}(u_\tau+\zeta_u)^{1/2}({u}I-
{\tilde\zz})^{-1/2}{\tilde y}\rrvert  ^{2}} \,ds
\nonumber
\\
&&\times E \bigl[e^{(1-\lambda)({u_\tau}+\zeta_u)g ((1-\lambda
)^{-1/2}(u_\tau+\zeta_u)^{-1/2}
({u}I-{\tilde\zz})^{-1/2}S^{\prime} )} \bigr], \nonumber
\end{eqnarray}
where $\mathcal S_{u} = \{s\dvtx |{(1-\lambda)^{-1/2}(u_\tau+\zeta
_u)^{-1/2}}({u}I-\tilde{\zz})^{-1/2}s|<u^{-1+\delta'}\}$, and
$S^{\prime}$ is a random variable taking values on the set $\mathcal
S_{u}$ with density proportional to
\begin{eqnarray*}
&&e^{-({1}/{2})\llvert  s-(1-\lambda)^{1/2}(u_\tau+\zeta_u)^{1/2}({u}I-
{\tilde\zz})^{-1/2}{\tilde y}\rrvert  ^{2}}.
\nonumber
\end{eqnarray*}
We use $\kappa$ to denote the last two terms of \eqref{inter0}, that is,
\begin{eqnarray}
\label{kappa} \kappa &=&\int_{\mathcal S_{u}} e^{-({1}/{2})\llvert  s-(1-\lambda)^{1/2}(u_\tau+\zeta_u)^{1/2}({u}I-
{\tilde\zz})^{-1/2}{\tilde y}\rrvert  ^{2}} \,ds
\nonumber
\\[-8pt]
\\[-8pt]
\nonumber
&&\hspace*{14pt}{}\times E \bigl[e^{(1-\lambda)({u_\tau}+\zeta_u)g ((1-\lambda
)^{-1/2}(u_\tau+\zeta_u)^{-1/2}
({u}I-{\tilde\zz})^{-1/2}S^{\prime} )} \bigr].
\end{eqnarray}
It is helpful to keep in mind that $\kappa$ is approximately $(2\pi)^{d/2}$.
We insert $\kappa$ back to the expression of $\mathcal{I}_2$. Together
with the fact that ${\tilde y}^{\top}({u}I-{\tilde\zz})^{-1}{\tilde y}
= |\tilde y|^{2}/u + o(u^{-1})$, we have
%
\begin{eqnarray}
\label{ii2} \mathcal{I}_2&\sim&\kappa(1-\lambda)^{-d/2}u^{-d+1}
\exp \biggl\{ \frac
{1}{2}u_{t_*}^2-\lambda
B_{t_{*}}- \frac{\mathbf1^\top\mu
_{22}\mathbf
1}{8\sigma^2} \biggr\}
\nonumber
\\[-8pt]
\\[-8pt]
\nonumber
&&{} \times\exp \biggl\{(1-\lambda) (u_\tau+\zeta_{u})
\biggl((1+\zeta _{u})w+\frac{|\tilde y|^{2}}{2u_\tau} + \frac{\mathbf1^\top
z}{2\sigma
u_\tau}
\biggr) \biggr\}.
\nonumber
\end{eqnarray}
Thus, we have that on the set $\{\mA_{\tau} > 0\}$,
%
\begin{eqnarray}
\label{KK} K&\geq&(2\pi)^{{(d+1)(d+2)}/{4}-{d}/{2}}\det(\Gamma )^{{1}/{2}} \cdot{
\det\bigl(-\Delta\mu_\sigma(t_*)\bigr)^{1/2}u_{t_*}^{d/2}}H_\lambda
\cdot \mathcal{I}_2
\nonumber\\
&=&\bigl(\kappa+o(1)\bigr) (2\pi)^{{(d+1)(d+2)}/{4}-{d}/{2}} \det(\Gamma)^{{1}/{2}}
\nonumber\\
&&{}\times{\det\bigl(-\Delta\mu_\sigma (t_*)\bigr)^{1/2}}H_\lambda
\cdot(1-\lambda)^{-d/2}u^{-d/2+1}
\\
&&{}\times \exp \biggl\{\frac{1}{2}u_{t_*}^2-\lambda
B_{t_{*}}- \frac{\mathbf
1^\top\mu_{22}\mathbf1}{8\sigma^2}\nonumber\\
&&\hspace*{33pt}{}+(1-\lambda) (u_\tau+\zeta
_{u}) \biggl((1+\zeta_{u})w+\frac{|\tilde y|^{2}}{2u_\tau} +
\frac{\mathbf
1^\top
z}{2\sigma u_\tau} \biggr) \biggr\}.
\nonumber
\end{eqnarray}
We further insert the $\mA_\tau$ defined in \eqref{mA} into \eqref{KK}
and obtain that
\begin{eqnarray}
\label{KKK} K&\geq&\bigl(\kappa+o(1)\bigr) (2\pi)^{{(d+1)(d+2)}/{4}-{d}/{2}}\det
(\Gamma
)^{{1}/{2}} \nonumber\\
&&{}\times {\det\bigl(-\Delta\mu_\sigma(t_*)
\bigr)^{1/2}}H_\lambda \cdot(1-\lambda)^{-d/2}u^{-d/2+1}
\nonumber
\\[-8pt]
\\[-8pt]
\nonumber
&&\times \exp \biggl\{\frac{1}{2}u_{t_*}^2-
B_{t_{*}}- \frac{\mathbf1^\top
\mu
_{22}\mathbf1}{8\sigma^2}+(1-\lambda)u_\tau\bigl(1+o(1)
\bigr)\mA_\tau\\
&&\hspace*{132pt}{}+ (1-\lambda )\zeta_{u}\cdot \bigl(|\tilde
y|^{2}+|z| \bigr) \biggr\}.
\nonumber
\end{eqnarray}
\subsubsection*{Part 2.2: The analysis of $dP/dQ$ when $\mA_{\tau}<0$}\label{part22}
In this part, we focus mostly on the $K_1$ term, whose handling is very
similar to that of $K$. Therefore, we only list out the key steps.
For some large constant $M$, let
\[
D= \bigl\{\bigl|t-\tau-({u}I-{\tilde\zz})^{-1}{{\tilde y}}
\bigr|<Mu^{-1}\bigr\}
\]
that is, the dominating region of the integral.
We split the set $D = (A^*\cap D) \cup((A^*)^c\cap D)$. There are two
situations: $\operatorname{mes}((A^*)^c\cap D)> \operatorname{mes} (A^*\cap D)$ and $\operatorname{mes}((A^*)^c\cap
D)\leq \operatorname{mes} (A^*\cap D)$.
For the first situation, the term $K_1$ is dominating; for the second
situation, the term $K$ (more precisely $\mathcal I_2$) is dominating.

To simplify $K_1$, we write it as
\begin{eqnarray*}
K_1 &=&(2\pi)^{{(d+1)(d+2)}/{4}-{d}/{2}}{\det(\Gamma )^{{1}/{2}}} \cdot{
\det\bigl(-\Delta\mu_\sigma(t_*)\bigr)^{1/2}u_{t_*}^{d/2}}H_{\lambda_1}
\\
&&{}\times \biggl[\int_{(A^{*})^c\cap D}+\cdots+\int_{(A^{*})^c\cap
D^{c}}\cdots\biggr]
\\
&\triangleq&(2\pi)^{{(d+1)(d+2)}/{4}-{d}/{2}}\det(\Gamma )^{{1}/{2}} \cdot{\det
\bigl(-\Delta\mu_\sigma(t_*)\bigr)^{1/2}u_{t_*}^{d/2}}H_{\lambda_1}\\
&&{}\times [\mathcal{I}_{1,2}+\mathcal{I}_{1,3} ].
\end{eqnarray*}
Note that the difference between $K_{1}$ and $K$ is that the term
``$-\lambda$'' has been replaced by ``$\lambda_{1}$.'' With exactly the
same derivation for \eqref{I2}, we obtain that $\mathcal{I}_{1,2}$
equals [by replacing ``$-\lambda$'' in \eqref{I2} by ``$\lambda_1$'']
%
\begin{eqnarray}
\label{I2'} &&\int_{(A^*)^c\cap D } \exp \biggl\{\frac{u_{t_*}(t-t_*)^{\top
}\Delta
\mu_\sigma(t_*)(t-t_*)}{2}+
\frac{1}{2}u_t^2 \biggr\}\times u_t
\nonumber\\
&&\qquad{}\times\exp \biggl\{(1+\lambda_1) u_t
\bigl[w_t+ u_\tau C(t-\tau )-u_t \bigr]\nonumber\\
&&\hspace*{56pt}{}+
\frac{(1+\lambda_1)}{2\sigma}\mathbf1^\top \bigl(z_t-\mu
_{02}u_t+\mu _2(t-\tau)u_\tau
\bigr)  +\lambda_1 B_t- \frac{\mathbf1^\top\mu
_{22}\mathbf1}{8\sigma^2} \biggr\}
\nonumber
\\
&&\qquad{}\times\exp \bigl\{\bigl( \bigl(w_t+u_\tau C(t-\tau)-u_t \bigr)^2
\\
&&\hspace*{57pt}{}-2 \bigl(w_t+\bigl(u_\tau C(t-\tau)-u_t\bigr) \bigr)\nonumber\\
&&\hspace*{57pt}{}\times \mu_{20}\mu_{22}^{-1}{
\bigl(z_t-\mu
_{02}u_t+\mu_2(t-\tau)u_\tau \bigr)}\bigr)\nonumber\\
&&\hspace*{135pt}{}/ \bigl(
2 \bigl(1-\mu_{20}\mu_{22}^{-1}\mu_{02} \bigr)\bigr) \bigr\}\,dt.
\nonumber
\end{eqnarray}
With a very similar derivation as in part 2.1, in particular, the
result in \eqref{LRD}, we have that
%
\begin{eqnarray}
\label{I121} \mathcal{I}_{1,2} &\sim&u \times\exp \biggl\{
\frac{1}{2}u_{t_*}^2+\lambda_1
B_{t_{*}}- \frac
{\mathbf1^\top\mu_{22}\mathbf1}{8\sigma^2} \biggr\}
\nonumber\hspace*{-25pt}\\
&&\times\exp \biggl\{(1+\lambda_1) (u_\tau+
\zeta_{u}) \biggl((1+\zeta _{u})w+\frac{1}{2}{\tilde
y}^{\top}({u}I-{\tilde\zz})^{-1}{\tilde y} +
\frac{\mathbf1^\top z}{2\sigma u} \biggr) \biggr\}
\nonumber\hspace*{-25pt}
\\
&& \times\int_{(A^{*})^c\cap D} \exp \biggl\{(1+\lambda_1)
(u_\tau+\zeta_{u}) \biggl[-\frac{1}{2} \bigl(t-\tau
-({u}I-{\tilde\zz})^{-1}{{\tilde y}} \bigr)^{\top} ({u}I-{
\tilde \zz} ) \hspace*{-25pt}\\
&&\hspace*{204pt}{}\times \bigl(t-\tau-({u}I-{\tilde\zz})^{-1}{{\tilde y}} \bigr)
\biggr] \biggr\}
\nonumber\hspace*{-25pt}
\\
&&\times\exp \bigl\{(1+\lambda_1) (u_\tau+
\zeta_{u}) \bigl[{u_\tau}C_{4}(t-\tau)+R(t-
\tau)+g(t-\tau) \bigr] \bigr\} \,dt.
\nonumber\hspace*{-25pt}
\end{eqnarray}
Furthermore, similarly to the results in \eqref{ii2}, we have that
%
\begin{eqnarray}
\mathcal{I}_{1,2} &\sim&
\kappa_{1,2}(1+
\lambda_1)^{-{d}/{2}}u^{-d+1} e^{({1}/{2})u_{t_*}^2+\lambda_1 B_{t_{*}}-
{\mathbf1^\top
\mu
_{22}\mathbf1}/{(8\sigma^2)}}
\nonumber
\\[-8pt]
\\[-8pt]
\nonumber
&&{}\times e^{(1+\lambda_1)(u_\tau+\zeta_{u}) ((1+\zeta
_{u})w+
({1}/{2}){\tilde y}^{\top}({u}I-{\tilde\zz})^{-1}{\tilde y}
+ {\mathbf1^\top z}/{(2\sigma u_\tau)} )},
\end{eqnarray}
where $\kappa_{1,2}$ is defined as
\begin{eqnarray*}
\label{k121def} \kappa_{1,2}&=&\int_{t_1(s)\in(A^*)^c\cap D}
e^{-{1}/{2}\llvert  s-(1+\lambda_1)^{1/2}(u_\tau+\zeta_u)^{1/2}({u}I-
{\tilde\zz})^{-1/2}{\tilde y}\rrvert  ^{2}}\,ds
\\
&&\times E \bigl[e^{(1+\lambda_1)({u_\tau}+\zeta_u)g
((1+\lambda
_1)^{-1/2}(u_\tau+\zeta_u)^{-1/2}
({u}I-{\tilde\zz})^{-1/2}S_{1,2} )} \bigr],
\end{eqnarray*}
the change of variable $t_1(s) =\tau+ {(1+\lambda_1)^{-1/2}(u_\tau
+\zeta
_u)^{-1/2}}({u}I-\tilde{\zz})^{-1/2}s$
and $S_{1,2}$ is a random variable taking values in the set $\{s\dvtx t(s)\in(A^*)^c\cap D\}$ with an appropriately chosen density function
similarly as in \eqref{inter0}.
In summary, the only difference between $\mathcal{I}_{1,2}$ and
$\mathcal{I}_2$ lies in that the multiplier $-\lambda$ is replaced by
$\lambda_{1}$.

We now proceed to providing a lower bound of $(1-\rho_1 - \rho_2) K +
\rho_1 K_1$.
Note that
\[
\max\bigl\{\operatorname{mes}\bigl(\bigl(A^{*}\bigr)^c\cap D\bigr), \operatorname{mes}
\bigl(A^{*}\cap D\bigr)\bigr\} \geq\tfrac{1}{2} \operatorname{mes}(D).
\]
Therefore at least one of $(A^{*})^c\cap D$ and $A^{*}\cap D$ is nonempty.
If $\operatorname{mes}((A^{*})^c\cap D)
\geq\frac{1}{2}\operatorname{mes}(D)$, we have the bound
\[
(1-\rho_1 - \rho_2) K + \rho_1
K_1\geq\rho_1K_1\geq\Theta(1)
\rho_1 u^{d/2} \mathcal I_{1,2}.
\]
Similarly, if $\operatorname{mes}(A^{*}\cap D)
\geq\frac{1}{2}\operatorname{mes}(D)$, we have that
\[
(1-\rho_1 - \rho_2) K + \rho_1
K_1 \geq\Theta(1) (1-\rho_1 - \rho_2)
u^{d/2} \mathcal I_{2}.
\]
We further split $\mathcal{I}_2$ in part 2.1 into two parts:
%
\begin{equation}
\mathcal{I}_2 = \int_{A^*\cap D}\cdots \,dt +\int
_{A^*\cap
D^c}\cdots \,dt \triangleq\mathcal{I}_{2,1}+
\mathcal{I}_{2,2}.
\end{equation}
Similarly to the derivation of $\mathcal I_{1,2}$, we have that
\begin{eqnarray*}
\label{211} \mathcal{I}_{2,1} &\sim&\kappa_{2,1}(1-
\lambda)^{-d/2}u^{-d+1} e^{
({1}/{2})u_{t_*}^2-\lambda B_{t_{*}}- {\mathbf1^\top\mu
_{22}\mathbf
1}/{(8\sigma^2)}}\\
&&{}\times
e^{(1-\lambda)(u_\tau+\zeta_{u}) ((1+\zeta
_{u})w+{|\tilde y|^2}/{(2u_\tau)}+ {\mathbf1^\top z}/{(2\sigma
u_\tau)} )},
\end{eqnarray*}
where
%
\begin{eqnarray}
\label{k211def} \kappa_{2,1} &=&\int_{t_2(s) \in A^*\cap D}
e^{-{1}/{2}\llvert  s-(1-\lambda)^{1/2}(u_\tau+\zeta_u)^{1/2}({u}I-
{\tilde\zz})^{-1/2}{\tilde y}\rrvert  ^{2}} \,ds
\nonumber
\\[-8pt]
\\[-8pt]
\nonumber
&&\times E \bigl[e^{(1-\lambda)({u_\tau}+\zeta_u)g ((1-\lambda
)^{-1/2}(u_\tau+\zeta_u)^{-1/2}
({u}I-{\tilde\zz})^{-1/2}S_{2,1} )} \bigr]. 
\end{eqnarray}
$S_{2,1}$ is a random variable taking values on the set $\{s\dvtx t_2(s)\in
A^{*}\cap D\}$ with an appropriate density function similarly as in
\eqref{inter0} and $t_2(s) =\tau+ (1-\lambda)^{-1/2}\times (u_\tau+\zeta
_u)^{-1/2}({u}I-\tilde{\zz})^{-1/2}s$.

Then combining the above results of $\mathcal{I}_{1,2}$ and $\mathcal
{I}_{2,1}$, we have that for the case in which $\mA_\tau<0$
\begin{eqnarray*}
&& \rho_1 K_1+(1-\rho_1-\rho_2)K
\\
&&\qquad\geq\Theta(1) u^{d/2} \bigl[\mathbb{I}_{C_{1}}
\rho_1 \mathcal{I}_{1,2}+\mathbb {I}_{C_{2}}(1-
\rho_1-\rho_2) \mathcal{I}_{2,1} \bigr]
\\
&&\qquad\geq \Theta(1) u^{-d/2+1} e^{({1}/{2})u_{t_*}^2} \\
&&\qquad\quad{}\times \bigl[
\mathbb{I}_{C_{1}}\cdot\rho_1\kappa_{1,2}
e^{(1+\lambda
_1)(u_\tau+\zeta_{u}) ((1+\zeta_{u})w+{|\tilde
y|^{2}}/{(2u_\tau)}+ {\mathbf1^\top z}/{(2\sigma u_\tau)} )}
\\
&&\hspace*{16pt}\qquad\quad{} + \mathbb{I}_{C_{2}}\cdot(1-\rho_1-\rho_2)
(1-\lambda)^{-d/2}\kappa _{2,1}\\
&&\hspace*{90pt}{}\times e^{(1-\lambda)(u_\tau+\zeta_{u})
 ((1+\zeta_{u})w+
{|\tilde y|^2}/{(2u_\tau)}+ {\mathbf1^\top z}/{(2\sigma u_\tau
)}
)} \bigr],
\end{eqnarray*}
where ${C_{1}}=\{f(\cdot)\dvtx \operatorname{mes}((A^{*})^c\cap D)\geq \operatorname{mes}(A^{*}\cap D)\}$
and ${C_{2}}=C_1^c$.
We further insert $\mA_\tau$ defined in \eqref{mA}. Note that on the
set $\{\mA_\tau<0\}$, $(1+\lambda_1)\mA_\tau<(1-\lambda)\mA_\tau
$ and
$B_t$ is bounded away from zero and infinity. Then
%
\begin{eqnarray}
\label{3star} &&(1-\rho_1-\rho_2)K+\rho_1
K_1 \nonumber\\
&&\qquad\geq \Theta(1) u^{-d/2+1} e^{({1}/{2})u_{t_*}^2} \cdot
e^{(1+\lambda
_1)(1+\zeta_{u})u_\tau\mA_\tau+ \zeta_{u}\cdot (|\tilde
y|^{2}+|z| )}
\\
&&\qquad\quad{}\times \bigl[\mathbb{I}_{C_{1}}\cdot\rho_1
\kappa_{1,2}+ \mathbb {I}_{C_{2}}\cdot (1-\rho_1-
\rho_2) (1-\lambda)^{-d/2}\kappa_{2,1} \bigr].\nonumber
\end{eqnarray}

\subsection*{Part 3}


We now put together the results in parts 1 and 2 and obtain an
approximation for \eqref{part2}. Recall that
%
\begin{eqnarray} \label{part3cal}
&&E^{Q} \biggl[ \biggl( \frac{dP}{dQ} \biggr) ^{2};
\mathcal{%
E}_{b},\mathcal{L} \biggr]\nonumber\\[-2pt]
&&\qquad
\leq E^{Q} \biggl[ \frac{1}{ [ (1-\rho_{1}-\rho_{2})K ]
^{2}};\mathcal{E}_{b},
\mathcal{L},\mA_{\tau}\geq0 \biggr]
\\[-2pt]
&&\qquad\quad{}+E^{Q} \biggl[ \frac{1}{ [
(1-\rho_{1}-\rho_{2})K+\rho_{1}K_{1} ] ^{2}};\mathcal{E}_{b},%
\mathcal{L},\mA_{\tau}<0 \biggr].\nonumber
\end{eqnarray}
We consider the two terms on the right-hand side of the above display
one by
one. We start with the first term
%
\begin{eqnarray}
\label{1st} && E^{Q} \biggl[ \frac{1}{ [ (1-\rho_{1}-\rho_{2})K ] ^{2}};
\mathcal{E}_{b},\mathcal{L},\mA_{\tau}\geq0 \biggr]
\nonumber\\[-2pt]
&&\qquad= E^{Q} \biggl[ \frac{1}{ [ (1-\rho_{1}-\rho_{2})K ] ^{2}}; \mathcal{E}_{b},
\mathcal{L},\mA_{\tau}\geq0, \imath=0 \biggr]
\\[-2pt]
&&\qquad\quad{}+ E^{Q} \biggl[ \frac{1}{ [ (1-\rho_{1}-\rho_{2})K ]
^{2}};\mathcal{E}_{b},
\mathcal{L} ,\mA_{\tau}\geq0, \imath=1 \biggr].
\nonumber
\end{eqnarray}
The index $\tau$ admits density $l(t)$ when $\imath=0$ and $\tau$ is
uniformly
distributed over $T$ if $\imath=1$.

Consider the first expectation in \eqref{1st}. Note that conditionally
on $\tau$ and $\imath=0$, on the set $\mathcal{L}\cap
\{\mA_{\tau}\geq0\}$, $(w,y,z)$ follows density $(1-\rho_1-\rho
_2)h_{0,\tau}^{\ast}(w,y,z)/(1-\rho_2)$ defined as in
(\ref{hall*}). Thus, according to \eqref{KKK}, we have that the
conditional expectation
%
\begin{eqnarray}\label{111}
&&E^{Q} \biggl[ \frac{1}{(1-\rho_{1}-\rho_{2})^{2}K^{2}};\mathcal{E}%
_{b},\mathcal{L},\mA_{\tau}\geq0 \Big\vert\imath=0,\tau \biggr]
\nonumber\\[-2pt]
&&\qquad\leq\bigl(1+o(1)\bigr) \biggl[ \frac{H_{\lambda}^{-1}\det(\Gamma)
^{-{1}/{2}}\det
(-\Delta\mu_{\sigma}(t_{\ast}))^{-1/2}}{(2\pi)^{
{(d+1)(d+2)}/{4}-%
{d}/{2}}}\nonumber\\[-2pt]
&&\hspace*{82pt}{}\times{(1-\lambda)^{d/2}u^{d/2-1}}e^{-({1}/{2})u_{t_*}^{2}+B_{t_{\ast}}
+{\mathbf{1}^{\top}\mu
_{22}\mathbf
{1}}/{%
(8\sigma^{2})}}
\biggr] ^{2}
\\[-2pt]
&&\qquad\quad{}\times\int_{{\cal A}_{\tau}>0,\mathcal{L}}e^{-2(1-\lambda
)u ( (1+o(1))\mA_{\tau}+o({|y|^{2}}/{(2u)}+{\mathbf
{1}^{\top
}z}/{(2\sigma u)}) ) }\cdot
\gamma_{u}(u\sigma\mA_\tau) \nonumber\\[-2pt]
&&\hspace*{81pt}{}\times \frac{1-\rho_1-\rho_2}{1-\rho_2}h_{0,\tau}^{\ast}
(w,y,z)\,dw\,dy\,dz,
\nonumber
\end{eqnarray}
where
\begin{eqnarray*}
&&\gamma_{u}(x)=
E \biggl[  \frac{1}{(1-\rho_{1}-{\rho_{2}})^{2}\kappa
^{2}};\\[-2pt]
&&\hspace*{53pt} x>
\bigl(1+o\bigl(u^{-1-\delta_0/4}\bigr)\bigr)\bigl[\xi_u + o
\bigl(u^{-\delta_0/8}\bigr)\bigr]\Big\rrvert \imath, \tau,w,{y},{z} \biggr],
\end{eqnarray*}
with the expectation taken with respect to the process $g(t)$. We insert
the analytic form of $h_{0,\tau}^{\ast}(w,y,z)$ into \eqref{hall*} and
obtain that
%
\begin{eqnarray}
\label{cc} &&\int_{{\cal A}_{\tau}>0,\mathcal{L}}e^{-2(1-\lambda)u (
(1+o(1))\mA_{\tau}+o({|y|^{2}}/{(2u)}+{\mathbf{1}^{\top}z}/{
(2\sigma u)}) )} \cdot
\gamma_{u}(u\sigma\mA_\tau)\nonumber\\
&&\hspace*{24pt}\quad{}\times \frac
{1-\rho
_1-\rho_2}{1-\rho_2}
h_{0,\tau}^{\ast}(w,y,z)\,dw\,dy\,dz
\nonumber
\\
&&\qquad=\frac{(1-\rho_{1}-{\rho_{2}})H_{\lambda}\cdot u_{\tau}}{1-\rho
_2}
\nonumber
\\
&&\qquad\quad{}\times \int_{\mA_{\tau
}>0}\gamma_{u}(u
\sigma\mA_{\tau})\exp \bigl\{ -2\bigl(1-\lambda+o(1) \bigr)u\mA
_{\tau
}+o\bigl(|z|+|y|^{2}\bigr) \bigr\}
\\
&&\qquad\quad\times\exp \biggl\{ -\lambda u_{\tau}\mA_{\tau}\nonumber\\
&&\hspace*{64pt}{}-
\frac
{1}{2} \biggl[ \frac{%
|\mu_{20}\mu_{22}^{-1}z|^{2}}{1-\mu_{20}\mu_{22}^{-1}\mu_{02}}%
+ \biggl\vert
\mu_{22}^{-1/2}z-\frac{\mu_{22}^{1/2}\mathbf
{1}}{2\sigma
}%
\biggr\vert^{2} \biggr]\nonumber \\
&&\hspace*{205pt}{}-\frac{1-\lambda}{2}{y}^{\top}y \biggr\} \,d\mA
_{\tau
}\,dy\,dz.\nonumber
\end{eqnarray}
Thanks to the Borel--TIS inequality (Lemma \ref{LemBorel}), Lemma \ref%
{LemRemainder} and the definition of~$\kappa$ in (\ref{kappa}), for
$x>0$, $%
\gamma_{u}(x)$ is bounded and as $b\rightarrow\infty$,
\[
E \biggl[ \frac{1}{\kappa^{2}};x>\bigl(1+o\bigl(u^{-1-\delta_0/4}\bigr)\bigr)
\bigl[\xi_u + o\bigl(u^{-\delta_0/8}\bigr)\bigr] \biggr]
\rightarrow (2\pi)^{-d}.
\]
Thus, by the dominated convergence theorem and with $H_\lambda$ defined
as in \eqref{HL}, as $u\rightarrow\infty$, we have that
\[
(\ref{cc})\sim\frac{(2\pi)^{-d}}{(1-\rho_{1}-\rho_{2})(1-\rho
_2)}\frac{e^{-\lambda\eta
}\lambda}{2-\lambda}.
\]
We insert it back to (\ref{111}) and obtain that%
%
\begin{eqnarray}
\label{E11} &&E^{Q} \biggl[ \frac{1}{(1-\rho_{1}-\rho_{2})^{2}K^{2}};
\mathcal{E}%
_{b},\mathcal{L},\mA_{\tau}\geq0 \Big\vert
\imath=0, \tau \biggr]
\nonumber
\\
&&\qquad\leq\bigl(1+o(1)\bigr)\frac{(2\pi)^{-d}}{(1-\rho_{1}-\rho_{2})(1-\rho
_2)}\frac{e^{-\lambda\eta}\lambda}{%
2-\lambda}
\nonumber
\\[-8pt]
\\[-8pt]
\nonumber
&&\qquad\quad{}\times \biggl[ \frac{H_{\lambda}^{-1}\det(\Gamma)^{-
{1}/{2}}\det
(\Delta
\mu_{\sigma}(t_{\ast}))^{-1/2}}{(2\pi)^{{(d+1)(d+2)}/{4}-{d}/{2}}%
}\\
&&\hspace*{49pt}{}\times {(1-\lambda)^{d/2}u^{d/2-1}}e^{-({1}/{2})u_{t_*}^{2}+B_{t_{\ast
}}+{\mathbf{1}^{\top}\mu_{22}\mathbf{1}}/{(8\sigma
^{2})}}
\biggr] ^{2}.
\nonumber
\end{eqnarray}
%
Using the asymptotic approximation of $v(b)$ given by Proposition \ref
{CorGRF}, we obtain that
%
\begin{eqnarray}
\label{11st} &&E^{Q} \biggl[ \frac{1}{ [ (1-\rho_{1}-\rho_{2})K ] ^{2}};%
\mathcal{E}_{b},\mathcal{L},\mA_{\tau}\geq0,\imath=0 \biggr]
\nonumber
\\[-8pt]
\\[-8pt]
\nonumber
&&\qquad\leq \frac
{1+o(1)}{1-\rho_{1}-\rho_{2}}\frac{e^{\lambda\eta
}}{\lambda(2-\lambda)}v^{2}(b).
\end{eqnarray}
We choose
\[
\rho_1=\rho_2=\eta= 1- \lambda= 1/\log\log b \sim1/\log
u.
\]
Then, the right-hand side of the above inequality is bounded by
$(1+\eps
)v^2(b)$ for $b$ sufficiently large.

The handling of the second term of (\ref{1st}) is similar except that $
(w,y,z)$ follows density $h_{\tau}^{\ast}(w,y,z)$. Thus, we only
mention the
key steps. Note that
%
\begin{eqnarray}
\label{E12} &&\hspace*{-2pt}E^{Q} \biggl[  \frac{1}{(1-\rho_{1}-\rho
_{2})^{2}K^{2}};\mathcal
{E}%
_{b},\mathcal{L},\mA_{\tau}\geq0 \Big\rrvert
\imath=1,\tau \biggr]
\nonumber
\\
&&\hspace*{-2pt}\qquad=\bigl(1+o(1)\bigr) \biggl[ \frac{H_{\lambda}^{-1}\det(\Gamma)
^{-{1}/{2}}\det
(-\Delta\mu_{\sigma}(t_{\ast}))^{-1/2}}{(2\pi)^{
{(d+1)(d+2)}/{4}-%
{d}/{2}}}\nonumber\\
&&\qquad\hspace*{59pt}{}\times{(1-\lambda)^{d/2}u^{d/2-1}}e^{-
({1}/{2})u_{t_*}^{2}+B_{t_{\ast}}+{\mathbf{1}^{\top}\mu
_{22}\mathbf
{1}}/{%
(8\sigma^{2})}}
\biggr] ^{2}
\\
&&\hspace*{-2pt}\qquad\quad{}\times\frac{\det(\Gamma)^{-{1}/{2}}}{(2\pi)^{
{(d+1)(d+2)}/{4}}}
\nonumber
\\
&&\hspace*{-2pt}\qquad\quad{}\times \int_{\mA_{\tau}\geq0,\mathcal{L}} \gamma_{u}(u
\sigma\mA_{\tau})\nonumber \\
&&\hspace*{78pt}{}\times\exp\biggl\{-2(1-\lambda)u\mathcal A_{\tau}\nonumber\\
&&\hspace*{111pt}{}-\frac{1+o(1)}{2}
\nonumber\\
&&\hspace*{122pt}{}\times\biggl[
{y}^{\top}y+\frac{|w-\mu_{20}\mu_{22}^{-1}z|^{2}}{1-\mu_{20}\mu
_{22}^{-1}\mu_{02}}+z^{\top}\mu_{22}^{-1}z \biggr] \biggr\}\,d\mA_{\tau}\,dy\,dz
\nonumber
\\
&&\hspace*{-2pt}\qquad=O(1) (1-\lambda)^{-1}u^{-1}\cdot{u^{d-2}}e^{-u_{t_*}^{2}}.\nonumber
\end{eqnarray}
%
According to the asymptotic form of $v(b)$ and with $\rho_2 =
1-\lambda
= 1/\log\log b$, we have that
%
\begin{eqnarray}
\label{22nd}&& E^{Q} \biggl[ \frac{1}{ [ (1-\rho_{1}-\rho_{2})K ]
^{2}};\mathcal{%
E}_{b},\mathcal{L},\mA_{\tau}\geq0,\imath=1 \biggr]
\nonumber
\\[-8pt]
\\[-8pt]
\nonumber
&&\qquad=O(1)
\rho_2(1-\lambda)^{-1}u^{-1}\cdot
{u^{d-2}}e^{-u_{t_*}^{2}} =o(1)v^{2}(b).
\end{eqnarray}
Therefore, combining the results in \eqref{11st} and \eqref{22nd}, we
have the first term in~\eqref{part3cal} is bounded by $(1+2\eps)v^2(b)$.

The last step is to show that the second term of (\ref{part3cal}) is
of a
smaller order of $v^{2}(b)$. First, we split the expectation%
%
\begin{eqnarray}
\label{E2}\label{222}
&&E^{Q} \biggl[ \frac{1}{ [ (1-\rho_{1}-\rho_{2})K+\rho
_{1}K_{1} ]
^{2}};
\mathcal{E}_{b},\mathcal{L},\mA_{\tau}<0 \biggr]
\nonumber\\
&&\qquad=E^{Q} \biggl[ \frac{1}{ [ (1-\rho_{1}-\rho_{2})K+\rho
_{1}K_{1} ]
^{2}};\mathcal{E}_{b},
\mathcal{L},\mA_{\tau}<0,\imath=1 \biggr]
\nonumber
\\[-8pt]
\\[-8pt]
\nonumber
&&\qquad\quad{}+E^{Q} \biggl[ \frac{1}{ [ (1-\rho_{1}-\rho_{2})K+\rho
_{1}K_{1} ]
^{2}};\mathcal{E}_{b},
\mathcal{L},-\eta/u_{\tau}<\mA_{\tau
}<0,\imath =0 \biggr]
\nonumber
\\
&&\qquad\quad{}+E^{Q} \biggl[ \frac{1}{ [ (1-\rho_{1}-\rho_{2})K+\rho
_{1}K_{1} ]
^{2}};\mathcal{E}_{b},
\mathcal{L},\mA_{\tau}\leq-\eta/u_{\tau
},\imath=0 \biggr].
\nonumber
\end{eqnarray}
We study these three terms one by one. Let%
%
\begin{eqnarray}
\label{gammaind} \gamma_{1,u}(x)&=& E \biggl[ \frac{1}{ [ \mathbb
{I}_{C_{1}}\cdot
\rho
_{1}\kappa_{1,2}+\mathbb{I}_{C_{2}}\cdot(1-\rho_{1}-\rho
_{2})(1-\lambda)^{-d/2}\kappa_{2,1}%
] ^{2}};
\nonumber
\\[-8pt]
\\[-8pt]
\nonumber
&&\hspace*{22pt} x>\bigl(1+o\bigl(u^{-1-\delta_0/4}\bigr)\bigr)\bigl[\xi_u + o
\bigl(u^{-\delta_0/8}\bigr)\bigr] \Big\vert \imath,\tau,w,{y},{z} \biggr]
\nonumber
.
\end{eqnarray}
We start with the first expectation in \eqref{E2}. Plugging in the
lower bound for
$(1-\rho_{1}-\rho_{2})K+\rho_{1}K_{1}$ derived in \eqref{3star}, we have
%
\begin{eqnarray}
\label{add1} &&E^{Q} \biggl[ \frac{1}{ [ (1-\rho_{1}-\rho_{2})K+\rho
_{1}K_{1} ]
^{2}};
\mathcal{E}_{b},\mathcal{L},\mA_{\tau}<0 \Big| \imath =1,\tau
\biggr]
\nonumber\\
&&\qquad=O(1){u^{d-2}}e^{-u_{t_*}^{2}}\nonumber\\
&&\qquad\quad{}\times \int_{\mA_{\tau}<0}\gamma
_{1,u}(u\sigma\mA_{\tau
})\\
&&\hspace*{73pt}{}\times\exp\biggl\{-2(1+\lambda_{1})u\mA_{\tau}\nonumber\\
&&\hspace*{106pt}{}-\frac{1}{2} \biggl[ {y}^{\top
}y+\frac
{%
|w-\mu_{20}\mu_{22}^{-1}z|^{2}}{1-\mu_{20}\mu_{22}^{-1}\mu_{02}}%
+z^{\top}\mu_{22}^{-1}z \biggr] \biggr\}\,d
\mA_{\tau}\,dy\,dz.
\nonumber
\end{eqnarray}
We deal with the $\gamma_{1,u}(u\sigma\mA_{\tau})$ term in the above
integration. On the set ${\cal L}$, $u\sigma\mA_{\tau
}>-u^{3/2+\epsilon}$. By Lemma \ref{LemGamma}, for $-u^{3/2+\epsilon
}<x<0$, there
exists a constant $\delta^{\ast}>0$ such that
\begin{eqnarray*}
&&E \biggl[ \frac{1}{\rho_{1}^{2}\kappa_{1,2}^{2}};x>\bigl(1+o\bigl(u^{-1-\delta
_0/4}\bigr)\bigr)
\bigl[\xi_u + o\bigl(u^{-\delta_0/8}\bigr)\bigr]\Big \vert\imath, \tau
,w,{y},{z},C_1 \biggr]\\
&&\qquad = O(1)\rho _{1}^{-2}e^{u^{\delta^{\ast}}x}
\end{eqnarray*}
and
\begin{eqnarray*}
&&E \biggl[ \frac{1}{(1-\rho_{1}-\rho_{2})^{2}\kappa_{2,1}^{2}}%
;x>\bigl(1+o\bigl(u^{-1-\delta_0/4}
\bigr)\bigr)\bigl[\xi_u + o\bigl(u^{-\delta_0/8}\bigr)\bigr] \Big\vert
\imath, \tau,w,{y},{z},C_2 \biggr]
\\
&&\qquad = O(1) (1-\rho_{1}-\rho_{2})^{-2}(1-
\lambda)^{-d}e^{u^{\delta
^{\ast}}x}.
\end{eqnarray*}
Therefore,
the above approximations and the dominated
convergence theorem imply that conditionally on ${\cal L}$,
\[
\int_{\mA_{\tau}<0}\gamma_{1,u}(u\sigma
\mA_{\tau
})e^{-2(1+\lambda_{1})u\mA_{\tau}}\,d\mA_{\tau}= O(1)\cdot\max\bigl\{
\rho _1^{-2}, (1-\lambda)^{-2d}\bigr\}\cdot
u^{-1-\delta^*}.
\]
Thus, \eqref{add1} equals
\begin{eqnarray*}
\eqref{add1} &=& O(1) \max\bigl\{\rho_1^{-2}, (1-
\lambda)^{-2d}\bigr\}\cdot u^{-1-\delta
^*}\cdot{u^{d-2}}e^{-u_{t_*}^{2}}.
\end{eqnarray*}
Taking expectation of the above equation with respect to $\imath$ and
$\tau$ and choosing $\rho_1, \rho_2$ and $1-\lambda$ be $(\log\log
b)^{-1}$, we have
%
\begin{equation}
\label{33rd} \quad E \biggl[ \frac{1}{ [ (1-\rho_{1}-\rho_{2})K+\rho
_{1}K_{1}
] ^{2}};%
\mathcal{E}_{b},
\mathcal{L},\mA_{\tau}<0,\imath=1 \biggr] 
=
o(1)v^2(b).
\end{equation}

For the second term in (\ref{222}), with the same bound of $\gamma
_{1,u}$, we have
\begin{eqnarray*}
&&E^{Q} \biggl[ \frac{1}{ [ (1-\rho_{1}-\rho_{2})K+\rho
_{1}K_{1}%
] ^{2}};\mathcal{E}_{b},
\mathcal{L},-\eta/u_{\tau}<\mA _{\tau}<0 \Big| \imath=0,\tau \biggr]
\\
&&\qquad=O(1){u^{d-2}}e^{-u_{t_*}^{2}}\\
&&\qquad\quad{}\times{%
u_{\tau}}
\int_{-{\eta}/{u_{\tau}}<\mA_{\tau}<0}\gamma _{1,u}(u\sigma
\mA_{\tau})e^{-2(1+\lambda_{1})u\mA_{\tau}}e^{-\lambda u_{\tau
}\mA
_{\tau
}}
\\
&&\qquad\quad\times\exp \biggl\{ -\frac{1}{2} \biggl[ \frac{%
|\mu_{20}\mu_{22}^{-1}z|^{2}}{1-\mu_{20}\mu_{22}^{-1}\mu_{02}}%
+ \biggl\vert\mu_{22}^{-1/2}z-\frac{\mu_{22}^{1/2}\mathbf
{1}}{2\sigma}%
\biggr\vert^{2} \biggr] \\
&&\hspace*{201pt}{}-\frac{1-\lambda}{2}{y}^{\top}y \biggr\} \,d\mA
_{\tau}\,dy\,dz
\\
&&\qquad=O(1)\cdot\max\bigl\{\rho_1^{-2}, (1-
\lambda)^{-2d}\bigr\}\cdot u^{-\delta
^*}\cdot{u^{d-2}}e^{-u_{t_*}^{2}}
\\
&&\qquad= o(1)v^2(b),
\end{eqnarray*}
and similarly for the third term in (\ref{222}),
%
\begin{eqnarray}
\label{44th} &&E^{Q} \biggl[ \frac{1}{ [ (1-\rho_{1}-\rho_{2})K+\rho
_{1}K_{1}%
] ^{2}};
\mathcal{E}_{b},\mathcal{L},\mA_{\tau}\leq-\eta
/u_{\tau} \Big| \imath=0,\tau \biggr]
\nonumber
\\
&&\qquad=O(1)\rho_1\cdot{u^{d-2}}e^{-u_{t_*}^{2}}u_\tau
\nonumber\\
&&\qquad\quad{}\times\int_{\mA_{\tau
}<-{\eta}/{u_{\tau}}}\gamma_{1,u}(u\sigma\mA_{\tau
})
\nonumber
\\
&&\qquad\quad{}\times e^{-2(1+\lambda_{1})u\mA_{\tau}}\nonumber\\
&&\qquad\quad{}\times\exp\biggl\{\lambda_{1}u_{\tau}\mA
_{\tau}\\
&&\hspace*{32pt}\qquad\quad{}-%
\frac{1}{2} \biggl[ \frac{|\mu_{20}\mu_{22}^{-1}z|^{2}}{1-\mu
_{20}\mu
_{22}^{-1}\mu_{02}}+ \biggl\vert\mu_{22}^{-1/2}z-\frac{\mu_{22}^{1/2}
\mathbf{1}}{2\sigma} \biggr\vert^{2} \biggr] \nonumber\\
&&\hspace*{200pt}{}-\frac{1+\lambda
_{1}}{2}{y}%
^{\top}y\biggr\}\,d\mA_{\tau}\,dy\,dz
\nonumber
\\
&&\qquad=O(1)\rho_1\cdot\max\bigl\{\rho_1^{-2}, (1-
\lambda)^{-2d}\bigr\}\cdot u^{-\delta^*}\cdot {u^{d-2}}e^{-u_{t_*}^{2}}
\nonumber
\\
&&\qquad= o(1)v^2(b).\nonumber
\end{eqnarray}

We put all the estimates in \eqref{11st}, \eqref{22nd}, \eqref{33rd}
and \eqref{44th} back to \eqref{part3cal}.
For any $\varepsilon>0$, if we choose $\eta= \rho_1 = \rho_2 = 1-
\lambda= 1/\log\log b $, then for $b$ sufficiently large we have that
\[
E^{Q} \biggl[ \biggl( \frac{dP}{dQ} \biggr)^{2};
\mathcal {E}_{b},\mathcal {L} \biggr] \leq(1+3\varepsilon)v^{2}(b).
\]
We complete the proof of Theorem \ref{main} for the case that $\mu(t)
\neq0$.

\subsection{Case 2: Constant mean function}
The proof when $\mu(t)\equiv0$ is very similar, except that we need to
consider two situations: first,
$\protect\tau$ is not close to the boundary of $T$ and otherwise.
More precisely, for a given $\delta'>0$ small enough, we consider the
case when $\tau\in\{t\dvtx |t-\tau|\leq u^{-1/2+\delta^{\prime}}\}
\subset T$ and
otherwise.

For the first situation, $\tau$ is ``far away'' from the
boundary of $T$, which is the important case, the derivation is same as
that of the case where $\mu(t)$ is not a constant.
For the case in which $\tau$ is within $u^{-1/2 + \delta'}$ distance
from the boundary of~$T$, the contribution of the boundary case is
$o(v^2(b))$. An intuitive interpretation is that the important region
of the integral $\mathcal{I}(T)$ might be cut off by the boundary of~$T$.
Therefore, in cases that $\tau$ is too close to the boundary, the
tail $\mathcal{I}(T)$ is not heavier than that of the interior case.
The rigorous analysis is basically repeating the parts~1, 2 and 3 on a
truncated region. Therefore, we omit the details.


\section{Proof of Theorem \texorpdfstring{\protect\ref{sup}}{7}}\label{proofsup}

The proof of Theorem \ref{sup} is analogous to that of Theorem \ref
{main}. According to Lemma \ref{LXsup}, we focus on the set (for some
small $\epsilon_0>0$)
%
\begin{equation}
\label{Lstar} \mL_* = \mL\cap\Bigl\{\sup_{|t-\tau|> 2u^{-1/2+\epsilon}} g(t) -
\epsilon _0 u |t|^2<0\Bigr\}.
\end{equation}
A similar three-part procedure is applied here.

In part 1, using the transformation from $f$ to the process $f_*$, we have
\begin{eqnarray*}
\beta_u(T)&=&\sup_{t\in T} \biggl\{f(t)+
\frac{\mathbf1^\top\bar
f''_t}{2\sigma u_t}+\frac{B_{t}}{u_t}+\mu_\sigma(t) \biggr\}
\\
&=&\sup_{t\in T} \biggl\{f_*(t)+u_\tau C(t-\tau)+
\frac{\mathbf1^\top
(z_t-u_t\mu_{02}+u_\tau\mu_2(t-\tau))}{2\sigma u_t}\\
&&\hspace*{193pt}{}+\frac
{B_{t}}{u_t}+\mu_\sigma(t) \biggr\}.
\end{eqnarray*}
We insert the expansions in (\ref{expf}), (\ref{expc}) and \eqref{expmu}
into the expression of $\beta_u(T)$ and obtain that $\beta_u(T) $ equals
\begin{eqnarray*}
&& \sup_{t\in T} \biggl\{ w+y^{\top}(t-\tau)+
\frac{1}{2}(t-\tau )^{\top}z(t-\tau) +R_{f}(t-\tau)+g(t-
\tau)
\\
&&\hspace*{21pt}{}+u_\tau \biggl(1-\frac{1}{2}(t-\tau )^{\top}(t-
\tau)+C_{4}(t-\tau)+R_{C}(t-\tau) \biggr)
\\
&&\hspace*{21pt}{}+\mu_\sigma(\tau)+\partial\mu_\sigma(\tau)^{\top}(t-
\tau )+\frac
{1}{2}(t-\tau)^{\top}\Delta\mu_\sigma(\tau)
(t-\tau)+\sigma ^{-1}R_{\mu}(t-\tau)
\\
&&\hspace*{171pt}{}+\frac{\mathbf1^\top(z_t-u_t\mu_{02}+u_\tau\mu_2(t-\tau
))}{2\sigma
u_t}+\frac{B_{t}}{u_t} \biggr\}
\\
&&\qquad= \sup_{t\in T} \biggl\{u+w+\frac{1 }{2}{\tilde
y}^{\top
}(uI-{\tilde \zz})^{-1}{\tilde y}\\
&&\hspace*{54pt}{} -\frac{1}{2}
\bigl(t-\tau-(uI-{\tilde\zz})^{-1}{\tilde y}\bigr)^{\top
}(uI-{
\tilde\zz}) \bigl(t-\tau-(uI-{\tilde\zz})^{-1}{\tilde y} \bigr)
\\
&&\hspace*{54pt}{}+ u_\tau C_{4} (t-\tau ) + R (t-\tau )+g(t-\tau )\\
&&\hspace*{140pt}{}+
\frac{\mathbf1^\top(z_t-u_t\mu_{02}+u_\tau\mu_2(t-\tau
))}{2\sigma
u_t}+\frac{B_{t}}{u_t} \biggr\}.
\end{eqnarray*}
Note that the above display is approximately a quadratic function of
$t-\tau$ and is maximized approximately at $t-\tau= (uI -\tilde{\zz
})^{-1} \tilde y$. In addition, on the set $\mL_*$, we have that
$|\tau
- t_*| < 2u^{-1/2+\epsilon}$ and thus $\tilde y = y
+O(u^{-1/2+\epsilon})$.
Therefore, on the set~$\mL_*$, we have the following approximation of
$\beta_u(T)$:
\begin{eqnarray*}
\mA_\tau+\inf_{|t-\tau|<2u^{-1/2+\epsilon}} g(t)& \leq&\beta_u(T)
- u +u^{-1-\delta_0/4}o\bigl(|w|+|y|+|z|\bigr) \\
&\leq&\mA_\tau+\sup
_{|t-\tau
|<2u^{-1/2+\epsilon}} g(t).
\end{eqnarray*}
%
Thus, we obtain the same representation as in part 1 in the proof of
Theorem~\ref{main}.

Since we use the same change of measure, the analysis of the likelihood
ratio is exactly the same as part 2 of Theorem \ref{main}. For part 3,
we compute the second moment of $dP/dQ$ on the set $\{\beta_u(T) > u\}$.
This is also identical to the proof of Theorem \ref{main}.
Thus, with the same choice of tuning parameters, we have that
\[
E^Q \biggl[ \biggl(\frac{dP}{dQ} \biggr)^2;
\beta_u(T) > u \biggr]\leq(1+ \eps)v^2(b).
\]
Additionally, Lemma \ref{LXsup} provides an approximation that
$P(\beta
_u(T)>u)\sim v(b)$. Thus, we use Lemma \ref{LemTV1} (presented at the
beginning of Section~\ref{SecProof}) and complete the proof.

\section{Proof of Theorem \texorpdfstring{\protect\ref{bias}}{10}}\label{proofbias}
For the bias control, we need the following result~\cite{LiuXudensity}.

\begin{proposition}{\label{theorem 4}}
Suppose that conditions \textup{C1--C6} are satisfied. Let $F'(x)$ be the
probability density function of $\log\mathcal{I}(T)=\log\int_T
e^{\sigma f(t)+\mu(t)}\,dt$. Then the following approximation holds as
$x\rightarrow\infty$:
\[
F'(x)\sim\sigma^{-2}x \cdot v\bigl(e^x
\bigr).
\]
\end{proposition}

Thus, for any small $\varepsilon$,
%
\begin{equation}
P \bigl( b<\mathcal{I}(T)<b(1+\varepsilon/\log b) | \mathcal {I}(T)>b \bigr) =
\Theta(\varepsilon). \label{den}
\end{equation}
Similar to the log-normal distribution, the overshoot of $\mathcal
I(T)$ is $\Theta(b/\log b)$. Note that
\[
\bigl|v_M(b)-v(b)\bigr|\leq P\bigl(\mathcal{I}(T)>b,\mathcal{I}_{M}(T)<b
\bigr)+P\bigl(\mathcal{I} 
(T)<b,\mathcal{I}_{M}(T)>b\bigr).
\]
Let
\[
\mL_{\varepsilon}=\Bigl\{\sup_{t\in T}\bigl|\partial f(t)\bigr|\leq2
\bigl(1-u^{-2}\log \varepsilon\bigr)u
\Bigr\}.
\]
Note that $\partial f(t)$ is a $d$-dimensional Gaussian process. Using
Borel--TIS lemma, we obtain that
\[
P \bigl( \mL_{\varepsilon}^{c} \bigr) =o(1) \varepsilon\cdot v(b).
\]
Therefore, it is sufficient to control $P(\mathcal{I}(T)>b,\mathcal
{I}_{M}(T)<b,\mL%
_{\varepsilon})$ and $P(\mathcal{I}(T)<b,\mathcal{I}_{M}(T)>b,\mL%
_{\varepsilon})$.

By the definition of $\mathcal I_M$ in \eqref{defIM}, there exists a
constant $c_{1}>0$ such that
\begin{eqnarray*}
\Delta&=&\bigl|\mathcal{I}(T)-\mathcal{I}_{M}(T)\bigr| \leq\sum
_{i=1}^M \biggl\vert\int_{T_{N}(t_i)}e^{\sigma f(t)+\mu
(t)}\,dt-\operatorname{mes}
\bigl(T_N(t_i)\bigr)\cdot e^{\sigma f(t_{\mathbf{i}})+\mu(t_{\mathbf{i}
})} \biggr\vert
\\
&\leq&c_{1}\min\bigl\{\mathcal{I}_{M}(T),\mathcal{I}(T)
\bigr\}\cdot\sup_{t\in
T}\bigl|\partial f(t)\bigr|/N.
\end{eqnarray*}
Then we have, on the set $\mL_{\varepsilon}$, $\Delta\leq2c_{1}\min
\{%
\mathcal{I}_{M}(T),\mathcal{I}(T)\}(1-u^{-2}\log\varepsilon)u/N$, which
implies that
\begin{eqnarray*}
P\bigl(\mathcal{I}(T)>b,\mathcal{I}_{M}(T)<b,\mL_{\varepsilon}
\bigr)&\leq&  P\bigl(b<\mathcal{I}(T) < b\bigl(1+2\bigl(1-u^{-2}\log
\varepsilon\bigr)u/N\bigr)\bigr)
\\
&=& O(1)\frac{u(1-u^{-2}\log
\varepsilon)\log b}{N}v(b).
\end{eqnarray*}
The last step is due to the result of Proposition \ref{theorem 4} and
further \eqref{den}. Thus, it is sufficient to choose $N=O(\varepsilon
^{-1-\varepsilon_0}u^{2+\varepsilon_{0}})$ so that the above
probability is bounded by $
\eps v(b)$. The bound of $P(\mathcal{I}(T)<b,\mathcal{I}_{M}(T)>b,\mL
_{\varepsilon})$ is completely analogous.

\section{Proof of Theorem \texorpdfstring{\protect\ref{variance}}{11}}\label{proofth5}

The proof of Theorem \ref{variance} is similar to that of Theorem \ref
{main}. Therefore, we only lay out the key steps.
The only difference is that we replace the integral by a finite sum
over $T_{N}$. Recall that the proof of Theorem~\ref{main} consists of
three parts: first, we write the event $\{\mathcal{I}_{M}(T) > b\}$ as
a function of $(w,y,z)$ (with an ignorable correction term); second, we
write the likelihood ratio as a function of $(w,y,z)$ (with an
ignorable correction term); third, we integrate the likelihood ratio
with respect to $(\imath, \tau, w,y,z)$. For the current proof, we also
have three similar parts.

\textit{Part} 1.
For the first step in the proof of Theorem \ref
{main}, we write $\mathcal{I}(T)> b$ if and only if
$\mA_\tau+\frac{o(|w|+|y|+|z|+1)}{u^{1+\delta_0/4}} >u^{-1}\sigma
^{-1}\xi_{u}$. With the current discretization size, as proved in
Theorem \ref{bias},
\[
\log\mathcal{I}(T) - \log\mathcal{I}_{M}(T) = o\bigl(u^{-1-\varepsilon_0/2}
\bigr).
\]
Thus, we reach the same result that $\mathcal{I}_{M}(T)> b$ if $\mA
_\tau
+\frac{o(|w|+|y|+|z|)}{u^{1+\delta_0/4}}+\break   o(u^{-1-\varepsilon_0/2})
>u^{-1}\sigma^{-1}\xi_{u}$.

\textit{Part} 2. Consider the likelihood ratio
\[
\frac{dQ}{dP}= \int_{T} \biggl[(1-
\rho_{1}-\rho_{2}) l(t) \operatorname{LR}(t)+\rho _{1}l(t)\operatorname{LR}_{1}(t)
+ \frac{\rho_{2}}{\operatorname{mes}(T)}\operatorname{LR}_{2}(t) \biggr] \,dt.
\]
Under the discretization setup, we have
\[
\label{dis} \frac{dQ_M}{dP}=\frac{1-{\rho_1}-{\rho_2}}{\kappa} \sum
_{i=1}^M l(t_i) \operatorname{LR}(t_i)+
\frac{\rho_1}{\kappa} \sum_{i=1}^M
l(t_i) \operatorname{LR}_1(t_i)+ {\rho_2}
\sum_{i=1}^M \frac{1}{M}
\operatorname{LR}_2(t_i),
\]
which is a discrete approximation of $dQ/dP$.
In the proof of Theorem \ref{main}, after taking all the terms not
consisting of $t$ out of the integral [such as that in \eqref{LRD}],
the discrete sum is essentially approximating the following integral:
\[
\int_{|t-\tau|<u^{-1+\delta'}} e^{-({(1-\lambda)(u_\tau+\zeta_{u})}/{2}) (t-\tau
-({u}I-{\tilde
\zz})^{-1}{{\tilde y}} )^{\top} ({u}I-{\tilde\zz} )
(t-\tau-({u}I-{\tilde\zz})^{-1}{{\tilde y}} )} \,dt.
\]
The above integral concentrates on a region of size $O(u^{-1})$. Given
that we choose $N> u^{2}$, the discretized likelihood ratio in
$dQ_M/dP$ approximate $dQ/dP$ up to a constant in the sense that
%
\begin{equation}
\label{discretize} \frac{dQ_M}{dP}= \Theta(1) \frac{dQ}{dP}. 
\end{equation}

\textit{Part} 3. With the results of parts~1 and 2, the analysis of
part 3 is completely analogous to part 3 in the proof of Theorem \ref
{main}. Thus, we conclude that
\[
E^{Q_{M}} \bigl(\tilde L^{2}_{b}\bigr)\leq
\kappa_1 v(b)^{2},
\]
where the constant $\kappa_1$ depends on the $\Theta(1)$ in \eqref
{discretize}.

\begin{appendix}
\section*{Appendix: The lemmas}\label{SecLem}

In this section, we state all the lemmas used in
the previous sections. To facilitate reading, we move several lengthy
proofs (Lemmas \ref{LX1}, \ref{LXsup}, \ref{LemExp}, \ref{LemI2},
\ref
{LemA} and \ref{LemGamma}) to the supplemental materials \cite{supp}, as those
proofs are not particularly related to the proof of the theorems and
mostly involve tedious elementary algebra.

The first lemma is known as the Borel--TIS lemma, which was proved
independently by \cite{Bor75,CIS}.

\begin{lemma}[(Borel--TIS)]\label{LemBorel}
Let $f(t)$, $t\in{\cal U}$, ${\cal U}$ is a parameter set,
be a mean zero Gaussian random field. $f$ is almost surely
bounded on ${\cal U}$. Then
$E(\sup_{{\cal U}}f(t) )<\infty$, and
\[
P \Bigl(\max_{t\in{\cal U}}f(t) -E \Bigl[\max_{t\in{\cal U}}f(t)
\Bigr]\geq b \Bigr)\leq e^{
-{b^{2}}/{(2\sigma_{{\cal U}}^{2})}},
\]
where
$\sigma_{{\cal U}}^{2}=\max_{t\in{\cal U}}\operatorname{Var}[f( t)]$.
\end{lemma}

\begin{lemma} \label{LX1}
Conditionally on the set ${\cal L}$ as defined in \eqref{LQ}, we have that
\[
E^{Q} \biggl[ \biggl(\frac{dP}{dQ} \biggr)^2;
\mathcal{I}(T)>b,{\cal L}^c \biggr]= o(1)v^2(b).
\]
\end{lemma}

\begin{lemma} \label{LXsup}
On the set ${\cal L}_*$ as defined in \eqref{Lstar}, we have that for
$k=1$ and $2$
\[
E^{Q} \biggl[ \biggl(\frac{dP}{dQ} \biggr)^k;
\beta_u(T)>u,\mL _*^c \biggr] = o(1)P \bigl(
\beta_u(T)>u \bigr)^k. 
\]
In addition, we have the approximation $P (\beta_u(T)>u
)\sim v(b)$.
\end{lemma}

\begin{lemma}
\label{LemRemainder} Let $\xi_{u}$ be as defined in (\ref{xi}),
then there exist small constants $\delta^*,\lambda',\lambda'' >0$ such
that for all $x>0$ and sufficiently large $u$
\[
P \bigl( |\xi_{u}|>x \bigr) \leq e^{- \lambda' u^{\delta^*}x^{2}}+ e^{-\lambda'' u^2}.
\]
\end{lemma}

\begin{pf}
For $\delta<\delta_0/10$, we split the expectation into two parts $\{
|\tilde S|\leq u^{\delta}\}$ and $\{|\tilde S|> u^\delta, \tau+
(uI-\zz
)^{-1/2}\tilde S \in T\}$.
Note that $|S|\leq\kappa u^{\delta}$ and $g(t)$ is a mean zero
Gaussian random field with $\operatorname{Var}(g(t))=O(|t|^{4+\delta_0})$. A direct
application
of the Borel--TIS inequality (Lemma \ref{LemBorel}) yields the result of
this lemma.
\end{pf}

\begin{lemma}
\label{LemExp} Let $S$ be a random variable taking values in
$ \{s\dvtx (uI-\zz)^{-1/2}s+\tau\in T \}$ with density
proportional to
\[
e^{ -({\sigma}/{2} )(s-(uI-\zz)^{-1/2}\tilde y )^{\top
} (s-(uI-\zz)^{-1/2}\tilde y )}.
\]
If $|y|\leq u^{1/2+\epsilon}$ and $|z|\leq u^{1/2+\epsilon}$ and
$\epsilon\ll\delta_0$, then
\begin{eqnarray*}
&&\log \bigl\{ Ee ^{\sigma
(u-\mu_{\sigma}(\tau) )C_{4} ((uI-\zz)^{-({1}/{2})}S )+\sigma R ((uI-\zz)^{-%
({1}/{2})}S )} \bigr\}
\\
&&\qquad= \frac{1}{8\sigma u}\sum_{i}\partial
_{\mathit{iiii}}^{4}C(0)+\frac{o (|w|+|y|+|z|+1 )}{u^{1+\delta_0/4}},
\end{eqnarray*}
where the expectation is taken with respect to $S$.
\end{lemma}


\begin{lemma}\label{LemDet}
\[
\log\bigl(\det\bigl(I-u^{-1}\zz\bigr)\bigr) = - u^{-1}\operatorname{Tr}(
\zz) + \tfrac{1} 2 u^{-2} \mathbf{I}_2(\zz) + o
\bigl(u^{-2}\bigr),
\]
where $\operatorname{Tr}$ is the trace
of a matrix, $\mathbf{I}_2 (\zz)=\sum_{i=1}^d \lambda_i^2$, and
$\lambda_i$'s are the eigenvalues of~$\zz$.
\end{lemma}
\begin{pf}
The result is immediate by
noting that
$
\det(I-u^{-1}\zz)=\prod_{i=1}^{d}(1-\lambda_{i}/u)$,
and $\operatorname{Tr}(\zz)= \sum_{i=1}^d \lambda_i$.
\end{pf}

\begin{lemma}\label{LemI2}
On the set $\mathcal L$,
$\mathcal{I}_2$ defined as in \eqref{I22} can be written as
\begin{eqnarray*}
&&\int_{A^{*},|t-\tau|<u^{-1+\delta'}} \exp \biggl\{\frac
{u_{t_*}(t-t_*)^{\top}\Delta\mu_\sigma(t_*)(t-t_*)}{2}+
\frac
{u^{2}_{t}}{2} \biggr\}\times u_{t}
\\
&&\qquad{}\times \exp\biggl\{(1-\lambda) u_t  \bigl[w_{t}+ u_\tau C(t-\tau)
-u_t \bigr]\\
&&\hspace*{56pt}{}+
\frac
{(1-\lambda)}{2\sigma}
\mathbf1^\top \bigl(z_t-\mu_{02}u_t+\mu_2(t-\tau)u_\tau \bigr)
-\lambda
B_{t}- \frac{\mathbf1^\top\mu_{22}\mathbf1}{8\sigma^2}\biggr\}
\\
&&\qquad{} \times \exp\bigl\{\bigl(\bigl (w_{t}+u_\tau C(t-\tau)-u_t \bigr)^2
-2\bigl (w_{t}+u_\tau C(t-\tau)-u_t \bigr)\\
&&\hspace*{122pt}{}\times\mu_{20}\mu_{22}^{-1}{
\bigl(z_t-\mu
_{02}u_t+\mu_2(t-\tau)u_\tau \bigr)}\bigr)\\
&&\hspace*{200pt}{}/\bigl(2 (1-\mu_{20}\mu
_{22}^{-1}\mu
_{02} )\bigr)\bigr\}\,dt.
\end{eqnarray*}
\end{lemma}


\begin{lemma}\label{LemA}
For $\eta=1/\log\log b$, on the set $\mL$, if $\mA_{\tau}\geq0$, then
\[
\bigl\{|t-\tau|\leq u^{-1+\delta'}\bigr\}\subseteq A^{*}.
\]
\end{lemma}

\begin{lemma}
\label{LemGamma} On the set $\mL$,
there exists some $\delta^{\ast}>0$ such that for all
$-u^{3/2+\epsilon}<x <0$,
\begin{eqnarray*}
&&E \biggl[ \frac{1}{\rho_1^2\kappa_{1,2}^2};x>\bigl(1+o\bigl(u^{-1-\delta
_0/4}\bigr)\bigr)
\bigl[\xi _u + o\bigl(u^{-\delta_0/8}\bigr)\bigr] \Big\vert \imath,
\tau,w,{y},{z}, {C_1} \biggr] \\
&&\qquad= O(1)\rho_1^{-2}e^{u^{\delta^{\ast}}x},
\\
&& E \biggl[\frac{1}{(1-\rho_1-\rho_2)^2\kappa
_{2,1}^2};x>\bigl(1+o\bigl(u^{-1-\delta_0/4}\bigr)\bigr)
\bigl[\xi_u + o\bigl(u^{-\delta_0/8}\bigr)\bigr] \Big\vert \imath,
\tau,w,{y},{z},{C_2} \biggr]
\\
&&\qquad= O(1) (1-\rho_1-\rho_2)^{-2}(1-
\lambda)^{-d}e^{u^{\delta^{\ast}}x},
\end{eqnarray*}
where ${C_1}=\{\operatorname{mes}(A^c\cap D)\geq \operatorname{mes}(A\cap D)\}$ and
${C_2}=\{\operatorname{mes}(A^c\cap D)< \operatorname{mes}(A\cap D)\}$.
\end{lemma}
\end{appendix}
%

\begin{supplement}[id=suppA]
\stitle{Supplement to ``On the conditional
distributions and the efficient simulations of
exponential integrals of gaussian random fields''\\}
\slink[doi]{10.1214/13-AAP960SUPP} 
\sdatatype{.pdf}
\sfilename{aap960\_supp.pdf}
\sdescription{Proofs of Proposition \ref{CorGRF} and
Lemmas \ref{LX1}, \ref{LXsup}, \ref{LemExp}, \ref{LemI2}, \ref{LemA} and~\ref{LemGamma}
are provided in the supplementary material.}
\end{supplement}

%
%

%


\printaddresses


\begin{thebibliography}{56}

\def\cprime{$'$}

\bibitem{Aberg08}
\begin{barticle}[mr]
\bauthor{\bsnm{Aberg},~\bfnm{Sofia}\binits{S.}} \AND
\bauthor{\bsnm{Guttorp},~\bfnm{Peter}\binits{P.}}
(\byear{2008}).
\btitle{Distribution of the maximum in air pollution fields}.
\bjournal{Environmetrics}
\bvolume{19}
\bpages{183--208}.
\bid{doi={10.1002/env.866}, issn={1180-4009}, mr={2420185}}
\end{barticle}
\bptok{imsref}%
\endbibitem

\bibitem{Adl81}
\begin{bbook}[mr]
\bauthor{\bsnm{Adler},~\bfnm{Robert~J.}\binits{R.~J.}}
(\byear{1981}).
\btitle{The Geometry of Random Fields}.
\bpublisher{Wiley},
\blocation{Chichester}.
\bid{mr={0611857}}
\end{bbook}
\bptok{imsref}%
\endbibitem

\bibitem{ABL09}
\begin{barticle}[mr]
\bauthor{\bsnm{Adler},~\bfnm{Robert~J.}\binits{R.~J.}},
\bauthor{\bsnm{Blanchet},~\bfnm{Jose~H.}\binits{J.~H.}} \AND
\bauthor{\bsnm{Liu},~\bfnm{Jingchen}\binits{J.}}
(\byear{2012}).
\btitle{Efficient {M}onte {C}arlo for high excursions of {G}aussian random fields}.
\bjournal{Ann. Appl. Probab.}
\bvolume{22}
\bpages{1167--1214}.
\bid{doi={10.1214/11-AAP792}, issn={1050-5164}, mr={2977989}}
\end{barticle}
\bptok{imsref}%
\endbibitem

\bibitem{AST09}
\begin{barticle}[mr]
\bauthor{\bsnm{Adler},~\bfnm{Robert~J.}\binits{R.~J.}},
\bauthor{\bsnm{Samorodnitsky},~\bfnm{Gennady}\binits{G.}} \AND
\bauthor{\bsnm{Taylor},~\bfnm{Jonathan~E.}\binits{J.~E.}}
(\byear{2013}).
\btitle{High level excursion set geometry for non-{G}aussian infinitely divisible random fields}.
\bjournal{Ann. Probab.}
\bvolume{41}
\bpages{134--169}.
\bid{doi={10.1214/11-AOP738}, issn={0091-1798}, mr={3059195}}
\bptnote{check year}%
\end{barticle}
\bptok{imsref}%
\endbibitem

\bibitem{AdlTay07}
\begin{bbook}[mr]
\bauthor{\bsnm{Adler},~\bfnm{Robert~J.}\binits{R.~J.}} \AND
\bauthor{\bsnm{Taylor},~\bfnm{Jonathan~E.}\binits{J.~E.}}
(\byear{2007}).
\btitle{Random Fields and Geometry}.
\bpublisher{Springer},
\blocation{New York}.
\bid{mr={2319516}}
\end{bbook}
\bptok{imsref}%
\endbibitem

\bibitem{ATW09}
\begin{bmisc}[auto:STB|2013/12/09|07:59:19]
\bauthor{\bsnm{Adler},~\bfnm{R.~J.}\binits{R.~J.}},
\bauthor{\bsnm{Taylor},~\bfnm{J.~E.}\binits{J.~E.}} \AND
\bauthor{\bsnm{Worsley},~\bfnm{K.~J.}\binits{K.~J.}}
(\byear{2009}).
\bhowpublished{Applications of random fields and geometry:
Foundations and case studies. Preprint}.
\end{bmisc}
\bptok{imsref}%
\endbibitem

\bibitem{Ahs78}
\begin{barticle}[auto:STB|2013/12/09|07:59:19]
\bauthor{\bsnm{Ahsan},~\bfnm{S.~M.}\binits{S.~M.}}
(\byear{1978}).
\btitle{Portfolio selection in a lognormal securities market}.
\bjournal{Zeitschrift f\"{u}r Nationalokonomie---Journal of Economics}
\bvolume{38}
\bpages{105--118}.
\end{barticle}
\bptok{imsref}%
\endbibitem

\bibitem{ASM00}
\begin{bbook}[mr]
\bauthor{\bsnm{Asmussen},~\bfnm{S{\o}ren}\binits{S.}}
(\byear{2000}).
\btitle{Ruin Probabilities}.
\bseries{Advanced Series on Statistical Science \& Applied Probability}
\bvolume{2}.
\bpublisher{World Scientific},
\blocation{River Edge, NJ}.
\bid{doi={10.1142/9789812779311}, mr={1794582}}
\end{bbook}
\bptok{imsref}%
\endbibitem

\bibitem{AsmussenKlupp}
\begin{barticle}[mr]
\bauthor{\bsnm{Asmussen},~\bfnm{S{\o}ren}\binits{S.}} \AND
\bauthor{\bsnm{Kl{\"u}ppelberg},~\bfnm{Claudia}\binits{C.}}
(\byear{1996}).
\btitle{Large deviations results for subexponential tails, with applications to insurance risk}.
\bjournal{Stochastic Process. Appl.}
\bvolume{64}
\bpages{103--125}.
\bid{doi={10.1016/S0304-4149(96)00087-7}, issn={0304-4149}, mr={1419495}}
\end{barticle}
\bptok{imsref}%
\endbibitem

\bibitem{AR08}
\begin{barticle}[mr]
\bauthor{\bsnm{Asmussen},~\bfnm{S{\o}ren}\binits{S.}} \AND
\bauthor{\bsnm{Rojas-Nandayapa},~\bfnm{Leonardo}\binits{L.}}
(\byear{2008}).
\btitle{Asymptotics of sums of lognormal random variables with {G}aussian copula}.
\bjournal{Statist. Probab. Lett.}
\bvolume{78}
\bpages{2709--2714}.
\bid{doi={10.1016/j.spl.2008.03.035}, issn={0167-7152}, mr={2465111}}
\end{barticle}
\bptok{imsref}%
\endbibitem

\bibitem{AW05}
\begin{barticle}[mr]
\bauthor{\bsnm{Aza{\"{\i}}s},~\bfnm{Jean-Marc}\binits{J.-M.}} \AND
\bauthor{\bsnm{Wschebor},~\bfnm{Mario}\binits{M.}}
(\byear{2005}).
\btitle{On the distribution of the maximum of a {G}aussian field with {$d$} parameters}.
\bjournal{Ann. Appl. Probab.}
\bvolume{15}
\bpages{254--278}.
\bid{doi={10.1214/105051604000000602}, issn={1050-5164}, mr={2115043}}
\end{barticle}
\bptok{imsref}%
\endbibitem

\bibitem{AW08}
\begin{barticle}[mr]
\bauthor{\bsnm{Aza{\"{\i}}s},~\bfnm{Jean-Marc}\binits{J.-M.}} \AND
\bauthor{\bsnm{Wschebor},~\bfnm{Mario}\binits{M.}}
(\byear{2008}).
\btitle{A general expression for the distribution of the maximum of a {G}aussian field and the approximation of the tail}.
\bjournal{Stochastic Process. Appl.}
\bvolume{118}
\bpages{1190--1218}.
\bid{doi={10.1016/j.spa.2007.07.016}, issn={0304-4149}, mr={2428714}}
\end{barticle}
\bptok{imsref}%
\endbibitem

\bibitem{AW09}
\begin{bbook}[mr]
\bauthor{\bsnm{Aza{\"{\i}}s},~\bfnm{Jean-Marc}\binits{J.-M.}} \AND
\bauthor{\bsnm{Wschebor},~\bfnm{Mario}\binits{M.}}
(\byear{2009}).
\btitle{Level Sets and Extrema of Random Processes and Fields}.
\bpublisher{Wiley},
\blocation{Hoboken, NJ}.
\bid{doi={10.1002/9780470434642}, mr={2478201}}
\end{bbook}
\bptok{imsref}%
\endbibitem

\bibitem{BasSha01}
\begin{barticle}[auto:STB|2013/12/09|07:59:19]
\bauthor{\bsnm{Basak},~\bfnm{S.}\binits{S.}} \AND
\bauthor{\bsnm{Shapiro},~\bfnm{A.}\binits{A.}}
(\byear{2001}).
\btitle{Value-at-risk-based risk management: Optimal policies and asset prices}.
\bjournal{Review of Financial Studies}
\bvolume{14}
\bpages{371--405}.
\end{barticle}
\bptok{imsref}%
\endbibitem

\bibitem{Berman85}
\begin{barticle}[mr]
\bauthor{\bsnm{Berman},~\bfnm{Simeon~M.}\binits{S.~M.}}
(\byear{1985}).
\btitle{An asymptotic formula for the distribution of the maximum of a {G}aussian process with stationary increments}.
\bjournal{J. Appl. Probab.}
\bvolume{22}
\bpages{454--460}.
\bid{issn={0021-9002}, mr={0789369}}
\end{barticle}
\bptok{imsref}%
\endbibitem

\bibitem{BL10}
\begin{barticle}[mr]
\bauthor{\bsnm{Blanchet},~\bfnm{Jose}\binits{J.}} \AND
\bauthor{\bsnm{Liu},~\bfnm{Jingchen}\binits{J.}}
(\byear{2010}).
\btitle{Efficient importance sampling in ruin problems for multidimensional regularly varying random walks}.
\bjournal{J. Appl. Probab.}
\bvolume{47}
\bpages{301--322}.
\bid{doi={10.1239/jap/1276784893}, issn={0021-9002}, mr={2668490}}
\end{barticle}
\bptok{imsref}%
\endbibitem

\bibitem{BL11}
\begin{barticle}[mr]
\bauthor{\bsnm{Blanchet},~\bfnm{Jose}\binits{J.}} \AND
\bauthor{\bsnm{Liu},~\bfnm{Jingchen}\binits{J.}}
(\byear{2012}).
\btitle{Efficient simulation and conditional functional limit theorems for ruinous heavy-tailed random walks}.
\bjournal{Stochastic Process. Appl.}
\bvolume{122}
\bpages{2994--3031}.
\bid{doi={10.1016/j.spa.2012.05.001}, issn={0304-4149}, mr={2931349}}
\end{barticle}
\bptok{imsref}%
\endbibitem

\bibitem{BL08}
\begin{barticle}[mr]
\bauthor{\bsnm{Blanchet},~\bfnm{Jose~H.}\binits{J.~H.}} \AND
\bauthor{\bsnm{Liu},~\bfnm{Jingchen}\binits{J.}}
(\byear{2008}).
\btitle{State-dependent importance sampling for regularly varying random walks}.
\bjournal{Adv. in Appl. Probab.}
\bvolume{40}
\bpages{1104--1128}.
\bid{issn={0001-8678}, mr={2488534}}
\end{barticle}
\bptok{imsref}%
\endbibitem

\bibitem{BLY}
\begin{bmisc}[auto:STB|2013/12/09|07:59:19]
\bauthor{\bsnm{Blanchet},~\bfnm{J.~H.}\binits{J.~H.}},
\bauthor{\bsnm{Liu},~\bfnm{J.}\binits{J.}} \AND
\bauthor{\bsnm{Yang},~\bfnm{X.}\binits{X.}}
(\byear{2010}).
\bhowpublished{Monte Carlo for large credit portfolios with
potentially high correlations. In \textit{Proceedings of the 2010
Winter Simulation Conference}, Baltimore, MD.}
\end{bmisc}
\bptok{imsref}%
\endbibitem

\bibitem{Bor75}
\begin{barticle}[mr]
\bauthor{\bsnm{Borell},~\bfnm{Christer}\binits{C.}}
(\byear{1975}).
\btitle{The {B}runn--{M}inkowski inequality in {G}auss space}.
\bjournal{Invent. Math.}
\bvolume{30}
\bpages{207--216}.
\bid{issn={0020-9910}, mr={0399402}}
\end{barticle}
\bptok{imsref}%
\endbibitem

\bibitem{Camp94}
\begin{barticle}[auto:STB|2013/12/09|07:59:19]
\bauthor{\bsnm{Campbell},~\bfnm{M.~J.}\binits{M.~J.}}
(\byear{1994}).
\btitle{Time-series regression for counts---An investigation into the relationship between sudden-infant-death-syndrome and environmental-temperature}.
\bjournal{J. Roy. Statist. Soc. Ser. A}
\bvolume{157}
\bpages{191--208}.
\end{barticle}
\bptok{imsref}%
\endbibitem

\bibitem{ChLe95}
\begin{barticle}[mr]
\bauthor{\bsnm{Chan},~\bfnm{K.~S.}\binits{K.~S.}} \AND
\bauthor{\bsnm{Ledolter},~\bfnm{Johannes}\binits{J.}}
(\byear{1995}).
\btitle{Monte {C}arlo {EM} estimation for time series models involving counts}.
\bjournal{J. Amer. Statist. Assoc.}
\bvolume{90}
\bpages{242--252}.
\bid{issn={0162-1459}, mr={1325132}}
\end{barticle}
\bptok{imsref}%
\endbibitem

\bibitem{CIS}
\begin{binproceedings}[mr]
\bauthor{\bsnm{Cirel'son},~\bfnm{B.~S.}\binits{B.~S.}},
\bauthor{\bsnm{Ibragimov},~\bfnm{I.~A.}\binits{I.~A.}} \AND
\bauthor{\bsnm{Sudakov},~\bfnm{V.~N.}\binits{V.~N.}}
(\byear{1976}).
\btitle{Norms of {G}aussian sample functions}.
In \bbooktitle{Proceedings of the {T}hird {J}apan--{USSR} {S}ymposium on
{P}robability Theory ({T}ashkent, 1975)}.
\bseries{Lecture Notes in Math.}
\bvolume{550}
\bpages{20--41}.
\bpublisher{Springer},
\blocation{Berlin}.
\bid{mr={0458556}}
\end{binproceedings}
\bptok{imsref}%
\endbibitem

\bibitem{COX55}
\begin{barticle}[mr]
\bauthor{\bsnm{Cox},~\bfnm{D.~R.}\binits{D.~R.}}
(\byear{1955}).
\btitle{Some statistical methods connected with series of events}.
\bjournal{J. R. Stat. Soc. Ser. B Stat. Methodol.}
\bvolume{17}
\bpages{129--157; discussion, 157--164}.
\bid{issn={0035-9246}, mr={0092301}}
\bptnote{check related}%
\end{barticle}
\bptok{imsref}%
\endbibitem

\bibitem{COIS80}
\begin{bbook}[mr]
\bauthor{\bsnm{Cox},~\bfnm{David~Roxbee}\binits{D.~R.}} \AND
\bauthor{\bsnm{Isham},~\bfnm{Valerie}\binits{V.}}
(\byear{1980}).
\btitle{Point Processes}.
\bpublisher{Chapman \& Hall},
\blocation{London}.
\bid{mr={0598033}}
\end{bbook}
\bptok{imsref}%
\endbibitem

\bibitem{DDW00}
\begin{barticle}[mr]
\bauthor{\bsnm{Davis},~\bfnm{Richard~A.}\binits{R.~A.}},
\bauthor{\bsnm{Dunsmuir},~\bfnm{William~T.~M.}\binits{W.~T.~M.}} \AND
\bauthor{\bsnm{Wang},~\bfnm{Ying}\binits{Y.}}
(\byear{2000}).
\btitle{On autocorrelation in a {P}oisson regression model}.
\bjournal{Biometrika}
\bvolume{87}
\bpages{491--505}.
\bid{doi={10.1093/biomet/87.3.491}, issn={0006-3444}, mr={1789805}}
\end{barticle}
\bptok{imsref}%
\endbibitem

\bibitem{Due04}
\begin{bbook}[auto:STB|2013/12/09|07:59:19]
\bauthor{\bsnm{Deutsch},~\bfnm{H.~P.}\binits{H.~P.}}
(\byear{2004}).
\btitle{Derivatives and Internal Models},
\bedition{3rd}~ed.
\bpublisher{Palgrave Macmillan},
\blocation{Basingstoke, UK}.
\end{bbook}
\bptok{imsref}%
\endbibitem

\bibitem{DufPan97}
\begin{barticle}[auto:STB|2013/12/09|07:59:19]
\bauthor{\bsnm{Duffie},~\bfnm{D.}\binits{D.}} \AND
\bauthor{\bsnm{Pan},~\bfnm{J.}\binits{J.}}
(\byear{1997}).
\btitle{An overview of value at risk}.
\bjournal{The Journal of Derivatives}
\bvolume{4}
\bpages{7--49}.
\end{barticle}
\bptok{imsref}%
\endbibitem

\bibitem{Duf01}
\begin{barticle}[mr]
\bauthor{\bsnm{Dufresne},~\bfnm{Daniel}\binits{D.}}
(\byear{2001}).
\btitle{The integral of geometric {B}rownian motion}.
\bjournal{Adv. in Appl. Probab.}
\bvolume{33}
\bpages{223--241}.
\bid{doi={10.1239/aap/999187905}, issn={0001-8678}, mr={1825324}}
\end{barticle}
\bptok{imsref}%
\endbibitem

\bibitem{DupEll97}
\begin{bbook}[mr]
\bauthor{\bsnm{Dupuis},~\bfnm{Paul}\binits{P.}} \AND
\bauthor{\bsnm{Ellis},~\bfnm{Richard~S.}\binits{R.~S.}}
(\byear{1997}).
\btitle{A Weak Convergence Approach to the Theory of Large Deviations}.
\bpublisher{Wiley},
\blocation{New York}.
\bid{doi={10.1002/9781118165904}, mr={1431744}}
\end{bbook}
\bptok{imsref}%
\endbibitem

\bibitem{FR10}
\begin{barticle}[mr]
\bauthor{\bsnm{Foss},~\bfnm{Serguei}\binits{S.}} \AND
\bauthor{\bsnm{Richards},~\bfnm{Andrew}\binits{A.}}
(\byear{2010}).
\btitle{On sums of conditionally independent subexponential random variables}.
\bjournal{Math. Oper. Res.}
\bvolume{35}
\bpages{102--119}.
\bid{doi={10.1287/moor.1090.0430}, issn={0364-765X}, mr={2676758}}
\end{barticle}
\bptok{imsref}%
\endbibitem

\bibitem{GadrichAdler93}
\begin{barticle}[mr]
\bauthor{\bsnm{Gadrich},~\bfnm{Tamar}\binits{T.}} \AND
\bauthor{\bsnm{Adler},~\bfnm{Robert~J.}\binits{R.~J.}}
(\byear{1993}).
\btitle{Slepian models for nonstationary {G}aussian processes}.
\bjournal{J. Appl. Probab.}
\bvolume{30}
\bpages{98--111}.
\bid{issn={0021-9002}, mr={1206355}}
\end{barticle}
\bptok{imsref}%
\endbibitem

\bibitem{GHS00}
\begin{barticle}[auto:STB|2013/12/09|07:59:19]
\bauthor{\bsnm{Glasserman},~\bfnm{P.}\binits{P.}},
\bauthor{\bsnm{Heidelberger},~\bfnm{P.}\binits{P.}} \AND
\bauthor{\bsnm{Shahabuddin},~\bfnm{P.}\binits{P.}}
(\byear{2000}).
\btitle{Variance reduction techniques for estimating value-at-risk}.
\bjournal{Management Science}
\bvolume{46}
\bpages{1349--1364}.
\end{barticle}
\bptok{imsref}%
\endbibitem

\bibitem{Grigoriu89}
\begin{barticle}[auto:STB|2013/12/09|07:59:19]
\bauthor{\bsnm{Grigoriu},~\bfnm{M.}\binits{M.}}
(\byear{1989}).
\btitle{Reliability of Daniels systems subject to quasistatic and dynamic non-stationary Gaussian load processes}.
\bjournal{Prob. Eng. Mech.}
\bvolume{4}
\bpages{128--134}.
\end{barticle}
\bptok{imsref}%
\endbibitem

\bibitem{Hu90}
\begin{barticle}[mr]
\bauthor{\bsnm{H{\"u}sler},~\bfnm{J.}\binits{J.}}
(\byear{1990}).
\btitle{Extreme values and high boundary crossings of locally stationary {G}aussian processes}.
\bjournal{Ann. Probab.}
\bvolume{18}
\bpages{1141--1158}.
\bid{issn={0091-1798}, mr={1062062}}
\end{barticle}
\bptok{imsref}%
\endbibitem

\bibitem{HPZ11}
\begin{barticle}[mr]
\bauthor{\bsnm{H{\"u}sler},~\bfnm{J{\"u}rg}\binits{J.}},
\bauthor{\bsnm{Piterbarg},~\bfnm{Vladimir}\binits{V.}} \AND
\bauthor{\bsnm{Zhang},~\bfnm{Yueming}\binits{Y.}}
(\byear{2011}).
\btitle{Extremes of {G}aussian processes with random variance}.
\bjournal{Electron. J. Probab.}
\bvolume{16}
\bpages{1254--1280}.
\bid{doi={10.1214/EJP.v16-904}, issn={1083-6489}, mr={2827458}}
\end{barticle}
\bptok{imsref}%
\endbibitem

\bibitem{LS70}
\begin{barticle}[mr]
\bauthor{\bsnm{Landau},~\bfnm{H.~J.}\binits{H.~J.}} \AND
\bauthor{\bsnm{Shepp},~\bfnm{L.~A.}\binits{L.~A.}}
(\byear{1970}).
\btitle{On the supremum of a {G}aussian process}.
\bjournal{Sankhy\=a Ser. A}
\bvolume{32}
\bpages{369--378}.
\bid{issn={0581-572X}, mr={0286167}}
\end{barticle}
\bptok{imsref}%
\endbibitem

\bibitem{Leadbetter83}
\begin{bbook}[mr]
\bauthor{\bsnm{Leadbetter},~\bfnm{M.~R.}\binits{M.~R.}},
\bauthor{\bsnm{Lindgren},~\bfnm{Georg}\binits{G.}} \AND
\bauthor{\bsnm{Rootz{\'e}n},~\bfnm{Holger}\binits{H.}}
(\byear{1983}).
\btitle{Extremes and Related Properties of Random Sequences and Processes}.
\bpublisher{Springer},
\blocation{New York}.
\bid{mr={0691492}}
\end{bbook}
\bptok{imsref}%
\endbibitem

\bibitem{LT91}
\begin{bbook}[mr]
\bauthor{\bsnm{Ledoux},~\bfnm{Michel}\binits{M.}} \AND
\bauthor{\bsnm{Talagrand},~\bfnm{Michel}\binits{M.}}
(\byear{1991}).
\btitle{Probability in {B}anach Spaces: Isoperimetry and Processes}.
\bseries{Ergebnisse der Mathematik und Ihrer Grenzgebiete (3)}
\bvolume{23}.
\bpublisher{Springer},
\blocation{Berlin}.
\bid{mr={1102015}}
\end{bbook}
\bptok{imsref}%
\endbibitem

\bibitem{Lindgren70}
\begin{barticle}[mr]
\bauthor{\bsnm{Lindgren},~\bfnm{Georg}\binits{G.}}
(\byear{1970}).
\btitle{Some properties of a normal process near a local maximum}.
\bjournal{Ann. Math. Statist.}
\bvolume{41}
\bpages{1870--1883}.
\bid{issn={0003-4851}, mr={0272051}}
\end{barticle}
\bptok{imsref}%
\endbibitem

\bibitem{Lindgren79}
\begin{barticle}[mr]
\bauthor{\bsnm{Lindgren},~\bfnm{Georg}\binits{G.}}
(\byear{1979}).
\btitle{Prediction of level crossings for normal processes containing deterministic components}.
\bjournal{Adv. in Appl. Probab.}
\bvolume{11}
\bpages{93--117}.
\bid{doi={10.2307/1426770}, issn={0001-8678}, mr={0517553}}
\end{barticle}
\bptok{imsref}%
\endbibitem

\bibitem{Liu10}
\begin{barticle}[mr]
\bauthor{\bsnm{Liu},~\bfnm{Jingchen}\binits{J.}}
(\byear{2012}).
\btitle{Tail approximations of integrals of {G}aussian random fields}.
\bjournal{Ann. Probab.}
\bvolume{40}
\bpages{1069--1104}.
\bid{doi={10.1214/10-AOP639}, issn={0091-1798}, mr={2962087}}
\end{barticle}
\bptok{imsref}%
\endbibitem

\bibitem{LiuXu11}
\begin{barticle}[mr]
\bauthor{\bsnm{Liu},~\bfnm{Jingchen}\binits{J.}} \AND
\bauthor{\bsnm{Xu},~\bfnm{Gongjun}\binits{G.}}
(\byear{2012}).
\btitle{Some asymptotic results of {G}aussian random fields with varying mean functions and the associated processes}.
\bjournal{Ann. Statist.}
\bvolume{40}
\bpages{262--293}.
\bid{doi={10.1214/11-AOS960}, issn={0090-5364}, mr={3014307}}
\end{barticle}
\bptok{imsref}%
\endbibitem

\bibitem{LiuXudensity}
\begin{barticle}[mr]
\bauthor{\bsnm{Liu},~\bfnm{Jingchen}\binits{J.}} \AND
\bauthor{\bsnm{Xu},~\bfnm{Gongjun}\binits{G.}}
(\byear{2013}).
\btitle{On the density functions of integrals of {G}aussian random fields}.
\bjournal{Adv. in Appl. Probab.}
\bvolume{45}
\bpages{398--424}.
\bid{doi={10.1239/aap/1370870124}, issn={0001-8678}, mr={3102457}}
\end{barticle}
\bptok{imsref}%
\endbibitem


\bibitem{supp}
\begin{bmisc}[auto]
\bauthor{\bsnm{Liu},~\bfnm{Jingchen}\binits{J.}} \AND
\bauthor{\bsnm{Xu},~\bfnm{Gongjun}\binits{G.}}
(\byear{2014}).
\bhowpublished{Supplement to ``On the conditional distributions and the efficient simulations
of exponential integrals of gaussian random fields.'' DOI:\doiurl{10.1214/13-AAP960SUPP}.}
\end{bmisc}
\bptok{imsref}%
\endbibitem


\bibitem{MS70}
\begin{barticle}[mr]
\bauthor{\bsnm{Marcus},~\bfnm{M.~B.}\binits{M.~B.}} \AND
\bauthor{\bsnm{Shepp},~\bfnm{L.~A.}\binits{L.~A.}}
(\byear{1970}).
\btitle{Continuity of {G}aussian processes}.
\bjournal{Trans. Amer. Math. Soc.}
\bvolume{151}
\bpages{377--391}.
\bid{issn={0002-9947}, mr={0264749}}
\end{barticle}
\bptok{imsref}%
\endbibitem

\bibitem{MitzUpf05}
\begin{bbook}[mr]
\bauthor{\bsnm{Mitzenmacher},~\bfnm{Michael}\binits{M.}} \AND
\bauthor{\bsnm{Upfal},~\bfnm{Eli}\binits{E.}}
(\byear{2005}).
\btitle{Probability and Computing: Randomized Algorithms and Probabilistic Analysis}.
\bpublisher{Cambridge Univ. Press},
\blocation{Cambridge}.
\bid{mr={2144605}}
\end{bbook}
\bptok{imsref}%
\endbibitem

\bibitem{NSY08}
\begin{barticle}[mr]
\bauthor{\bsnm{Nardi},~\bfnm{Yuval}\binits{Y.}},
\bauthor{\bsnm{Siegmund},~\bfnm{David~O.}\binits{D.~O.}} \AND
\bauthor{\bsnm{Yakir},~\bfnm{Benjamin}\binits{B.}}
(\byear{2008}).
\btitle{The distribution of maxima of approximately {G}aussian random fields}.
\bjournal{Ann. Statist.}
\bvolume{36}
\bpages{1375--1403}.
\bid{doi={10.1214/07-AOS511}, issn={0090-5364}, mr={2418661}}
\end{barticle}
\bptok{imsref}%
\endbibitem

\bibitem{Pit96}
\begin{bbook}[mr]
\bauthor{\bsnm{Piterbarg},~\bfnm{Vladimir~I.}\binits{V.~I.}}
(\byear{1996}).
\btitle{Asymptotic Methods in the Theory of {G}aussian Processes and Fields}.
\bseries{Translations of Mathematical Monographs}
\bvolume{148}.
\bpublisher{Amer. Math. Soc.},
\blocation{Providence, RI}.
\bid{mr={1361884}}
\end{bbook}
\bptok{imsref}%
\endbibitem

\bibitem{ST74}
\begin{barticle}[auto:STB|2013/12/09|07:59:19]
\bauthor{\bsnm{Sudakov},~\bfnm{V.~N.}\binits{V.~N.}} \AND
\bauthor{\bsnm{Tsirelson},~\bfnm{B.~S.}\binits{B.~S.}}
(\byear{1974}).
\btitle{Extremal properties of half spaces for spherically invariant measures}.
\bjournal{Zap. Nauchn. Sem. LOMI}
\bvolume{45}
\bpages{75--82}.
\end{barticle}
\bptok{imsref}%
\endbibitem

\bibitem{Sun93}
\begin{barticle}[mr]
\bauthor{\bsnm{Sun},~\bfnm{Jiayang}\binits{J.}}
(\byear{1993}).
\btitle{Tail probabilities of the maxima of {G}aussian random fields}.
\bjournal{Ann. Probab.}
\bvolume{21}
\bpages{34--71}.
\bid{issn={0091-1798}, mr={1207215}}
\end{barticle}
\bptok{imsref}%
\endbibitem

\bibitem{TA96}
\begin{barticle}[mr]
\bauthor{\bsnm{Talagrand},~\bfnm{Michel}\binits{M.}}
(\byear{1996}).
\btitle{Majorizing measures: The generic chaining}.
\bjournal{Ann. Probab.}
\bvolume{24}
\bpages{1049--1103}.
\bid{doi={10.1214/aop/1065725175}, issn={0091-1798}, mr={1411488}}
\end{barticle}
\bptok{imsref}%
\endbibitem

\bibitem{TTA05}
\begin{barticle}[mr]
\bauthor{\bsnm{Taylor},~\bfnm{Jonathan}\binits{J.}},
\bauthor{\bsnm{Takemura},~\bfnm{Akimichi}\binits{A.}} \AND
\bauthor{\bsnm{Adler},~\bfnm{Robert~J.}\binits{R.~J.}}
(\byear{2005}).
\btitle{Validity of the expected {E}uler characteristic heuristic}.
\bjournal{Ann. Probab.}
\bvolume{33}
\bpages{1362--1396}.
\bid{doi={10.1214/009117905000000099}, issn={0091-1798}, mr={2150192}}
\end{barticle}
\bptok{imsref}%
\endbibitem

\bibitem{TraubWW88}
\begin{bbook}[mr]
\bauthor{\bsnm{Traub},~\bfnm{J.~F.}\binits{J.~F.}},
\bauthor{\bsnm{Wasilkowski},~\bfnm{G.~W.}\binits{G.~W.}} \AND
\bauthor{\bsnm{Wo{\'z}niakowski},~\bfnm{H.}\binits{H.}}
(\byear{1988}).
\btitle{Information-based Complexity}.
\bpublisher{Academic Press},
\blocation{Boston, MA}.
\bid{mr={0958691}}
\end{bbook}
\bptok{imsref}%
\endbibitem

\bibitem{Wos96}
\begin{bincollection}[mr]
\bauthor{\bsnm{Wo{\'z}niakowski},~\bfnm{Henryk}\binits{H.}}
(\byear{1997}).
\btitle{Computational complexity of continuous problems}.
In \bbooktitle{Nonlinear Dynamics, Chaotic and Complex Systems ({Z}akopane, 1995)}
\bpages{283--295}.
\bpublisher{Cambridge Univ. Press},
\blocation{Cambridge}.
\bid{mr={1481726}}
\bptnote{check year}%
\end{bincollection}
\bptok{imsref}%
\endbibitem

\bibitem{Yor92}
\begin{barticle}[mr]
\bauthor{\bsnm{Yor},~\bfnm{Marc}\binits{M.}}
(\byear{1992}).
\btitle{On some exponential functionals of {B}rownian motion}.
\bjournal{Adv. in Appl. Probab.}
\bvolume{24}
\bpages{509--531}.
\bid{doi={10.2307/1427477}, issn={0001-8678}, mr={1174378}}
\end{barticle}
\bptok{imsref}%
\endbibitem

\bibitem{Zeger88}
\begin{barticle}[mr]
\bauthor{\bsnm{Zeger},~\bfnm{Scott~L.}\binits{S.~L.}}
(\byear{1988}).
\btitle{A regression model for time series of counts}.
\bjournal{Biometrika}
\bvolume{75}
\bpages{621--629}.
\bid{doi={10.1093/biomet/75.4.621}, issn={0006-3444}, mr={0995107}}
\end{barticle}
\bptok{imsref}%
\endbibitem

\end{thebibliography}
\end{document}